\documentclass[11pt,reqno]{amsart}

\usepackage{amsfonts,amsthm,amsmath,amssymb,amscd}
\usepackage{mathtools}
\usepackage{esint}
\usepackage{graphics}
\usepackage{indentfirst}
\usepackage{cite}
\usepackage{latexsym}
\usepackage[dvips]{epsfig}
\usepackage{mathrsfs}
\usepackage{hyperref}
\usepackage{color}
\date{}
\allowdisplaybreaks \allowdisplaybreaks[2]

\numberwithin{equation}{section}

\usepackage{graphics}
\usepackage{epsfig}
\newtheorem{theorem}{Theorem}[section]

\newtheorem{lemma}[theorem]{Lemma}
\newtheorem{proposition}[theorem]{Proposition}

\newtheorem{rmk}{Remark}
\numberwithin{rmk}{section}

\newenvironment{pf}{{\noindent \it \bf Proof:}}{{\hfill$\Box$}\\}
\providecommand{\abs}[1]{\left\vert#1\right\vert}
\providecommand{\norm}[1]{\left\Vert#1\right\Vert}

\def\bA{\boldsymbol{A}}
\def\boe{\boldsymbol{\eta}}
\def\bn{\boldsymbol{n}}
\def\bu{\boldsymbol{u}}
\def\bv{\boldsymbol{v}}
\def\bF{\boldsymbol{F}}
\def\D{\partial}
\def\DD{\nabla}
\def\Div{\text{div}}
\def\MH{\mathcal{H}}
\def\hal{\frac{1}{2}}
\def\ddt{\frac{d}{dt}}
\def\vep{\varepsilon}
\def\ls{\lesssim}

\begin{document}

\title[Inviscid limits of Compressible Viscoelastic surface waves]{Vanishing viscosity limits for the free boundary problem of compressible viscoelastic fluids with surface tension}

\author{Xumin~Gu}
\address[Xumin Gu]{School of Mathematics, Shanghai University of finance and economics, Shanghai 200433, China}
\email{\href{mailto:gu.xumin@shufe.edu.cn}{gu.xumin@shufe.edu.cn}}

\author{Yu~Mei}
\address[Yu Mei]{School of Mathematics and Statistics, Northwestern Polytechnical University, Xi'an 710129, China}
\email{\href{mailto:yu.mei@nwpu.edu.cn}{yu.mei@nwpu.edu.cn}}

\keywords{Free boundary; Viscoelastic fluids; Vanishing viscosity; Compressible fluids; Elastodynamics.}

\subjclass[2010]{35Q35, 35R35, 76A10, 76N10, 76N20}
\maketitle
\begin{abstract}
We consider the free boundary problem of compressible isentropic neo-Hookean viscoelastic fluid equations with surface tension. Under the physical kinetic and dynamic conditions proposed on the free boundary, we investigate regularities of classical solutions to viscoelastic fluid equations in Sobolev spaces which are uniform in viscosity and justify the corresponding vanishing viscosity limits. The key ingredient of our proof is that the deformation gradient tensor in Lagrangian coordinates can be represented as a parameter in terms of flow map  so that  the inherent structure of the elastic term improves the uniform regularity of normal derivatives in the limit of vanishing viscosity. This result indicates that the boundary layer does not appear in the free boundary problem of compressible viscoelastic fluids which is different to the case studied by the second author for the free boundary compressible Navier-Stokes system.

\end{abstract}


\section{Introduction}
\subsection{Formulation in Eulerian coordinates}
Consider the motion of  compressible isentropic neo-Hookean viscoelastic fluids in a time-dependent domain $\Omega^\vep(t)$, whose exterior is assumed to be the atmosphere of a given constant pressure $p_e>0$. The boundary of domain $\Gamma^\vep(t)$ is free to move and subject to physical kinetic and dynamic conditions under the effect of surface tension.  The governing equations of such a free boundary problem can be written as follows
\begin{equation}\label{FVEE}
	\begin{cases}
		\D_t \rho^\varepsilon+\Div(\rho^\varepsilon \bu^\varepsilon)=0,&\text{in } \Omega^\vep(t),\\
		\D_t(\rho^\varepsilon \bu^\varepsilon)+\Div(\rho^\varepsilon \bu^\varepsilon\otimes \bu^\varepsilon)+\DD p(\rho^\varepsilon)-\Div\boldsymbol{T}^\varepsilon=\Div(\rho^\vep \bF^\vep {\bF^\vep}^\mathsf{T}),&\text{in } \Omega^\vep(t),\\
		\D_t\bF^\vep+\bu^\vep\cdot \DD \bF^\vep=\DD \bu^\vep\bF^\vep,&\text{in } \Omega^\vep(t),
	\end{cases}
\end{equation}
where $\rho^\vep$ is the density, $\bu^\varepsilon=(u_1^\varepsilon,\cdots,u_d^\varepsilon)^\mathsf{T}\in\mathbb{R}^d$ is the velocity, $\bF^\varepsilon=(F^\varepsilon_{ij})\in\mathbb{M}^{d\times d}$ is the deformation gradient. Here $\mathsf{T}$ denotes the transpose of matrix. The pressure $p(\rho^\varepsilon)$ only depends on the density $\rho^\varepsilon$ in the isentropic case which is assumed here to satisfy
\begin{equation}\label{gamma-law}
	p(\rho^\vep)=A(\rho^\varepsilon)^\gamma,\quad\gamma>1.
\end{equation} 
For simplicity, we take $A=1$ here and afterwards.  The viscous stress tensor $\boldsymbol{T}^\vep$ is given by 
\begin{equation}\label{vis-tensor}
	\boldsymbol{T}^\vep=2\mu\varepsilon \boldsymbol{Su}^\varepsilon+\lambda\varepsilon \Div\bu^\varepsilon\boldsymbol{I},
\end{equation}
where $\boldsymbol{Su}^\varepsilon=(\DD\bu^\varepsilon+(\DD\bu^\vep)^\mathsf{T})/2$ is the symmetric part of $\DD\bu^\varepsilon$, and $\mu\vep$, $\lambda\vep$ are viscosity coefficients satisfying the physical constrain $\mu>0,2\mu+d\lambda>0$. To study this free boundary problem, we impose the following two boundary conditions on $\Gamma^\vep(t)$. On the one hand, the kinetic boundary condition, which states that the free
boundary moves along the fluid particles, reads
 \begin{equation}\label{VKBC}
	\boldsymbol{\mathcal{V}}^\vep(\Gamma^\vep(t))=\bu^\vep\cdot \bn^\vep,\text{ on }\Gamma^\vep(t),
\end{equation}
where $\boldsymbol{\mathcal{V}}^\vep(\Gamma^\vep(t))$ is the normal velocity of $\Gamma^\vep(t)$ and $\bn^\vep$ is the outward unit normal vector of $\Gamma(t)$. On the other hand, the dynamic boundary condition, when the surface tension is considered, can be written as
 \begin{equation}\label{VDBC}
	(-p(\rho^\vep)\boldsymbol{I}+\boldsymbol{T}^\vep+ (\rho^\vep\bF^\vep{\bF^\vep}^\mathsf{T}-\boldsymbol{I}))\bn^\vep+p_e\bn^\vep=\sigma H\bn^\vep, \text{ 
		on }\Gamma^\vep(t),
\end{equation}
which represents the balance of stress tensors on both sides of the free boundary. Here $H$ is the twice mean curvature of $\Gamma^\vep(t)$, $\sigma>0$ is the surface tension coefficient and $p_e>0$ is a given constant pressure outside $\Omega(t)$ which can guarantee the non-vacuum density of viscoelastic fluid if initially. We also impose the following initial data for \eqref{FVEE}
\begin{equation}\label{V-Ini-D}
	(\rho^\vep,\bu^\vep,\bF^\vep)(x,t=0)=(\rho_0^\vep,\bu_0^\vep,\bF_0^\vep), \quad x\in \Omega^\vep(0).
\end{equation}
When the fluid is inviscid, the free boundary problem \eqref{FVEE}-\eqref{V-Ini-D} can reduce to the one of compressible isentropic inviscid elastic fluid, through formally taking $\vep\rightarrow 0$, which reads as
\begin{equation}\label{FEE}
\begin{cases}
\D_t \rho+\Div(\rho \bu)=0,&\text{in } \Omega(t),\\
\D_t(\rho \bu)+\Div(\rho \bu\otimes \bu)+\DD p(\rho)=\Div(\rho \bF\bF^\mathsf{T}), &\text{in }\Omega(t),\\
\D_t\bF+\bu\cdot \DD \bF=\DD \bu\bF, &\text{in }\Omega(t),\\
\boldsymbol{\mathcal{V}}(\Gamma(t))=\bu\cdot \bn,&\text{on }\Gamma(t),\\
 (-p\boldsymbol{I}+ (\rho\bF\bF^\mathsf{T}-\boldsymbol{I}))\bn+p_e\bn=\sigma H\bn, &\text{on }\Gamma(t),\\
 (\rho,\bu,\bF)(x,t=0)=(\rho_0,\bu_0,\bF_0),&\text{in }\Omega(0).
\end{cases}
\end{equation}
In this paper, we will investigate the vanishing viscosity limits of the free bounday problem of compressible isentropic viscoelastic fluid equations \eqref{FVEE}-\eqref{V-Ini-D}, that is, whether solutions of \eqref{FVEE}-\eqref{V-Ini-D} converge to ones of the free boundary elastodynamic equations \eqref{FEE} modeling inviscid compressible isentropic flows as viscosity tends to zero.
\subsection{Motivations and Related Results}

The vanishing viscosity limit for viscous fluid is one of the most fundamental problems in the mathematical
theory of fluid mechanics. It has been extensively studied by many mathematicians in various settings of domains and physical boundary conditions as well as different kinds of fluids. For the incompressible homogeneous Newtonian fluids, the justification of vanishing viscosity limits in the whole space have been proved in  \cite{mcgrath1968,swann1971,kato1972} for smooth solutions and in \cite{constantin1995,gallay2011,constantina, ciampa2021}  for irregular ones. However, in the presence of  physical boundaries, such a problem becomes much more complicated due to the possible appearance of boundary layers. When the Dirichlet boundary condition is imposed, the boundary layer, as illustrated by the Prandtl theory, must appear and be very strong so that the verification of vanishing viscosity limit for incompressible Newtonian fluids becomes one of the most challenging open problems in mathematical fluid mechanics. Some important progress on, but still far from completely solving, this issue are made in \cite{sammartino1998,sammartino1998a,maekawa2014,wang2017,fei2018,nguyen2018} for the analytic initial data at least near the boundary. Whereas, when the Navier-slip boundary condition is imposed, the boundary layer becomes weaker so that the vanishing viscosity limit holds even if vorticity is produced at the boundary. We can refer to \cite{clopeau1998,filho2005,iftimie2006,xiao2007,iftimie2011,masmoudi2012a} and references therein for considerable progress. For the compressible Newtonian fluids, the vanishing viscosity limit problems are discussed in \cite{xin1999,wang2005,wang2012a,sueur2014,wang2015,paddick2016,wang2016a}. Precisely, Xin-Yanagisawa \cite{xin1999} studied the vanishing viscosity limit of the linearized compressible Navier-Stokes system  with the no-slip boundary condition in the 2-D half plane. Wang-Williams \cite{wang2012a} constructed a boundary layer solution of the compressible Navier-Stokes equations with Navier-slip boundary conditions in 2-D half plane. Wang-Xin-Yong \cite{wang2015} obtained an uniform regularity for the solutions of the compressible Navier-Stokes with general Navier-slip boundary conditions in 3-D domains with curvature, especially, the vanishing viscosity limit of viscous solution to the corresponding inviscid one was also obtained with  rate of convergence in  $L^\infty$. It is also shown that the boundary layer for  density is weaker than the one for velocity fields. Very recently, Wang-Xie \cite{wang2021a} justified the vanishing viscosity limit of solutions for the compressible viscoelastic flows under the no-slip boundary condition governed by the viscoelastic equations, based on the uniform conormal regularity estimates.

Although much literature exists for the vanishing viscosity limit problems on a fixed domain, less discussions are available for viscous surface waves.  To our knowledge, the first result of the vanishing viscosity limit for viscous surface waves due to Masmoudi and Rousset \cite{masmoudi2017}, who proved the local existence of solutions to the free boundary incompressible Navier-Stokes system with uniform in viscosity regularity in conormal Sobolev and Lipschitz spaces. Later, the vanishing viscosity and surface tension limits were established in \cite{elgindi2018,wang2021} for the incompessible viscous surface waves with surface tension. For compressible surface waves, the second author, Wang and Xin in \cite{mei2018} established uniform regularities of both density and velocity in conormal Sobolev and Lipschitz spaces and justified the vanishing viscosity and surface tension limits. 

In the results mentioned above for the vanishing viscosity limits of viscous surface waves, we can only expect uniform Lipschitz bounds, but not higher Sobolev norms $H^{k}$ ($k\geq2$), since a boundary layer generally appears near the free boundary for viscous surface waves and uniform bounds of many normal derivatives break down as viscosity tends to zero. This motivate us to study whether there exists a mechanism to prevent the appearance of boundary layers for viscous surface waves. Recently, the first author and Lei in \cite{gu2020a} proved the local well-posedness of free-bounary incompressible elastodynamics with surface tension by establishing uniform in viscosity Sobolev regularity of a modified viscoelastic surface waves. This result indicates that the inherent structure of the elastic term on the boundary may prevent the appearance of  boundary layers and enable us to justify the vanishing viscosity limit in standard Sobolev spaces. Our aim in this paper is further to discuss this smoothing effect of elasticity for compressible isentropic viscoelastic in the limit of vanishing viscosity. Precisely, we will strictly prove  that solutions of free boundary problem \eqref{FVEE}-\eqref{V-Ini-D} converges to ones of \eqref{FEE} in the $H^3$-sense. 

Let us also review here the well-posedness results of viscoelastic and elastic fluids. For the Cauchy problem, we refer to \cite{lin2005,chen2006,lin2008,lei2008} for the local and global well-posedness of  the incompressible viscoelastic fluid system. Sideris and Thomases\cite{sideris2005,sideris2007} established the global well-posedness of the three-dimensional incompressible neo-Hookean elastodynamic system with small initial data. Taking advantage of the strong null structure in Lagrangian coordinate, Lei \cite{lei2016a} obtained the two-dimensional result by using Klainerman vector field and Alinhac's ghost weight methods.  We also refer to Cai et al.\cite{cai2019} for the vanishing viscosity limit for global in time solutions to incompressible viscoelasticity in $\mathbb{R}^2$. For free boundary problems, the local  and global well-posedness of incompressible viscoelastic flows in an infinite strip with surface tension are respectively established in \cite{lemeur2011} for large initial data and \cite{xu2013} for small ones. Recently, Di Iorio, Marcati and Spirito \cite{diiorio2017,diiorio2020} proved the existence of splash singularities for incompressible viscoelastic flows. The well-posedness theory of the free boundary problem for incompressible neo-Hookean elastodynamics is more complicated. A prior estimates and local well-posedness are obtained in \cite{hao2016,hu2019,li2019a,gu2020} under some special boundary conditions, such as $\bF^\mathsf{T}\bn=0, p=0$ or $p=0,(\bF{\bF}^\mathsf{T}-\boldsymbol{I})\bn=({\bF}^\mathsf{T}-\boldsymbol{I})\bn=0$. Recently, the first author and Lei in \cite{gu2020a} established the local well-posedness result for free boundary problems being subject to the natural force balance law, which is consistent with  $\eqref{FEE}_5$ in the compressible case here. Regarding to the compressible neo-Hookean elastodynamics,  under the special boundary condition $\bF^\mathsf{T}\bn=0, p=0$, Trakhinin \cite{trakhinin2018} proved the local in time existence of a unique smooth solution of the free boundary problem by the Nash-Moser iteration approach, if the non-collinearity of $\bF$ or the Rayleigh–Taylor sign condition is satisfied. Very recently, Zhang \cite{zhang2021a} obtained the local well-posedness in Sobolev spaces by the combination of classical energy method and hyperbolic approach and also established the incompressible limit.

\subsection{Reformulation in Lagrangian coordinates}
To study the free boundary problem \eqref{FVEE}-\eqref{V-Ini-D}, we use the Lagrangian flow map to transform it  into the corresponding problem in a fixed reference domain $\Omega$. Let $\boe^\vep(x,t)\in\Omega^\varepsilon(t)$ be the ``position" of compressible viscoelastic fluid particle $x$ at time $t$, i.e.
\begin{equation}\label{flow map}
\begin{cases}
	\partial_t\boe^\vep(x,t)=\bu(\boe^\vep(x,t),t)& \text{ for }t>0, x\in \Omega\\
	\boe^\vep(x,0)=\boe_0^\vep(x), &\text{ for }x\in \Omega,
\end{cases}
\end{equation}
where $\boe^\vep_0$ is a diffeomorphism from the reference domain $\Omega$ to the initial domain $\Omega^\varepsilon(0)$ satisfying
\begin{equation}\label{Ini-flow map}
\bF^\vep_0(\boe^\vep_0)=\nabla\boe^\varepsilon_0.
\end{equation}
We set the following Lagrangian variables 
\begin{align*}
f^\vep(x,t)=\rho^\varepsilon(\boe^\vep(x,t),t),~~\bv^\vep(x,t)=\bu^\vep(\boe^\vep(x,t),t), ~~\boldsymbol{G}^\vep(x,t)=\bF^\varepsilon(\boe^\vep(x,t),t)
\end{align*}
and introduce
\begin{align*}
\boldsymbol{A}^\vep=(\DD\boe^\vep)^{-\mathsf{T}},~~J^\vep=\text{det} \DD\boe^\vep,~~ \boldsymbol{a}^\vep=J\boldsymbol{A}^\vep,~~q(f^\vep)=p(f^\vep)+1-p_e.
\end{align*}
 Then, utilizing the chain rule and Einstein's summation convention for repeated indexes, the free boundary compressible viscoelastic fluid equations \eqref{FVEE}-\eqref{V-Ini-D} can be written by components of Lagrangian variables in $\Omega$ as follows
\begin{equation}\label{FEL}
\begin{cases}
\D_t\eta^\vep_i=v^\varepsilon_i,&\text{in }\Omega,\\
\D_tf^\vep+f^\vep A^\vep_{kj}\D_jv^\vep_k=0,&\text{in }\Omega,\\
f^\vep\D_tv^\vep_i+A^\vep_{ik}\D_kq-2\mu\vep A^\vep_{kl}\D_l(S_Av)^\vep_{ik}-\lambda\vep A^\vep_{ij}\D_j(\Div_A v)^\vep-A^\vep_{kl}\D_l(f^\vep G^\vep_{ij}G^\vep_{kj})=0,&\text{in }\Omega,\\
\D_tG^\vep_{ij}=A^\vep_{kl}\D_lv^\vep_iG^\vep_{kj},&\text{in }\Omega,\\
-qa^\vep_{ij}N_j+f^\vep G^\vep_{ij}{G^\vep}^\mathsf{T}_{jk}a^\vep_{kl}N_l+2\mu\vep S_A(v)^\vep_{ik}a^\vep_{kl}N_l+\lambda\vep \Div_A (v)^\vep a^\vep_{ij} N_j=\sigma \sqrt{g}\Delta_g(\eta^\vep_i)&\text{on }\Gamma,\\
(f^\vep,\bv^\vep,\boldsymbol{G}^\vep,\boe^\vep)(x)|_{t=0}=(\rho^\vep_0(\boe^\vep_0(x)), \bu^\vep_0(\boe^\vep_0(x)), \bF^\vep_0(\boe^\vep_0(x)),\boe^\vep_0(x))&\text{in }\Omega,
\end{cases}
\end{equation}
where $(S_Av)^\vep_{ik}=(A^\varepsilon_{kj}\D_jv^\vep_i+A_{ij}^\vep\D_jv^\vep_k)/2$, $(\Div_A v)^\vep=A^\vep_{kl}\D_lv^\vep_k$, $\boldsymbol{N}=(N_1,\cdots,N_d)^\mathsf{T}$ is the unit outward normal vector to $\Gamma:=\D \Omega$, and $\Delta_g$ is the Laplacian-Beltrami operator on the curve $\boe^\vep(t,\Gamma)$ given by 
\begin{equation*}
\Delta_g(\eta^\vep_i)=\frac{1}{\sqrt{ g}}\D_\alpha(\sqrt{ g}g^{\alpha\beta}\D_\beta\eta^\vep_i),\quad g^{\alpha\beta}=\D_\alpha\eta^\vep_j\D_\beta\eta^\vep_j,\quad g=\text{det} (g^{\alpha\beta}).
\end{equation*}
Without loss of generality, we consider the two-dimensional problem (i.e. $d=2$) and the reference domain is given by $$\Omega=\mathbb{T}\times(0,1),$$
where $\mathbb{T}$ denotes the $1$-torus. The boundary of $\Omega$ is then given by the horizontally flat bottom and top
$$\Gamma=\mathbb{T}\times(\{0\}\cup\{1\}).$$ 
In such a case, $\sqrt{g}\Delta_g(\eta^\vep_i)=\D_1\left(\frac{\D_1\eta_i^\vep}{|\D_1\boe^\vep|}\right)$ and
\begin{equation*}
	\boldsymbol{N}=\begin{cases}
		\boldsymbol{e}_2, &x_2=1\\
		-\boldsymbol{e}_2,&x_2=0.
	\end{cases}
\end{equation*}
Furthermore, from Jacobi's formula
\begin{equation}\label{Jacobi's formula}
J^\vep_t=J^\vep A^\vep_{kj}\D_jv^\vep_k=a^\vep_{kj}\D_jv^\vep_k,
\end{equation}
 one can solve the first equation in \eqref{FEL} to obtain $\D_t(f^\vep J^\vep)=0$, which implies
\begin{equation}\label{f-eqn}
f^\vep=\tilde{\rho}^\vep_0(J^\vep)^{-1},\text{ with }\tilde{\rho}^\vep_0:=\rho^\vep_0(\boe^\vep_0)J^\vep_0.
\end{equation}
Next, since $\D_t A^\vep_{ij}=-A^\vep_{ik}\D_kv^\vep_lA^\vep_{lj}$, it is not difficult to derive from the fourth equation of \eqref{FEL} that $\D_t({\bA^\vep}^\mathsf{T}\boldsymbol{G}^\vep)=0,$ 
which yields
\begin{equation}\label{G-eqn}
\boldsymbol{G}^\vep=\DD \boe^\vep(\DD\boe^\vep_0)^{-1}\bF^\vep_0(\boe_0)=\DD\boe^\vep
\end{equation}
Note that $\boldsymbol{a}^\vep$ is the cofactor of  $\nabla\boe^\vep$, we have the following Piola identity 
\begin{equation}\label{Piola}
\D_l a^\vep_{kl}=\D_l(J^\vep A^\vep_{kl})=0
\end{equation}
which implies that
\begin{equation*}
-A^\vep_{kl}\D_l(f^\vep G^\vep_{ij}G^\vep_{kj})=-(J^\vep)^{-1}\D_l(A^\vep_{kl}f^\vep J^\vep G^\vep_{ij}G^\vep_{kj})
\end{equation*}
Therefore, substituting \eqref{f-eqn} and \eqref {G-eqn} into remainder equations in \eqref{FEL}, we obtain
\begin{equation}\label{FEL-final}
\begin{cases}
\D_t\eta^\vep_i=v^\vep_i,&\text{in }\Omega\\
\tilde{\rho}_0^\vep\D_tv^\vep_i+a_{ik}^\vep\D_kq-2\mu\vep a^\vep_{kl}\D_l(S_Av)^\vep_{ik}-\lambda\vep a^\vep_{ij}\D_j(\Div_A v)^\vep-\D_j(\tilde{\rho}_0^\vep\D_j\eta^\vep_i)=0,&\text{in }\Omega\\
-qa_{i2}^\vep+2\mu\vep (S_Av)^\vep_{ik}a^\vep_{k2}+\lambda\vep (\Div_Av)^\vep a^\vep_{i2}+\tilde{\rho}^\vep_0\D_2\eta^\vep_i=\sigma \D_1\left(\frac{\D_1\eta^\vep_i}{|\D_1\boe^\vep|}\right),&\text{on }\Gamma\\
(\bv^\vep,\boe^\vep)(x)|_{t=0}=(\bu_0^\vep(\boe^\vep_0(x)),\boe^\vep_0(x)),&\text{in }\Omega.
\end{cases}
\end{equation}
If formally taking $\vep=0$, then we have the following free boundary elastodynamic equations in Lagrangian coordinates
\begin{equation}\label{FEL-final-2}
	\begin{cases}
		\D_t\eta_i=v_i,&\text{ in }\Omega\\
		\tilde{\rho}_0\D_tv_i+a_{ik}\D_kq-\D_j(\tilde{\rho}_0\D_j\eta_i)=0,&\text{ in }\Omega\\
		-qa_{i2}+\tilde{\rho}_0\D_2\eta_i=\sigma \D_1\left(\frac{\D_1\eta_i}{|\D_1\eta|}\right),&\text{ on }\Gamma\\
		(\bv,\boe)(x)|_{t=0}=(\bu_0(\boe_0(x)),\boe_0(x))&\text{ in }\Omega,
	\end{cases}
\end{equation}
From now on, we will focus on the free boundary viscoelastic and elastic fluid systems \eqref{FEL-final}, \eqref{FEL-final-2} in the Lagrangain coordinates.
\subsection{Notation Conventions and Function spaces}
Before stating our main results, we make the following notation conventions. We use Einstein's summation convention for repeated indexes throughout the paper. The bold font, such as $\bu,\bv,\bF,\boldsymbol{G}, \bA$ ect., are used to denote a vector or matrix, while the non-bold font, $u_i,v_i,F_{ij},G_{ij}, A_{ij}$ ect., stand for corresponding each elements. We use $\bar{\D}:=\D_1,\D_t$ to denote the tangential derivatives and $\D$ also including $\D_2$. The usual $L^p$ spaces, Sobolev spaces $W^{m,p}$ and $H^m=W^{m,2}$ on both the domain $\Omega$ and the boundary $\Gamma$ are used. For notational simplifications, the norms of these spaces defined on $\Omega$ are denoted by $\|{\cdot}\|_{L^p}, \|{\cdot}\|_{W^{m,p}}$ and $\|{\cdot}\|_{m}$, and the norms of these spaces defined on $\Gamma$ are denoted by $|{\cdot}|_{L^p}, |{\cdot}|_{W^{m,p}}$ and $|{\cdot}|_{m}$.  For real $s\geq 0$, the Hilbert space $H^s(\Gamma)$ and the boundary norm $|{\cdot}|_s$ (or $|{\cdot}|_{H^s(\Gamma)}$) is defined by interpolation. The negative-order Sobolev space $H^{-s}(\Gamma)$ are defined via duality: for real $s\geq 0, H^{-s}(\Gamma):=[H^{s}(\Gamma)]^\prime.$ 
Moreover, $\|u\|_{H_{tan}^m}^2:=\sum\limits_{\ell\leq m}\|\D_1^\ell u\|_0^2$ is introduced for a function $u\in L^2$.
We also introduce the spatial-temporal Sobolev norms of any function $u\in L^2$ on $\Omega$ at the instantaneous time $t$ as:
\begin{equation*}
	\|u(t)\|_{\MH^m}^2:=\sum_{\ell\leq m}\|\partial_t^{\ell}u\|_{m-\ell}^2,\quad\|u\|_{\MH^{m,k}_{tan}}:=\sum_{\ell\leq m-k}\|\D_t^\ell\D_1^{m-\ell}u\|_0^2,\quad \|u\|_{\MH^{m}_{tan}}:=\|u\|_{\MH^{m,0}_{tan}},
\end{equation*}
and the spatial-temporal Sobolev norms on $\Gamma$ as:
$$|u(t)|_{\MH^m}^2:=\sum_{\ell\leq m}|\partial_t^{\ell}u|_{m-\ell}^2.$$
$\|\cdot\|_{L^p_t(\textbf{X})}$ is used as the norm of the space $L^p([0, t];\textbf{X})$. $C$ is used to denote generic constants, which only depends on the domain $\Omega$ and the boundary $\Gamma$, and $f\lesssim g$ is used to denote $f\leq Cg$.  $P(\cdot)$ denotes a generic polynomial function of its arguments, and the polynomial coefficients are generic constants $C$.

\subsection{Main Theorems}
The aim of this paper is to get a local well-posedness result for classical solutions to \eqref{FEL-final}  in a temporal interval independent of the viscosity $\vep\in(0,1]$. This result will also imply the local existence of strong solutions to the free boundary elastodynamic equations \eqref{FEL-final-2} with surface tension. To obtain the uniform regularities, we define the energy functional for the free boundary viscoelastic fluid equations \eqref{FEL-final} as
\begin{align}
	\mathfrak{E}^\vep(t)=&\|\boe^\vep\|_{\MH^m}^2(t)+\|\DD\boe^\vep\|_{\MH^{m,1}_{tan}}^2(t)+|\bar{\D}^{m-1}\D_1^2\boe^\vep\cdot\bn^\vep|_0^2(t)+\vep\|\DD^2\boe^\vep\|_{\MH^{m-1}}^2(t)\nonumber\\
	&+\int_{0}^{t}\|\DD\boe^\vep\|_{\MH^m}^2+\|\D_t^m\bv^\vep\|_0^4+\|\D_t^m\DD \boe^\vep\|_0^4+|\D_1\D_t^m\boe^\vep\cdot\bn^\vep|_0^4\,d\tau\\
	&+\int_{0}^{t}\vep^2\|\DD\bv^\vep\|_{\MH^m}^2+\vep\|\DD \bv^\vep\|_{\MH^{m,1}_{tan}}^2\,d\tau+\vep^2\left(\int_{0}^{t}\|\D_t^m\DD\bv^\vep\|_0^2\,d\tau\right)^2.\nonumber
\end{align}
The corresponding energy for the free boundary elastodynamics equations \eqref{FEL-final-2} are
\begin{align}
	\mathfrak{E}(t)=&\|\boe\|_{\MH^m}^2(t)+\|\DD\boe\|_{\MH^{m,1}_{tan}}^2(t)+|\bar{\D}^{m-1}\D_1^2\boe\cdot\bn|_0^2(t)\nonumber\\
	&+\int_{0}^{t}\|\DD\boe\|_{\MH^m}^2+\|\D_t^m\bv\|_0^4+\|\D_t^m\DD \boe\|_0^4+|\D_1\D_t^m\boe\cdot\bn|_0^4\,d\tau.
\end{align}
We also require that the initial data satisfy the compatibility condition on the boundary. That is, the high order temporal initial data $(\D_t^\ell \bv^\vep(0),\D_t^\ell\boe(0), \D_t^\ell q(0))$ can be defined by
\begin{align*}
	&\D_t^\ell v_i^\vep(0)=\tilde{\rho}_0^{-1}\D_t^{\ell-1}(-a_{ik}^\vep\D_kq+2\mu\vep a^\vep_{kl}\D_lS_A(v)^\vep_{ik}+\lambda\vep a^\vep_{ij}\D_j\Div_A (v)^\vep+\D_j(\tilde{\rho}_0^\vep\D_j\eta^\vep_i))\big|_{t=0}\\
	&\D_t^\ell\eta^\vep_i(0)=\D_t^{\ell-1}v_i^\vep(0),\quad \D_t^\ell q^\vep(0)=\tilde{\rho}_0^\gamma\D_t^\ell(J^\vep)^{-\gamma}\big|_{t=0}
\end{align*}
inductively for $\ell=1,\cdots, m$. And these data should satisfy 
\begin{equation}\label{compat-cond}
        \D_t^\ell\left(-qa_{i2}^\vep+2\mu\vep (S_Av)^\vep_{ik}a^\vep_{k2}+\lambda\vep (\Div_Av)^\vep a^\vep_{i2}+\tilde{\rho}^\vep_0\D_2\eta^\vep_i\right)\Bigg|_{t=0}=\sigma \D_t^\ell\left(\D_1\left(\frac{\D_1\eta^\vep_i}{|\D_1\boe^\vep|}\right)\right)\Bigg|_{t=0}\text{on }\Gamma
\end{equation}
for $\ell=1,\cdots,m-1$.

The uniform regularity theorem states as follows.
\begin{theorem}\label{Theorem 1}
	Let $m\geq 4$. Suppose that the initial data $(\rho_0^\vep, \boe_0^\vep, \bv_0^\vep)$ satisfy the compatibility conditions \eqref{compat-cond}, and have the following uniform bounds
	\begin{align}
		&c_0\leq \tilde{\rho}_0^\vep\leq C_0,\label{rho-uni-bdd}\\ \|\boe_0^\vep\|_{\MH^m}^2+\|\DD\boe^\vep_0\|_{\MH^{m,1}_{tan}}^2&+\|\D_t^m\bv^\vep_0\|_0^2+\|\D_t^m\DD\boe^\vep_0\|_0^2\nonumber\\
		+|\bar{\D}^{m-1}\D_1^2\boe^\vep_0\cdot\bn^\vep_0|_0&+|\bar{\D}^{m}\D_1\boe^\vep_0\cdot\bn^\vep_0|_0+\vep\|\DD^2\boe^\vep_0\|_{\MH^{m-1}}^2\leq C_0,\label{v-uni-bdd}
	\end{align}
for some generic constants $c_0,C_0$. Then there exists a $T_0>0$, independent of $\vep$, and a unique solution $(\boe^\vep,\bv^\vep)$ to the free boundary problem \eqref{FEL-final} on the time interval $[0,T_0]$ such that 
	\begin{equation}
		\sup_{t\in[0,T_0]}\mathfrak{E}^\vep(t)\leq C_1,
	\end{equation}
	where $C_1$ is a generic constant depending only on $c_0,C_0$.
\end{theorem}
\begin{rmk}
	The regularities of solutions implies that the flow map $\boe$ is at least Lipschitzian so that we can get the corresponding classical solutions in Eulerian coordinates.
\end{rmk}
\begin{rmk}
	The inherent structure of the elastic term enable us to perform uniform regularities in standard Sobolev spaces, which indicates that the boundary layer does not appear in the free boundary problem of compressible viscoelastic fluids. This is different to the case studied by the second author for the free boundary compressible Navier-Stokes system.
\end{rmk}
Based on the uniform regularities in Theorem \ref{Theorem 1}, we can justify the vanishing viscosity of limit for the free boundary problem of compressible viscoelastic fluid and obtain the local existence of classical solutions to the free boundary elastodynamic equations.
\begin{theorem}\label{vanishing viscosity}
	Under the assumptions of Theorem \ref{Theorem 1}, if we further assume  that there exist $(\rho_0,\boe_0,\bv_0)$ such that
	\begin{equation}
		\lim\limits_{\vep\rightarrow 0}\|\tilde{\rho}^\vep_0-\rho_0\|_0+\|\boe_0^\vep-\boe_0\|_0+\|\bv_0^\vep-\bv_0\|_0=0.
	\end{equation}
Then, there exists $(\boe,\bv)(t,\cdot)$ on the time interval $[0,T_0]$ such that
\begin{equation}
	\sup_{t\in[0,T_0]}\mathfrak{E}(t)\leq C_1
\end{equation}
and
\begin{equation}
	\lim\limits_{\vep\rightarrow 0}\sup_{t\in[0,T_0]}(\|\boe^\vep(t)-\boe(t)\|_{\MH^{m-1}}+\|\bv^\vep(t)-\bv(t)\|_{\MH^{m-2}})=0.
\end{equation}
Moreover, $(\boe,\bv)$ is the unique classical solution to the free boundary elastodynamic equations \eqref{FEL-final-2}.
\end{theorem}
\subsection{Sketch of proof}
We give some comments on the difficulties and ideas in the proof of main theorems. The crucial step is to derive uniform \textit{a priori} estimates of local classical solutions to \eqref{FEL-final} in a small time interval independent of $\vep$. Although the idea in \cite{gu2020a} can be adopted here to deal with elastic stress in some sense, we have to propose some new ideas to handle the compressibility of fluid.  

First, the estimate of $J$ now depends on the regularity of $\boe$, and the pressure $q$ in compressible fluid is  no longer a Lagrangian multiplier, whose estimates follows from a elliptic equation with Dirichlet or Neumann boundary condition for incompressible fluid. We derive the $L^2_T(\MH^m)$ estimates of $J$ and $q$ from a geometric identity, the Gagliardo-Nirenberg-M\"{o}ser-type inequality and the equation $q=\tilde{\rho}_0^{-\gamma}J^{-\gamma}$, which further gives the $L^\infty_T(\MH^{m-1})$ ones by using the fundamental theorem of calculus. (c.f. Lemma \ref{est-J-q}).
 
Next, for the tangential derivative estimates, we use the $L^2$-type energy estimate to get
\begin{align*}
	&\hal\ddt\int_\Omega(\tilde{\rho}_0|\bar{\D}^\alpha\bv|^2+\tilde{\rho}_0|\bar{\D}^\alpha\DD\boe|^2)dx-\int_\Omega\bar{\D}^\alpha q\bar{\D}^\alpha(a_{ik}\D_kv_i)dx+\text{dissipative terms}\\
	&+\int_\Gamma\left(\bar{\D}^\alpha qa_{i2}-\bar{\D}^\alpha\left(2\mu\vep a_{k2} (S_Av)_{ik}+\lambda\vep a_{i2}(\Div_Av)+\tilde{\rho}_0\D_2\eta_i\right)\right)\bar{\D}^\alpha v_i=\cdots
\end{align*}
Substituting the transport equation satisfied by pressure in compressible fluid, we can obtain 
\begin{align*}
	-\int_\Omega\bar{\D}^\alpha q\bar{\D}^\alpha(a_{ik}\D_kv_i)dx=\frac{1}{2\gamma }\frac{d}{dt}\int_\Omega\tilde{\rho}_0^{-\gamma}J^{\gamma+1}|\bar{\D}^\alpha q|^2dx+\cdots
\end{align*}
which provides the $L^2$-type energy of pressure. For the boundary integral in the tangential derivative estimates, we plug in $\eqref{FEL-final}_3$ and adopt the similar idea as \cite{gu2020a} to have
\begin{align*}
	&\int_\Gamma\left(\bar{\D}^\alpha qa_{i2}-\bar{\D}^\alpha\left(2\mu\vep a_{k2} (S_Av)_{ik}+\lambda\vep a_{i2}(\Div_Av)+\tilde{\rho}_0\D_2\eta_i\right)\right)\bar{\D}^\alpha v_i\\
	&=-\sigma\int_\Gamma\bar{\D}^\alpha\left(\frac{\D_1^2\eta_ka_{k2}}{|\D_1\boe|^3}\right)a_{i2}\bar{\D}^\alpha v_i-\int_\Gamma\mathfrak{B}\bar{\D}^\alpha a_{i2}\bar{\D}^\alpha v_i+\cdots,
\end{align*}
where $\mathfrak{B}=q+\D_1\eta_ka_{k2}/|\D_1\boe|^3$. Then, the first term on the right hand side of the above equation, which is related to the surface tension, provides the regularity of the boundary i.e.
\begin{align*}
	-\sigma\int_\Gamma\bar{\D}^\alpha\left(\frac{\D_1^2\eta_ka_{k2}}{|\D_1\boe|^3}\right)a_{i2}\bar{\D}^\alpha v_i=\frac{\sigma}{2}\ddt\int_\Gamma\frac{|\bar{\D}^\alpha\D_1\boe\cdot\bn|}{|\D_1\boe|}+\cdots
\end{align*}
While, by using the fact $a_{\cdot 2}=(-\D_1\eta_2,\D_1\eta_1)^\mathsf{T}=\D_1\boe^\perp$ and the anti-symmetric property of $\boe^\perp$ and $\boe$, the second term can reduce to
\begin{align*}
-\int_\Gamma\mathfrak{B}\bar{\D}^\alpha a_{i2}\bar{\D}^\alpha v_i=-\frac{1}{2}\frac{d}{dt}\int_\Gamma\mathfrak{B}\bar{\D}^\alpha a_{i2}\bar{\D}^\alpha\eta_i+\frac{1}{2}\int_\Gamma\D_1\mathfrak{B}(\bar{\D}^\alpha v_2\bar{\D}^\alpha \eta_1-\bar{\D}^\alpha v_1\bar{\D}^\alpha \eta_2)+\cdots	
\end{align*}
For the second term on the right hand side above in the case $\bar{\D}^\alpha=\D_t^m$, by using the decomposition \eqref{decomp-2}, we can write it into a troublesome term of form $\int_\Gamma fa_{i2}\D_t^mv_i$. Whereas, we can not use the duality argument as in \cite{gu2020a} to control it, because  the term $\|a_{ij}\D_j\D_t^mv_i\|_0$,  which follows from the normal trace estimate $|\D_t^mv_ia_{i2}|_{-\hal}\ls\|\DD\boe\|_{L^\infty}(\|\D_t^m\bv\|_0+\|a_{ij}\D_j\D_t^mv_i\|_0)$, can not be bounded for compressible fluids. Instead, we use integration by parts to get 
	\begin{align*}
	\int_{0}^{t}\int_\Gamma fa_{i2}\D_t^mv_idx_1d\tau&=\int_{0}^{t}\int_\Omega fa_{ij}\D_j\D_t^mv_idxd\tau+\int_{0}^{t}\D_jfa_{ij}\D_t^mv_idxd\tau\\
	&=\int_\Omega fa_{ij}\D_j\D_t^mv_idxd\tau\bigg|_0^t+\cdots
\end{align*} 
and bound it by $M_0+\delta\|\DD\D_t^m\boe(t)\|_0^2+P\left(\sup\limits_{t\in[0,T]}\mathfrak{E}(t)\right)$.

Finally, for the normal derivative estimates,  the approach used in \cite{zhang2021a}, based on the Hodge-type estimate and structure of wave equations of enthalpy, will be invalid for the physical dynamic boundary conditions \eqref{VDBC} studied here, because the boundary conditions for the equations of vorticity and enthalpy become much worse due to the elastic term $\rho \bF\bF^\mathsf{T}$ on the boundary. Moreover, we also can not directly use the same argument as in \cite{gu2020a},  in which the elastic term provides a Laplacian equation of $\boe$ with source terms involving only tangential derivatives. In our case, the elastic term provides a semi-linear elliptic equation of $\boe$ which further reduce to the following equation
\begin{align*}
	-\mathcal{A}_{ij}\D_2^2\eta_j-\mu\vep a_{k2}a_{k2}\D_2^2v_i-(\mu+\lambda)\vep a_{i2}a_{j2}\D_2^2v_j=\text{l.o.t},
\end{align*}
where $	\mathcal{A}_{ij}=\tilde{\rho}_0J\delta_{ij}+\gamma (\tilde{\rho}_0J^{-1})^\gamma a_{i2}a_{j2}$ is positively symmetric. Then, we can perform $L^2$-type energy estimates and control the normal derivatives inductively.
\section{Preliminary}
In this section, we recall some inequalities and derive some identities and elementary estimates 

\subsection{General inequalities}
The following Gagliardo-Nirenberg-M\"oser-type inequality will be used repeatedly in the paper.
\begin{lemma}
	The following product and commutator estimates hold:
	\begin{itemize}
	\item[(i)]For any  $m\in \mathbb{N}$, $|\alpha|+|\beta|=m$, and $g,h\in L^\infty(\Omega\times[0,t])\cap L^2([0,t];\MH^m)$, it holds that
	\begin{align}\label{pro-est}
	\|\D^\alpha g\D^\beta h\|_{L^2_t(L^2(\Omega))}&\lesssim\|g\|_{L^\infty_t(L^\infty)}\|h\|_{L^2_t(\MH^m)}+\|h\|_{L^\infty_t(L^\infty)}\|g\|_{L^2_t(\MH^m)}
	\end{align}
	\item[(ii)]For any $m\in \mathbb{N}$, $1\leq|\alpha|\leq m$, $h\in L^\infty(\Omega\times[0,t])\cap L^2([0,t];\MH^{m-1})$, and $g\in L^2([0,t];\MH^m)$ such that $\D g\in L^\infty(\Omega\times[0,t])$, it holds that
	\begin{align}
	\label{co1}
	\norm{\left[\D^{\alpha}, g\right]h}_{L^2_t(L^2)}&\lesssim\|h\|_{L^\infty_t(L^\infty)}\|\D g\|_{L^2_t(\MH^{m-1})}+\|\D g\|_{L^\infty_t(L^\infty)}\|h\|_{L^2_t(\MH^{m-1})},
	\end{align}
where $\left[\D^{\alpha}, g\right]h:= \D^{\alpha}(gh)-g\D^{\alpha} h.$ 
	\item[(iii)]For any $m\in\mathbb{N}$, $2\leq\abs{\alpha}\leq 3$, and $g, h\in L^\infty(\Omega\times[0,t])\cap L^2([0,t];\MH^{m-1})$, it holds that
	\begin{equation}
	\label{co2}
	\norm{\left[\D^{\alpha}, g, h\right]}_{L^2_t(L^2)}\lesssim\|\D g\|_{L^\infty_t(L^\infty)}\| \D h\|_{L^2_t(\MH^{m-2})}+\|\D h\|_{L^\infty_t(L^\infty)}\| \D g\|_{L^2_t(\MH^{m-2})}
	\end{equation}
where $\left[\D^{\alpha}, g, h\right]:= \D^{\alpha}(gh)-\D^{\alpha}g h-g\D^{\alpha} h.$
\end{itemize}
\end{lemma}
\begin{lemma}
	Let $g\in H^1([0,t];L^2)$. Then, we have
	\begin{equation}\label{L-infL-2}
		\|g(t)\|_0^2\ls t\|\D_tg\|_{L^2_t(L^2)}^2+\|g(0)\|_0^2.
	\end{equation}
\end{lemma}
\begin{proof}
	Since $g\in H^1([0,t];L^2)$, then $g(t,x)\in C([0,t];L^2)$. The fundamental theorem of calculus gives
	\begin{equation*}
		g(t,x)=g(0, x)+\int_{0}^{t}g_t(\tau,x)d\tau
	\end{equation*}
Then, by applying Minkowski and H\"{o}lder's inequalities, one can get \eqref{L-infL-2}.
\end{proof}
\begin{lemma}
	We have the anisotropic Sobolev embedding:
	\begin{equation}\label{ani-Sob-emb1}
		\|g\|_{L^\infty(\Omega)}^2\ls\|\DD g\|_{H^1_{tan}}^2+\|g\|_{H^2_{tan}}^2 
	\end{equation}
As a consequence, we also have
\begin{equation}\label{ani-Sob-em2}
	\|g\|_{L^\infty_{t,x}}^2\ls\|\DD g\|_{L^\infty_t(H^1_{tan})}^2+\|g\|_{L^\infty_t(H^2_{tan})}^2.
\end{equation}
\end{lemma}
We will also use the following lemma.
\begin{lemma}
	For $g\in H^1(\Gamma)$, and $h\in H^{\frac{1}{2}}(\Gamma)$ or $h\in H^{-\frac{1}{2}}(\Gamma)$, it holds that
	\begin{equation}
	\label{co123}
	\abs{gh}_{\frac{1}{2}} \ls \abs{g}_{1}\abs{h}_{\hal}, \abs{gh}_{-\hal}\ls \abs{g}_{1}\abs{h}_{-\hal}.
	\end{equation}
\end{lemma}
\begin{proof}
	It is direct to check that $\abs{gh}_{s} \ls \abs{g}_{1}\abs{h}_{s}$ for $s=0,1$ with the help of the Sobolev embedding $\abs{f}_{L^\infty}\ls \abs{f}_1$. Then the estimate \eqref{co123} follows by the interpolation. The second inequality follows by the dual estimate.
\end{proof}
The next lemma recall the embedding of fractional Sobolev spaces, refer to \cite{dinezza2012} for the proof.
\begin{lemma}
	Let $\Gamma$ be a bounded domain in $\mathbb{R}^d$. Then, for any $g\in H^s(\Gamma)$, we have the following Sobolev embedding:
	\begin{equation}
		|g|_{L^r}\ls |g|_s,
	\end{equation}
for any $r\in[1,\frac{2d}{d-2s}]$, if $2s<d$; for any $r\in[1,\infty)$, if $2s=d$. In particular, for $d=1,2$, we have
\begin{equation}\label{fra-emb}
	|g|_{L^4}\ls |g|_\frac{1}{2}.
\end{equation}
\end{lemma}

%
%

\subsection{Trace estimates}
First, we have the trace estimates, whose proof is given in \cite{gu2020a}
\begin{lemma}
	For $g\in H^1(\Omega)$, it holds that
	\begin{equation}\label{tre}
	\abs{g}_0^2\ls \norm{g}_0^2+\norm{g}_0\norm{\nabla g}_0
	\end{equation}
\end{lemma}
Next, we have the following normal trace estimates.
\begin{lemma}\label{normal trace}
	It holds that
	\begin{equation}\label{gga}
	\abs{\omega_ia_{i2}}_{-\hal}\ls \norm{\nabla\eta}_{L^\infty}\norm{\omega}_0+\norm{a_{ij}\D_j\omega_i}_0.
	\end{equation}
	
\end{lemma}
\begin{proof}
	We refer the reader to \cite[Section 5.9]{wang2021}.
\end{proof}
We also need the following lemma 
\begin{lemma}
	For any $g\in H^\hal(\Gamma)$, it holds  that
	\begin{equation}\label{tre-2}
		|\D_1g|_{-\hal}\ls |g|_\hal.
	\end{equation}
\end{lemma}
\subsection{Korn's inequality}
We refer to \cite{masmoudi2017,gu2020a} for the following Korn's type inequality.
\begin{lemma}
	For any $\boldsymbol{f}\in H^1(\Omega)$, it holds that
	\begin{equation}\label{Korn's ineq}
		\|\DD\boldsymbol{f}\|_0^2\ls P(\|\DD\boe\|_{2}^2)(\|\boldsymbol{S}_{\bA}(\boldsymbol{f})\|_0^2+\|\boldsymbol{f}\|_0^2).
	\end{equation}
\end{lemma}

\subsection{Geometric identities}
Differentiating $J$, $\bA$ and $\boldsymbol{a}$, we obtain
\begin{equation}\label{Geo-iden-1}
\D J=a_{ij}\D_j\D\eta_i,\quad \D A_{kj}=-A_{kl}\D_l \D\eta_i A_{ij},\quad \D a_{kj}=a_{li}\D_i\D\eta_lA_{kj}-a_{kl}\D_l\D\eta_iA_{ij}
\end{equation}
\section{Viscosity-independent A priori estimates}
In this section, we derive the $\vep$-independent estimates of smooth solutions to \eqref{FEL-final}, which is stated in the following proposition.
\begin{proposition}\label{uniform estimates}
Let $(\boe^\vep,\bv^\vep)$ be s solution to \eqref{FEL-final}. Then there exists a time $T$ independent of $\vep$ such that
\begin{equation}
	\sup_{t\in[0,T]}\mathfrak{E}^\vep(t)\leq 2M_0,
\end{equation}
where $M_0=P(\mathfrak{E}(0))$.
\end{proposition}
The proof of proposition can be divided into proofs of the following lemmas. For notational simplification, we only address the superscript $\vep$ of $\bv^\vep, \boe^\vep,\boldsymbol{a}^\vep,$ etc. in the statements of lemmas but omit it in the proof without causing any confusion.
Since
\begin{equation}\label{H1}
\tilde{\rho}_0^\vep\geq c_0,~~ \frac{1}{c_0}>J^\vep_0\geq c_0>0,
\end{equation}
for some $c_0>0$, we can assume that there exist a sufficiently small $T_\vep$ such that, for $t\in[0,T_\varepsilon]$, 
\begin{equation}\label{A-priori-assum}
|J^\vep(t)-J^\vep_0|\leq \frac{1}{8}c_0,\quad|\D_j\eta_i^\vep(t)-\D_j\eta_{0i}^\vep|\leq\frac{1}{8}c_0.
\end{equation}
We now derive $\vep$-independent estimates of smooth solutions $(\bv^\vep, \boe^\vep)$ to \eqref{FEL-final} under the a priori assumption \eqref{A-priori-assum}. We start from the basic energy estimates as follows.
\subsection{Basic energy estimates}
\begin{lemma}\label{basic-energy-est}
	For any $t\in[0,T_\vep]$, it holds that
	\begin{equation}\label{Basic energy}
	\|\bv^\vep(t)\|_0^2+\|\DD\boe^\vep(t)\|_0^2+\|Q(f^\vep)(t)\|_{L^1}+\sigma|\D_1\boe^\vep(t)|_{L^1}\lesssim M_0+T_\vep
	\end{equation}
\end{lemma}
\begin{proof}
	Taking the $L^2(\Omega)$ inner product of $\eqref{FEL-final}_2$ with $v_i$ gives
	\begin{align*}
		&\frac{1}{2}\frac{d}{dt}\int_\Omega\tilde{\rho}_0|\bv|^2 dx+\int_\Omega a_{ik}\D_kqv_i dx-2\mu\vep\int_\Omega a_{kl}\D_l(S_Av)_{ik}v_idx\\
	&-\lambda\vep\int_\Omega a_{ij}\D_j(A_{kl}\D_l v_k)v_i-\int_\Omega\D_j(\tilde{\rho}_0\D_j\eta_i) v_i dx=0.
	\end{align*}
	By the integration by parts and using \eqref{Piola}, $\eqref{FEL-final}_3$, \eqref{Jacobi's formula}, we have
	\begin{align*}
	&\int_\Omega a_{ik}\D_kqv_idx-2\mu\vep\int_\Omega a_{kl}\D_l(S_Av)_{ik}v_idx-\lambda\vep\int_\Omega a_{ij}\D_j(A_{kl}\D_l v_k)v_i-\int_\Omega\D_j(\tilde{\rho}_0\D_j\eta_i) v_idx\\
	&=-\int_\Omega a_{ik}\D_kv^iqdx+2\mu\vep\int_\Omega(S_Av)_{ik}a_{kl}\D_lv_idx+\lambda\vep\int_\Omega J(A_{ij}\D_j v_i)^2dx+\int_\Omega\tilde{\rho}_0\D_j\eta_i\D_jv_i dx\\
	&\quad+\int_\Gamma (qa_{i2}-2\mu\vep(S_Av)^\vep_{ik}a_{k2}-\lambda\vep A^\vep_{kl}\D_lv^\vep_ka_{i2}-\tilde{\rho}_0\D_2\eta_i)v_i\\
	&=-\int_\Gamma\sigma \D_1(\frac{\D_1\eta_i}{|\D_1\boldsymbol{\eta}|})v_i-\int_\Omega J_tqdx+\frac{1}{2}\frac{d}{dt}\int_\Omega\tilde{\rho}_0|\nabla \boldsymbol{\eta}|^2dx\\
	&\quad+2\mu\vep\int_\Omega(S_Av)_{ik}a_{kl}\D_lv_idx+\lambda\vep\int_\Omega J(A_{ij}\D_j v_i)^2dx\\
	&=\sigma\int_\Gamma \frac{\D_1\eta_i\D_1v_i}{|\D_1\boldsymbol{\eta}|}+\int_\Omega\tilde{\rho}_0f_tf^{-2}q(f)dx-\frac{d}{dt}\int_\Omega Jdx+\frac{1}{2}\frac{d}{dt}\int_\Omega\tilde{\rho}_0|\nabla \boldsymbol{\eta}|^2dx\\
	&\quad+2\mu\vep\int_\Omega(S_Av)_{ik}a_{kl}\D_lv_idx+\lambda\vep\int_\Omega J(A_{ij}\D_j v_i)^2dx\\
	&=\frac{d}{dt}\int_\Gamma \sigma|\D_1\boldsymbol{\eta}|+\frac{d}{dt}\int_\Omega\tilde{\rho}_0Q(f)dx+\frac{1}{2}\frac{d}{dt}\int_\Omega\tilde{\rho}_0|\nabla \boldsymbol{\eta}|^2dx-\frac{d}{dt}\int_\Omega Jdx\\
	&\quad+2\mu\vep\int_\Omega(S_Av)_{ik}a_{kl}\D_lv_idx+\lambda\vep\int_\Omega J(A_{ij}\D_j v_i)^2dx,
	\end{align*}
	where $Q(f)=\int_{1}^{f}q(\mu)\mu^{-2}d\mu$. Therefore, we obtain
	\begin{align*}
	&\frac{d}{dt}\left(\frac{1}{2}\int_\Omega\tilde{\rho}_0(|\bv|^2+|\nabla\boldsymbol{\eta}|^2+2Q(f))+\int_\Gamma\sigma|\D_1\boldsymbol{\eta}|-\int_\Omega Jdx\right)\\
	&+2\mu\vep\int_\Omega J\left|S_Av\right|^2\,dx+\lambda\vep\int_\Omega J\left|\Div_A v\right|^2dx=0
	\end{align*}
	which, following integration in time, yields
	\begin{align*}
	&\left(\frac{1}{2}\int_\Omega\tilde{\rho}_0(|\bv|^2+|\nabla\boldsymbol{\eta}|^2+2Q(f))+\int_\Gamma\sigma|\D_1\boldsymbol{\eta}|\right)(t)\\
	&\quad+2\mu\vep\int_\Omega J\left|S_Av\right|^2\,dx+\lambda\vep\int_\Omega J\left|\Div_A v\right|^2dx\\
	&=\frac{1}{2}\int_\Omega\tilde{\rho}_0(|\bv_0|^2+|\nabla\boldsymbol{\eta}_0|^2+2Q(f_0))+\int_\Gamma\sigma|\D_1\boldsymbol{\eta}_0|+\int_\Omega(J(t)-J_0)dx.\nonumber
	\end{align*}
    Noticting that by straightforward calculation, 
    \begin{equation*}
            \begin{split}
                    2\mu\left|S_A(\boldsymbol{f})\right|^2+\lambda \left|\Div_A (\boldsymbol{f})\right|^2 = &(2\mu+2\lambda)\left(S_A(\boldsymbol{f})_{11}^2+S_A(\boldsymbol{f})_{22}^2\right)+2\mu\left(S_A(\boldsymbol{f})_{12}^2+S_A(\boldsymbol{f})_{21}^2\right)\\&-\lambda\left(A_{1\ell}\D_\ell \boldsymbol{f}_1-A_{2\ell}\D_\ell \boldsymbol{f}_2\right)^2,
    \end{split}
    \end{equation*}
    by $2\mu+2\lambda>0$, $\mu>0$, we have
    \begin{equation}
            \label{kon}
            2\mu\left|S_A(\boldsymbol{f})\right|^2+\lambda \left|\Div_A (\boldsymbol{f})\right|^2\geq c\left|S_A(\boldsymbol{f})\right|^2.
    \end{equation}
    Thus, in view of \eqref{H1}, \eqref{A-priori-assum} and Korn's inequality \eqref{Korn's ineq}, \eqref{Basic energy} holds true.
\end{proof}

\subsection{Estimates of $J$ and $q$}
Before performing the higher-order estimates, we first prove the following key lemma for the estimate of  velocity divergence and pressure. 
\begin{lemma}\label{est-J-q}
	For any $t\in[0,T_\vep]$, $m\geq 3$, it holds that
	\begin{align}
		&\|J^\vep\|_{L^2_t(\MH^m)}^2+\|q^\vep\|_{L^2_t(\MH^m)}^2\ls T_\vep P\left(\sup_{t\in[0,T_\vep]}\mathfrak{E}^\vep(t)\right),\label{m-J-q}\\
		&\|J^\vep\|_{L^\infty_t(\MH^{m-1})}^2+\|q^\vep\|_{L^\infty_t(\MH^{m-1})}^2\ls M_0+T_\vep^2P\left(\sup_{t\in[0,T_\vep]}\mathfrak{E}^\vep(t)\right).\label{m-1-J-q2}
	\end{align}
\end{lemma}
\begin{proof}
	By using \eqref{Geo-iden-1}, \eqref{A-priori-assum} and \eqref{pro-est}, one has
	\begin{align}\label{J-est}
		\|J\|_{L^2_t(\MH^m)}^2&\ls \|J\|_{L^2_t(L^2)}^2+\sum_{|\alpha|\leq m-1}\|\D^\alpha(a_{ij}\D_j\D\eta_i)\|_{L^2_t(L^2)}^2\nonumber\\
		&\ls c_0t+\sum\limits_{|\beta|+|\nu|=|\alpha|\leq m-1}\|\D^\beta a_{ij}\D^\nu\D\D_j\eta_i\|_{L^2_t(L^2)}^2\nonumber\\
		&\ls c_0t+\|a_{ij}\|_{L^\infty_{t,x}}^2\|\D_j\eta_i\|_{L^2_t(\MH^m)}^2+\|\D_j\eta_i\|_{L^\infty_{t,x}}^2\|a_{ij}\|_{L^2_t(\MH^m)}^2\\
		&\ls c_0t+\|\DD\boe\|_{L^\infty_{t,x}}^2\|\DD\boe\|_{L^2_t(\MH^m)}^2\ls TP\left(\sup_{t\in[0,T]}\mathfrak{E}(t)\right),\nonumber
	\end{align}
	where we have used
	\begin{equation}\label{gra-eta-est}
		\|\DD\boe\|_{L^\infty_{t,x}}\ls \|\DD\boe\|_{L^2_t(\MH^{m-1})}\ls T^\hal\|\boe\|_{L^\infty_t(\MH^m)}
	\end{equation}
	Since $J$ is bounded from below and above in \eqref{A-priori-assum}, direct calculations yield that, for any $s\in \mathbb{R}$,
	\begin{equation*}
		|\D^\alpha(J^s)|\ls \sum\limits_{|\beta_1|+\cdots+|\beta_k|=|\alpha|}\left|\D^{\beta_1}J\cdots\D^{\beta_k}J\right|
	\end{equation*}
	Then, it follows from Minkowski, H\"{o}lder and Gagliardo-Nirenberg's inequalities in $\Omega\times[0,t]$ that
	\begin{align}\label{Js-est}
		\|J^s\|_{L^2_t(\MH^m)}^2&\ls \sum\limits_{|\alpha|\leq m}\|\D^\alpha(J^s)\|_{L^2_t(L^2)}^2\ls\sum\limits_{|\alpha|\leq m}\sum\limits_{|\beta_1|+\cdots+|\beta_k|=|\alpha|}\|\D^{\beta_1}J\cdots\D^{\beta_k}J\|_{L^2_t(L^2)}^2\nonumber\\
		&\ls \sum\limits_{|\alpha|\leq m}\sum\limits_{|\beta_1|+\cdots+|\beta_k|=|\alpha|}\|\D^{\beta_1}J\|_{L^{r_1}_{t,x}}^2\cdots\|\D^{\beta_k}J\|_{L^{r_k}_{t,x}}^2\nonumber\\
		&\ls\sum\limits_{|\alpha|\leq m}\sum\limits_{|\beta_1|+\cdots+|\beta_k|=|\alpha|}\prod\limits_{i=1}^k\left(\|\D^{|\alpha|}J\|_{L^2_{t,x}}^{2\theta_i}\|J\|_{L^\infty_{t,x}}^{2(1-\theta_i)}+\|J\|_{L^2_{t,x}}^2\right)\\
		&\ls \sum\limits_{|\alpha|\leq m}\|\D^{|\alpha|}J\|_{L^2_{t,x}}^2+P\left(\|J\|_{L^2_{t,x}}^2\right)\ls \|J\|_{L^2_t(\MH^m)}^2+
		P\left(\|J\|_{L^2_{t,x}}^2\right)\nonumber\\
		&\ls TP\left(\sup_{t\in[0,T]}\mathfrak{E}(t)\right) .\nonumber
	\end{align}
	where $r_i$ and $\theta_i$, $i=1,\cdots,k$ satisfy 
	\begin{equation*}
		\sum_{i=1}^{k}\frac{1}{r_i}=\frac{1}{2},\quad\frac{1}{r_i}=\frac{|\beta_i|}{3}+\left(\hal-\frac{|\alpha|}{3}\right)\theta_i+\frac{1-\theta_i}{\infty},\text{ so that }\sum_{i=1}^{k}\theta_i=1.
	\end{equation*}
	By using \eqref{L-infL-2}, one can get
	\begin{equation*}
		\|J\|_{L^\infty(\MH^{m-1})}^2\ls \|J_0\|_{\MH^m}^2+t\|\D_t J\|_{L^2_t(\MH^{m-1})}^2\ls M_0+T^2P\left(\sup_{t\in[0,T]}\mathfrak{E}(t)\right).  
	\end{equation*}
	For the pressure $q$, we get from \eqref{gamma-law} and \eqref{f-eqn} that
	\begin{equation*}
		q=Af^\gamma+1-p_e=A\tilde{\rho}_0^\gamma J^{-\gamma}+1-p_e.
	\end{equation*}
	Then, it follows from \eqref{pro-est} and \eqref{Js-est} that
	\begin{align*}
		\|q\|_{L^2_t(\MH^m)}^2&\ls \|\tilde{\rho}^\gamma\|_{L^\infty_{t,x}}^2\|J^{-\gamma}\|_{L^2_t(\MH^m)}^2+\|J^{-\gamma}\|_{L^\infty_{t,x}}^2\|\tilde{\rho}^\gamma\|_{L^2_t(\MH^m)}^2+T\ls TP\left(\sup_{t\in[0,T]}\mathfrak{E}(t)\right).  
	\end{align*}
	By using \eqref{L-infL-2} that
	\begin{equation*}
		\|q\|_{L^\infty_t(\MH^{m-1})}^2\ls \|q_0\|_{\MH^{m-1}}+t\|\D_tq\|_{L^2_t(\MH^{m-1})}^2\ls M_0+T^2P\left(\sup_{t\in[0,T]}\mathfrak{E}(t)\right).
	\end{equation*}
\end{proof}

\subsection{Tangential derivative estimates with at least one spatial derivative}
We now derive the high order estimates of tangential derivatives. We have to separate the estimate of fully temporal derivative from other tangential ones since we lose some key estimates for the fully temporal derivative case. The following lemma states the estimate of tangential derivatives with at least one spatial one.
\begin{lemma}\label{Tan-est}
	For any $t\in[0,T_\vep]$, $m\geq 4$, and any $\beta=(\beta_0,\beta_1)$ with $0\leq\beta_0+\beta_1\leq m-1$,   it holds that
	\begin{align}\label{tan-est-1}
	&\|\bv^\vep\|_{L^\infty_t(\MH_{tan}^{m,1})}^2+\|q^\vep\|_{L^\infty_t(\MH_{tan}^{m,1})}^2+\|\nabla\boe^\vep\|_{L^\infty_t(\MH_{tan}^{m,1})}^2+\sum_{|\beta|\leq m-1}|\bar{\D}^\beta\D_1^2\eta_i^\vep a_{i2}|_0^2+\vep\|\DD\bv^\vep\|_{L^2_t(\MH_{tan}^{m,1})}^2\nonumber\\
	&\ls M_0+\delta\sup_{t\in[0,T_\vep]}\mathfrak{E}(t)+T_\vep^\frac{1}{4}P\left(\sup_{t\in[0,T]}\mathfrak{E}(t)\right).
	\end{align}
\end{lemma}
\begin{proof}
	For $1\leq|\alpha|\leq m$, where $\alpha=(\alpha_0,\alpha_1)$ with $\alpha_1\geq 1$,  applying $\bar{\D}^\alpha$ to the second equation of \eqref{FEL-final}, and then taking the $L^2(\Omega)$ inner product with $\bar{\D}^\alpha v_i$ yield
	\begin{align*}
	&\frac{1}{2}\frac{d}{dt}\int_\Omega\tilde{\rho}_0|\bar{\D}^\alpha \bv|^2dx+\int_\Omega a_{ik}\D_k \bar{\D}^\alpha q\bar{\D}^\alpha v_idx-2\mu\vep\int_\Omega \bar{\D}^\alpha(a_{kl}\D_l(S_Av)_{ik})\bar{\D}^\alpha v_idx\\
	&-\lambda\vep\int_\Omega \bar{\D}^\alpha(a_{ij}\D_j\Div_Av)\bar{\D}^\alpha v_idx-\int_\Omega \D_j(\bar{\D}^\alpha(\tilde{\rho}_0\D_j\eta_i))\bar{\D}^\alpha v_idx\\
	&=-\underbrace{\int_\Omega[\bar{\D}^\alpha,\tilde{\rho}_0]\D_tv_i\bar{\D}^\alpha v_idx}_{R_\eta^1}-\underbrace{\int_\Omega[\bar{\D}^\alpha,a_{ik}]\D_kq\bar{\D}^\alpha v_idx}_{R_q^1}.\nonumber
	\end{align*}
	By integration by parts and using \eqref{Piola}, the third equation of \eqref{FEL-final}, \eqref{Jacobi's formula}, we have
	\begin{align*}
	&\quad\int_\Omega a_{ik}\D_k \bar{\D}^\alpha q\bar{\D}^\alpha v_idx-2\mu\vep\int_\Omega \bar{\D}^\alpha(a_{kl}\D_l(S_Av)_{ik})\bar{\D}^\alpha v_idx\\
	&\quad-\lambda\vep\int_\Omega \bar{\D}^\alpha(a_{ij}\D_j\Div_A (v))\bar{\D}^\alpha v_idx-\int_\Omega \D_j(\bar{\D}^\alpha(\tilde{\rho}_0\D_j\eta_i))\bar{\D}^\alpha v_idx\\
	&=-\int_\Omega a_{ik}\bar{\D}^\alpha q\bar{\D}^\alpha \D_kv_idx+2\mu\vep\int_\Omega\bar{\D}^\alpha((S_Av)_{ik}a_{kl})\D_l\bar{\D}^\alpha v_idx\\
	&\quad+\lambda\vep\int_\Omega \bar{\D}^\alpha(a_{ij}\Div_Av)\D_j\bar{\D}^\alpha v_idx+\int_\Omega \bar{\D}^\alpha(\tilde{\rho}_0\D_j\eta_i)\D_j\bar{\D}^\alpha v_idx\\
	&\quad+\int_\Gamma\left(\bar{\D}^\alpha qa_{i2}-\bar{\D}^\alpha\left(2\mu\vep a_{k2} (S_Av)_{ik}+\lambda\vep a_{i2}(\Div_Av)+\tilde{\rho}_0\D_2\eta_i\right)\right)\bar{\D}^\alpha v_i\\
	&=-\int_\Omega\bar{\D}^\alpha q\bar{\D}^\alpha(a_{ik}\D_kv_i)dx+\underbrace{\int_\Omega\bar{\D}^\alpha q[\bar{\D}^\alpha,a_{ik}]\D_kv_idx}_{R_\eta^2}+2\mu\vep\int_\Omega J|\boldsymbol{S}_{\bA}(\bar{\D}^\alpha \bv)|^2dx\\
	&\quad+\lambda\vep\int_\Omega J(\Div_A(\bar{\D}^\alpha v))^2dx+\underbrace{2\mu\vep\int_\Omega[\bar{\D}^\alpha, a_{kl}] (S_Av)_{ik}\D_l\bar{\D}^\alpha v_idx}_{R_\vep^1}\\
	&\quad+\underbrace{\lambda\vep\int_\Omega[\bar{\D}^\alpha,a_{ij}](\Div_A v)\D_j\bar{\D}^\alpha v_idx}_{R_\vep^2}+\underbrace{2\mu\vep\int_\Omega[\bar{\D}^\alpha, A_{ij}]\D_jv_kJS_A(\bar{\D}^\alpha v)_{ik}}_{R_\vep^3}\\
	&\quad+\underbrace{\lambda\vep\int_\Omega[\bar{\D}^\alpha, A_{ij}]\D_j v_iJ\Div_A(\bar{\D}^\alpha v)}_{R_\vep^4}+\int_\Omega\tilde{\rho}_0\D_j\bar{\D}^\alpha\eta_i\D_j\bar{\D}^\alpha v_idx+\underbrace{\int_\Omega[\bar{\D}^\alpha,\tilde{\rho}_0]\D_j\eta_i\D_j\bar{\D}^\alpha v_idx}_{R_\eta^3}\\
	&\quad+\int_\Gamma\left(\bar{\D}^\alpha qa_{i2}-\bar{\D}^\alpha\left(2\mu\vep a_{k2} (S_Av)_{ik}+\lambda\vep a_{i2}(\Div_Av)+\tilde{\rho}_0\D_2\eta_i\right)\right)\bar{\D}^\alpha v_i\\
	&=\frac{1}{2}\frac{d}{dt}\int_\Omega\tilde{\rho}_0|\bar{\D}^\alpha\nabla\boe|^2dx-\int_\Omega\bar{\D}^\alpha q\bar{\D}^\alpha\D_tJdx+2\mu\vep\int_\Omega J|\boldsymbol{S}_{\bA}(\bar{\D}^\alpha \bv)|^2dx+\lambda\vep\int_\Omega J(\Div_A(\bar{\D}^\alpha v))^2dx\\
	&\quad+\int_\Gamma\left(\bar{\D}^\alpha qa_{i2}-\bar{\D}^\alpha\left(2\mu\vep a_{k2} (S_Av)_{ik}+\lambda\vep a_{i2}(\Div_Av)+\tilde{\rho}_0\D_2\eta_i\right)\right)\bar{\D}^\alpha v_i+R_\eta^2+R_\eta^3+\sum_{i=1}^4R_\vep^i.
	\end{align*}
Moreover, \eqref{f-eqn} implies
\begin{equation}\label{q-eqn}
	\D_tJ=\D_t(\tilde{\rho}_0f^{-1})=-\tilde{\rho}_0\frac{\D_tf}{f^2}=-\tilde{\rho}_0\frac{\D_tq}{q'(f)f^2}=-\frac{\tilde{\rho}_0\D_tq}{\gamma Af^{\gamma+1}}=-\frac{J^{\gamma+1}}{\gamma A\tilde{\rho}_0^\gamma}\D_tq.
\end{equation}
so that
\begin{align*}\label{I-est}
	&-\int_\Omega\bar{\D}^\alpha q\bar{\D}^\alpha\D_tJdx=\frac{1}{\gamma }\int_\Omega\bar{\D}^\alpha q\bar{\D}^\alpha(\tilde{\rho}_0^{-\gamma}J^{\gamma+1}\D_t q)dx\nonumber\\
	&=\frac{1}{\gamma }\int_\Omega\tilde{\rho}_0^{-\gamma}J^{\gamma+1}\bar{\D}^\alpha q\D_t\bar{\D}^\alpha qdx+\frac{1}{\gamma }\int_\Omega\bar{\D}^\alpha q[\bar{\D}^\alpha,\tilde{\rho}_0^{-\gamma}J^{\gamma+1}]\D_tqdx\\
	&=\frac{1}{2\gamma }\frac{d}{dt}\int_\Omega\tilde{\rho}_0^{-\gamma}J^{\gamma+1}|\bar{\D}^\alpha q|^2dx\underbrace{-\frac{\gamma+1}{2\gamma }\int_\Omega\tilde{\rho}_0^{-\gamma}J^\gamma J_t|\bar{\D}^\alpha q|^2dx}_{R_q^2}+\underbrace{\frac{1}{\gamma }\int_\Omega\bar{\D}^\alpha q[\bar{\D}^\alpha,\tilde{\rho}_0^{-\gamma}J^{\gamma+1}]\D_tqdx}_{R_q^3}.\nonumber
\end{align*}
	Therefore, combining the above equations, we can obtain that
	\begin{equation}\label{tang-est}
	\begin{aligned}
		&\frac{1}{2}\frac{d}{dt}\int_\Omega\tilde{\rho}_0(|\bar{\D}^\alpha \bv|^2+|\bar{\D}^\alpha\nabla\boe|^2)dx+\frac{1}{2\gamma }\frac{d}{dt}\int_\Omega\tilde{\rho}_0^{-\gamma}J^{\gamma+1}|\bar{\D}^\alpha q|^2dx\\
		&\quad+2\mu\vep\int_\Omega J|\boldsymbol{S}_{\bA}(\bar{\D}^\alpha \bv)|^2dx+\lambda\vep\int_\Omega J(\Div_A(\bar{\D}^\alpha v))^2dx\\
		&\quad+\underbrace{\int_\Gamma\left(\bar{\D}^\alpha qa_{i2}-\bar{\D}^\alpha\left(2\mu\vep a_{k2} (S_Av)_{ik}+\lambda\vep a_{i2}(\Div_Av)+\tilde{\rho}_0\D_2\eta_i\right)\right)\bar{\D}^\alpha v_i}_{R_b}\\
		&=-\left(\sum_{i=1}^{3}R_\eta^i+\sum_{i=1}^{3}R_q^i+\sum_{i=1}^{4}R_\vep^i\right).
	\end{aligned}	
	\end{equation}	
We now estimate the terms $R_\eta^i$, $R_q^i$ and $R_\vep^i$ one by one.
	It follows from \eqref{co1} that
	\begin{equation}\label{R1-est}
		\begin{aligned}
				\left|\int_{0}^{t}R_\eta^1\right|&\lesssim\|[\bar{\D}^\alpha,\tilde{\rho}_0]\D_t v_i\|_{L^2_{t,x}}\|\bar{\D}^\alpha v_i\|_{L^2_{t,x}}\\
				&\ls(\|\D\tilde{\rho}_0\|_{L^\infty_{t,x}}\|\D_t\bv\|_{L^2_t(\MH^{m-1})}+\|\D_t\bv\|_{L^\infty_{t,x}}\|\D\tilde{\rho}_0\|_{L^2_t(\MH^{m-2})})\|\bar{\D}^\alpha \bv\|_{L^2_{t,x}}\\
				&\lesssim\|\bv\|_{L^2_t(\MH^m)}t^\frac{1}{2}\|\bar{\D}^\alpha \bv\|_{L^\infty_tL^2_{x}}\lesssim T^\frac{1}{2}\sup_{t\in[0,T]}\mathfrak E^\epsilon(t),
		\end{aligned}
	\end{equation}
where we have used $\|\D_t\bv\|_{L^\infty_{t,x}}\ls \|\D_t\bv\|_{L^2_t(\MH^{m-1})}$ from Sobolev's embedding.
	Similarly,  we have
		\begin{equation}\label{R2-est}
		\left|\int_{0}^{t}R_q^1\right|\lesssim\|\bar{\D}a_{ik}\|_{L^2_t(\MH^{m-1})}\|\D_kq\|_{L^2_t(\MH^{m-1})}t^\frac{1}{2}\|\bar{\D}^\alpha v_i\|_{L^\infty_tL^2_{x}}\lesssim TP\left(\sup_{t\in[0,T]}\mathfrak E^\epsilon(t)\right)
	\end{equation}
		\begin{equation}\label{R3-est}
		\left|\int_{0}^{t}R_\eta^2\right|\lesssim\|\bar{\D}a_{ik}\|_{L^2_t(\MH^{m-1})}\|\D_kv_i\|_{L^2_t(\MH^{m-1})}t^\frac{1}{2}\|\bar{\D}^\alpha q\|_{L^\infty_tL^2_{x}}\lesssim TP\left(\sup_{t\in[0,T]}\mathfrak E^\epsilon(t)\right),
	\end{equation}
and
	\begin{align}\label{R4-est}
			\left|\int_{0}^{t}R_\eta^3\right|&\lesssim\left|\int_\Omega[\bar{\D}^\alpha,\tilde{\rho}_0]\D_j\eta_i\D_j\bar{\D}^\alpha \eta_i\Big|_0^t\right|+\left|\int_{0}^{t}\int_\Omega[\bar{\D}^\alpha,\tilde{\rho}_0]\D_t\D_j\eta_i\D_j\bar{\D}^\alpha \eta_i\,dxd\tau\right|\nonumber\\
			&\ls M_0+\delta\|\bar{\D}^\alpha\nabla\boe(t)\|_{0}^2+T\|\D_t\DD\boe\|_{L^2_t(\MH^{m-1})}^2+\|\D_t\nabla\boe\|_{L^2_t(\MH^{m-1})}\|\nabla\boe\|_{L^2_t(\MH^{m,1}_{tan})}\\
			&\ls M_0+\delta\|\bar{\D}^\alpha\nabla\boe(t)\|_{0}^2+T^\frac{1}{2}\sup_{t\in[0,T]}\mathfrak E^\epsilon(t),\nonumber
	\end{align}	
where, in \eqref{R4-est}, we have used $\|\tilde{\rho}_0\|_{H^m}\leq C$ and
\begin{equation*}
	\begin{aligned}
		&\left|\int_\Omega[\bar{\D}^\alpha,\tilde{\rho}_0]\D_j\eta_i(t)\D_j\bar{\D}^\alpha \eta_i(t)\right|=\left|\int_\Omega\sum\limits_{|\beta|\geq 1,|\beta|+|\nu|=m}C_\alpha^\beta\bar{\D}^\beta\tilde{\rho}_0\bar{\D}^\nu\D_j\eta_i\bar{\D}^\alpha\D_j\eta_i\right|\\
		&\ls \delta\|\bar{\D}^\alpha\nabla\boe(t)\|_{0}^2+C_\delta\left(\|\D\tilde{\rho}_0\|_{L^\infty}^2\|\bar{\D}^{m-1}\DD\boe(t)\|_0^2+\|\D^m\tilde{\rho}_0\|_0^2\|\DD\boe(t)\|_{L^\infty}^2\right)\\
		&\quad+C_\delta\sum\limits_{|\beta|\geq 2,|\beta|+|\nu|=m}\|\D^{\beta}\tilde{\rho}_0\|_{L^4}^2\|\bar{\D}^\nu\DD \boe(t)\|_{L^4}^2\nonumber\\
		&\ls\delta\|\bar{\D}^\alpha\nabla\boe(t)\|_{0}^2+C_\delta\|\DD\boe(t)\|_{\MH^{m-1}}^2\\
		&\ls \delta\|\bar{\D}^\alpha\nabla\boe(t)\|_{0}^2+M_0+T\|\D_t\DD\boe\|_{L^2_t(\MH^{m-1})}^2.
	\end{aligned}
\end{equation*}
In view of \eqref{H1}, \eqref{A-priori-assum}, it is obvious that
	\begin{equation}\label{Rq1-est}
	\left|\int_{0}^{t}R_{q}^2\right|\lesssim T\sup_{t\in[0,T]}\mathfrak E^\epsilon(t)
	\end{equation} 
It follows from \eqref{co1} and \eqref{m-J-q} that
	\begin{align}\label{Rq2-est}
	\left|\int_{0}^{t}R_{q}^3\right|&\ls \|\bar{\D}^\alpha q\|_{L^2_t(L^2)}\|[\bar{\D}^\alpha,\tilde{\rho}_0^{-\gamma}J^{\gamma+1}]\D_tq\|_{L^2_t(L^2)}\nonumber\\
	&\ls \|\bar{\D}^\alpha q\|_{L^2_t(L^2)}(\|\bar{\D}(\tilde{\rho}_0^{-\gamma}J^{\gamma+1})\|_{L^\infty}\|\D_tq\|_{L^2_t(\MH^{m-1})}+\|\D_tq\|_{L^\infty_{t,x}}\|\bar{\D}(\tilde{\rho}_0^{-\gamma}J^{\gamma+1})\|_{L^2_t(\MH^{m-1})})\nonumber\\
	&\lesssim\|\bar{\D}^\alpha q\|_{L^2_t(L^2)}\|\D_tq\|_{L^2_t(\MH^{m-1})}\|\bar{\D}(\tilde{\rho}_0^{-\gamma}J^{\gamma+1})\|_{L^2_t(\MH^{m-1})}\lesssim T^\frac{1}{2}P\left(\sup_{t\in[0,T]}\mathfrak E^\epsilon(t)\right).
	\end{align}
In view of \eqref{A-priori-assum}, Sobolev interpolation inequalities and the following estimate from \eqref{ani-Sob-em2}
\begin{equation}\label{ubd}
	\|\bar{\D}\DD\boe\|_{L^\infty_{t,x}}\ls \|\bar{\D}\DD\boe\|_{L^\infty_t(H^2_{tan})}+\|\bar{\D}\DD^2\boe\|_{L^\infty_t(H^1_{tan})}
\end{equation}
one has
\begin{align}\label{Rvep1}
	\left|\int_{0}^{t}R_\vep^1+R_\vep^2\right|&\ls\vep\|[\bar{\D}^\alpha,\boldsymbol{a}](\bA\DD\bv)\|_{L^2_t(L^2)}\|\DD\bar{\D}^\alpha\bv\|_{L^2_t(L^2)}\nonumber\\
	&\ls \vep\left(\|\bA\DD\bv\|_{L^\infty}\|\D\boldsymbol{a}\|_{L^2_t(\MH^{m-1})}+\|\D\boldsymbol{a}\|_{L^\infty_{t,x}}\|\bA\DD\bv\|_{L^2_t(\MH^{m-1})}\right)\|\DD\bar{\D}^\alpha\bv\|_{L^2_t(L^2)}\\
	&\ls \vep\left(\|\DD\boe\|_{L^\infty_{t,x}}\|\D\DD\boe\|_{L^\infty_{t,x}}\|\DD\boe\|_{L^2_t(\MH^m)}+\|\D\DD\boe\|_{L^\infty_{t,x}}^2\|\bA\|_{L^2_t(\MH^{m-1})}\right)\|\DD\bar{\D}^\alpha\bv\|_{L^2_t(L^2)}\nonumber\\
	&\ls \delta\vep\|\DD\bar{\D}^\alpha \bv\|_{L^2_t(L^2)}^2+TP\left(\sup_{t\in[0,T]}\mathfrak{E}(t)\right).\nonumber
\end{align}
where we have used \eqref{gra-eta-est}
and
\begin{equation}\label{2nd-gra-eta-est}
	\|\D\DD\boe\|_{L^\infty_{t,x}}\ls \|\D\DD\boe\|_{L^2_t(\MH^{m-1})}\ls\|\DD\boe\|_{L^2_t(\MH^m)}
\end{equation}
Similarly, it follows from \eqref{co1} and \eqref{ubd} that
\begin{align}\label{Rvep2}
	\left|\int_{0}^{t}R_\vep^3+R_\vep^4\right|&\ls \vep\|[\bar{\D}^\alpha, \bA]\DD\bv\|_{L^2_t(L^2)}(\|\sqrt{J}\boldsymbol{S}_{\bA}(\bar{\D}^\alpha\bv)\|_{L^2_t(L^2)}+\|\sqrt{J}\Div_A(\bar{\D}v)\|_{L^2_t(L^2)})\nonumber\\
	&\ls\delta\vep(\|\sqrt{J}\boldsymbol{S}_{\bA}(\bar{\D}^\alpha\bv)\|_{L^2_t(L^2)}^2+\|\sqrt{J}\Div_A(\bar{\D}^\alpha v)\|_{L^2_t(L^2)}^2)+TP\left(\sup_{t\in[0,T]}\mathfrak{E}(t)\right).
\end{align}
Next, we estimate the boundary term $R_b$. Since
\begin{equation*}
\sigma\D_1\left(\frac{\D_1\eta_i}{|\D_1\boe|}\right)=\sigma\frac{\D_1^2\eta_ka_{k2}}{|\D_1\boe|^3}a_{i2},
\end{equation*}
then one gets from $\eqref{FEL-final}_3$ that
\begin{equation}\label{BB}
	\mathfrak{B}:=\sigma\frac{\D_1^2\eta_ka_{k2}}{|\D_1\boe|^3}+q=2\mu\vep\frac{a_{k2}A_{il}\D_lv_ka_{i2}}{|\D_1\boe|^2}+\lambda\vep\Div_Av+\frac{\tilde{\rho}_0J}{|\D_1\boe|^2}.
\end{equation}
then, we can obtain from $\eqref{FEL-final}_3$ and integration by parts that
\begin{align}\label{II-est}
	R_b=&-\int_\Gamma\sigma\bar{\D}^\alpha\left(\frac{\D_1^2\eta_ka_{k2}}{|\D_1\boe|^3}\right)a_{i2}\bar{\D}^\alpha v_i\underbrace{-\int_\Gamma\mathfrak{B}\bar{\D}^\alpha a_{i2}\bar{\D}^\alpha v_i}_{R_b^1}\underbrace{-\int_\Gamma\left[\bar{\D}^\alpha,\mathfrak{B},a_{i2}\right]\bar{\D}^\alpha v_i}_{R_b^2}\nonumber\\
	=&-\sigma\int_\Gamma\frac{1}{|\D_1\boe|^3}\bar{\D}^\alpha\D_1^2\eta_ka_{k2}a_{i2}\bar{\D}^\alpha v_i\underbrace{-\sigma\int_\Gamma\bar{\D}^\alpha\left(\frac{a_{k2}}{|\D_1\boe|^3}\right)\D_1^2\eta_k a_{i2}\bar{\D}^\alpha v_i}_{R_b^3}\nonumber\\
	&\underbrace{-\sigma\int_\Gamma\left[\bar{\D}^\alpha,\D_1^2\eta_k,\frac{a_{k2}}{|\D_1\boe|^3}\right]a_{i2}\bar{\D}^\alpha v_i}_{R_b^4}+R_b^1+R_b^2\\
	=&\sigma\int_\Gamma\frac{1}{|\D_1\boe|^3}\bar{\D}^\alpha\D_1\eta_k a_{k2}a_{i2}\bar{\D}^\alpha\D_t\D_1\eta_i\underbrace{+\sigma\int_\Gamma\D_1\left(\frac{a_{k2}a_{i2}}{|\D_1\boe|^3}\right)\bar{\D}^\alpha\D_1\eta_k\bar{\D}^\alpha v_i}_{R_b^5}+R_b^1+\cdots+R_b^4\nonumber\\
	=&\frac{\sigma}{2}\frac{d}{dt}\int_\Gamma\frac{|\bar{\D}^\alpha\D_1\eta_k a_{k2}|^2}{|\D_1\boe|^3}\underbrace{-\frac{\sigma}{2}\int_\Gamma\D_t\left(\frac{a_{k2}a_{i2}}{|\D_1\boe|^3}\right)\bar{\D}^\alpha\D_1\eta_k \bar{\D}^\alpha\D_1\eta_i}_{R_b^6}+R_b^1+\cdots+R_b^5.\nonumber
\end{align}
Before controlling the terms $R_b^i$, $i=1,\cdots,6$, we first derive estimate of $\mathfrak{B}$ as follows.
By using \eqref{BB}, \eqref{pro-est} and \eqref{A-priori-assum}, one has
\begin{align}\label{est-B}
	\|\mathfrak{B}\|_{L^2_t(\MH^{m})}^2&\ls \vep^2\left(\left\|\frac{a_{k2}A_{il}\D_{l}v_ka_{i2}}{|\D_1\boe|^2}\right\|_{L^2_t(\MH^{m})}^2+\|\Div_Av\|_{L^2_t(\MH^{m})}^2\right)+\left\|\frac{\tilde{\rho}_0J}{|\D_1\boe|^2}\right\|_{L^2_t(\MH^{m})}^2\nonumber\\
	&\ls \vep^2\left(\|\bA\DD\bv\|_{L^2_t(\MH^{m})}^2+\|\DD\bv\|_{L^\infty_{t,x}}\left(\|\D_1\boe\|_{L^2_t(\MH^{m})}+\left\|\frac{1}{|\D_1\boe|^2}\right\|_{L^2_t(\MH^{m})}^2\right)\right),\nonumber\\
	&\quad+\|J\|_{L^2_t(\MH^{m})}+\left\|\frac{1}{|\D_1\boe|^2}\right\|_{L^2_t(\MH^{m})}^2+\|\tilde{\rho}_0\|_{L^2_t(\MH^{m})}\\
	&\ls P\left(\|\D_1\boe\|_{L^2_t(\MH^{m})}^2\right)\left(\vep^2\|\DD\bv\|_{L^\infty}+1\right)+\vep^2\|\DD\bv\|_{L^2_t(\MH^{m})}\nonumber\\
	&\ls P\left(\|\D_1\boe\|_{L^2_t(\MH^{m})}^2\right)\left(\vep^2\|\DD\bv\|_{L^2_t(\MH^m)}+1\right)\ls P\left(\sup_{t\in[0,T]}\mathfrak{E}(t)\right).\nonumber
\end{align}
where we have used 
\begin{equation*}
	\||\D_1\boe|^{-2}\|_{L^2_t(\MH^{m})}^2\ls \|\D_1\boe\|_{L^2_t(\MH^{m})}^2+P\left(\|\D_1\boe\|_{L^2_{t,x}}^2\right)
\end{equation*}
Then, it follows from \eqref{L-infL-2} that
\begin{equation}\label{est-B-2}
	\|\mathfrak{B}\|_{L^\infty_t(\MH^{m-1})}\ls M_0+T\|\D_t\mathfrak{B}\|_{L^2_t(\MH^{m-1})}^2\ls M_0+TP\left(\sup_{t\in[0,T]}\mathfrak{B}(t)\right).
\end{equation}
As a consequence, the terms $R_b^i$, $i=1,\cdots,6$ can be bounded as follows. For $R_b^1$, direct calculations yield that
\begin{align*}
R_b^1&=-\frac{d}{dt}\int_\Gamma\mathfrak{B}\bar{\D}^\alpha a_{i2}\bar{\D}^\alpha\eta_i+\int_\Gamma\D_t\mathfrak{B}\bar{\D}^\alpha a_{i2}\bar{\D}^\alpha\eta_i+\int_\Gamma\mathfrak{B}\bar{\D}^\alpha \D_ta_{i2}\bar{\D}^\alpha\eta_i\\
&=-\frac{d}{dt}\int_\Gamma\mathfrak{B}\bar{\D}^\alpha a_{i2}\bar{\D}^\alpha\eta_i+\int_\Gamma\D_t\mathfrak{B}\bar{\D}^\alpha a_{i2}\bar{\D}^\alpha\eta_i+\int_\Gamma\D_1\mathfrak{B}(\bar{\D}^\alpha v_2\bar{\D}^\alpha \eta_1-\bar{\D}^\alpha v_1\bar{\D}^\alpha \eta_2)\\
&\quad+\underbrace{\int_\Gamma\mathfrak{B}(\bar{\D}^\alpha v_2\bar{\D}^\alpha \D_1\eta_1-\bar{\D}^\alpha v_1\bar{\D}^\alpha \D_1\eta_2)}_{-R_b^1}
\end{align*}
so that
\begin{align}\label{Rb1}
R_b^1=\underbrace{-\frac{1}{2}\frac{d}{dt}\int_\Gamma\mathfrak{B}\bar{\D}^\alpha a_{i2}\bar{\D}^\alpha\eta_i}_{R_b^{1,1}}+\underbrace{\frac{1}{2}\int_\Gamma\D_t\mathfrak{B}\bar{\D}^\alpha a_{i2}\bar{\D}^\alpha\eta_i+\frac{1}{2}\int_\Gamma\D_1\mathfrak{B}(\bar{\D}^\alpha v_2\bar{\D}^\alpha \eta_1-\bar{\D}^\alpha v_1\bar{\D}^\alpha \eta_2)}_{R_{b}^{1,2}}.
\end{align}
To estimate $R_b^{1,1}$, we use the following identities
\begin{equation}\label{decomp}
	a_{i2}a_{j2}+\D_1\eta_i\D_1\eta_j=|\D_1\boe|^2\delta_{ij},\quad\bar{\D}^\alpha a_{i2}a_{i2}=\bar{\D}^\alpha\D_1\eta_i\D_1\eta_i,\quad\bar{\D}^\alpha a_{i2}\D_1\eta_i=-\bar{\D}^\alpha\D_1\eta_ia_{i2}.
\end{equation}
to write
\begin{align*}
		R_b^{1,1}&=-\frac{1}{2}\frac{d}{dt}\int_\Gamma\mathfrak{B}\frac{1}{|\D_1\boe|^2}\bar{\D}^\alpha a_{i2}(a_{i2}a_{j2}+\D_1\eta_i\D_1\eta_j)\bar{\D}^\alpha\eta_j\\
		&=-\frac{1}{2}\frac{d}{dt}\int_\Gamma\mathfrak{B}\frac{\bar{\D}^\alpha\D_1\eta_i\D_1\eta_i}{|\D_1\boe|^2}a_{j2}\bar{\D}^\alpha\eta_j+\frac{1}{2}\frac{d}{dt}\int_\Gamma\mathfrak{B}\frac{\bar{\D}^\alpha\D_1\eta_ia_{i2}}{|\D_1\boe|^2}\D_1\eta_j\bar{\D}^\alpha\eta_j\\
		&=\frac{d}{dt}\int_\Gamma\mathfrak{B}\frac{\bar{\D}^\alpha\D_1\eta_ia_{i2}}{|\D_1\boe|^2}\D_1\eta_j\bar{\D}^\alpha\eta_j+\frac{1}{2}\frac{d}{dt}\int_\Gamma\D_1\left(\frac{\mathfrak{B}\D_1\eta_ia_{j2}}{|\D_1\boe|^2}\right)\bar{\D}^\alpha\eta_i\bar{\D}^\alpha\eta_j
\end{align*}
Thus, by using \eqref{fra-emb}, \eqref{tre} and \eqref{est-B},  we have
\begin{align}\label{Rb11}
\left|\int_{0}^{t}R_b^{1,1}\right|&\ls\left|\int_\Gamma\mathfrak{B}\frac{\bar{\D}^\alpha\D_1\eta_ia_{i2}}{|\D_1\boe|^2}\D_1\eta_j\bar{\D}^\alpha\eta_j\Bigg|_0^t\right|+\left|\int_\Gamma\D_1\left(\frac{\mathfrak{B}\D_1\eta_ia_{j2}}{|\D_1\boe|^2}\right)\bar{\D}^\alpha\eta_i\bar{\D}^\alpha\eta_j\Bigg|_0^t\right|\nonumber\\
&\ls M_0+|\bar{\D}^\alpha\D_1\eta_ia_{i2}|_0|\bar{\D}^\alpha\boe|_0|\mathfrak{B}|_1+|\D_1\mathfrak{B}|_{L^4}|\bar{\D}^\alpha\boe|_{L^4}|\bar{\D}^\alpha\boe|_0+|\mathfrak{B}|_1|\D_1^2\boe|_{L^4}|\bar{\D}^\alpha\boe|_{L^4}|\bar{\D}^\alpha\boe|_0\nonumber\\
&\ls M_0+\delta(|\bar{\D}^\alpha\D_1\eta_ia_{i2}|_0^2+\|\bar{\D}^\alpha\DD\boe\|_0^2)+C_\delta(\|\boe\|_{\MH^m}^2+\|\mathfrak{B}\|_{H^2}^2)^2\\
&\ls M_0+\delta\sup_{t\in[0,T]}\mathfrak{E}(t)+TP\left(\sup_{t\in[0,T]}\mathfrak{E}(t)\right).\nonumber
\end{align}
For $R_b^{1,2}$, since $\alpha_1\geq 1$, it follows from the dual estimate, \eqref{co123} and \eqref{est-B} that
	\begin{align}\label{Rb12}
		\left|\int_{0}^{t}R_b^{1,2}\right|&\ls\int_{0}^{t}|\bar{\D}\mathfrak{B}\bar{\D}^\alpha\boe|_{\frac{1}{2}}|\bar{\D}^\alpha \bar{\D}\boe|_{-\frac{1}{2}}\ls\int_{0}^{t}|\bar{\D}\mathfrak{B}|_1|\bar{\D}^\alpha\boe|_\frac{1}{2}|\bar{\D}^m\boe|_\frac{1}{2}\\
		&\ls \int_{0}^{t}(\|\mathfrak{B}\|_{\MH^2}+\|\DD\mathfrak{B}\|_{\MH^2}^\frac{1}{2}\|\mathfrak{B}\|_{\MH^2}^\frac{1}{2})\|\bar{\D}^\alpha\boe\|_1\|\bar{\D}^m\boe\|_1\nonumber\\
		&\ls (\|\mathfrak{B}\|_{L^2_t(\MH^2)}+\|\DD\mathfrak{B}\|_{L^2_t(\MH^2)})\|\bar{\D}^\alpha\boe\|_{L^4_t(H^1)}\|\bar{\D}^m\boe\|_{L^4_t(H^1)}\nonumber\\
		&\ls T^\frac{1}{4}P\left(\sup_{t\in[0,T]}\mathfrak{E}(t)\right).\nonumber
	\end{align}
which combining with \eqref{Rb11} implies that
\begin{equation}\label{Rb1-est}
	\left|\int_{0}^{t}R_b^1\right|\ls M_0+\delta\sup_{t\in[0,T]}\mathfrak{E}(t)+T^\frac{1}{4}P(\sup_{t\in[0,T]}\mathfrak{E}(t)).
\end{equation}
For $R_b^2$, it follows from duality and trace estimate that
\begin{align}\label{Rb2-est}
	\left|\int_{0}^{t}R_b^2\right|&\ls \int_{0}^{t}|[\bar{\D}^\alpha,\mathfrak{B},a_{i2}]|_\frac{1}{2}|\bar{\D}^\alpha v_i|_{-\frac{1}{2}}\nonumber\\
	&\ls \|[\bar{\D}^\alpha,\mathfrak{B},a_{i2}]\|_{L^2_t(H^1)}\|\bar{\D}^m\boe\|_{L^4_t(H^1)}t^\frac{1}{4}\ls T^\frac{1}{4}P\left(\sup_{t\in[0,T]}\mathfrak{E}(t)\right).
\end{align}
where we have used the following estimate, which is obtained by using \eqref{pro-est} and \eqref{est-B}
\begin{align*}
	\|\DD[\bar{\D}^\alpha,\mathfrak{B},a_{i2}]\|_{L^2_t(L^2)}^2&\ls \sum\limits_{1\leq|\beta|\leq m-1,|\beta|+|\nu|=m}\int_{0}^{t}\|\bar{\D}^\beta\DD\mathfrak{B}\bar{\D}^\nu \D_1\boe\|_0^2+\|\bar{\D}^\beta\mathfrak{B}\bar{\D}^\nu\DD\D_1\boe\|_0^2\nonumber\\
	&\ls \|\DD\mathfrak{B}\|_{L^\infty_{t,x}}\|\DD\D_1\boe\|_{L^2_t(\MH^{m-1})}+\|\DD\D_1\boe\|_{L^\infty_{t,x}}\|\DD\mathfrak{B}\|_{L^2_t(\MH^{m-1})}\\
	&\ls \|\DD\mathfrak{B}\|_{L^2_t(\MH^{m-1})}\|\DD\D_1\boe\|_{L^2_t(\MH^{m-1})}\ls P\left(\sup_{t\in[0,T]}\mathfrak{E}(t)\right).\nonumber
\end{align*} 
For $R_b^3$, noticing $\bar{\D}^\alpha a_{k2}\D_1^2\eta_k=-\bar{\D}^\alpha \D_1\eta_k\D_1 a_{k2}$, we can write it as 
\begin{align*}
R_b^3&=\sigma\int_\Gamma\frac{1}{|\D_1\boe|^3}\bar{\D}^\alpha \D_1\eta_k\D_1 a_{k2}a_{i2}\bar{\D}^\alpha v_i+\sigma\int_\Gamma\left[\bar{\D}^\alpha, \frac{1}{|\D_1\boe|^3}\right]\D_1\eta_k\D_1a_{k2}a_{i2}\bar{\D}^\alpha v_i\\
&=\underbrace{\sigma\frac{d}{dt}\int_\Gamma\frac{1}{|\D_1\boe|^3}\bar{\D}^\alpha \D_1\eta_k\D_1a_{k2}a_{i2}\bar{\D}^\alpha \eta_i}_{R_b^{3,1}}+\underbrace{\sigma\int_\Gamma\frac{1}{|\D_1\boe|^3}\bar{\D}^\alpha v_k\D_1a_{k2}a_{i2}\bar{\D}^\alpha \D_1\eta_i}_{R_b^{3,2}}\\
&\quad\underbrace{+\sigma\int_\Gamma\left(\bar{\D}^\alpha v_k\D_1\left(\frac{1}{|\D_1\boe|^3}\D_1a_{k2}a_{i2}\right)-\bar{\D}^\alpha \D_1\eta_k\D_t\left(\frac{1}{|\D_1\boe|^3}\D_1a_{k2}a_{i2}\right)\right)\bar{\D}^\alpha \eta_i}_{R_b^{3,3}}\\
&\quad+\underbrace{\sigma\int_\Gamma\left[\bar{\D}^\alpha, \frac{1}{|\D_1\boe|^3}\right]\D_1\eta_k\D_1a_{k2}a_{i2}\bar{\D}^\alpha v_i}_{R_b^{3,4}}.
\end{align*}
Then, utilizing integration by parts in $x_1$, \eqref{fra-emb}, \eqref{tre} and \eqref{L-infL-2}, we have
\begin{align}\label{Rb31}
	\left|\int_{0}^{t}R_b^{3,1}\right|&\ls\left|\sigma\int_\Gamma\frac{1}{|\D_1\boe|^3}\bar{\D}^\alpha \D_1\eta_k\D_1a_{k2}a_{i2}\bar{\D}^\alpha \eta_i\Bigg|_0^t\right|=\left|\sigma\int_\Gamma\bar{\D}^\alpha \eta_k\D_1\left(\frac{1}{|\D_1\boe|^3}\D_1a_{k2}a_{i2}\bar{\D}^\alpha \eta_i\right)\Bigg|_0^t\right|\nonumber\\
	&\ls M_0+\sigma|\bar{\D}^\alpha\eta_k|_0\left|\frac{\D_1a_{k2}}{|\D_1\boe|^3}\right|_{L^\infty}|\bar{\D}^\alpha\D_1\eta_ia_{i2}|_0(t)\nonumber\\
	&\quad+\sigma|\bar{\D}^\alpha\eta_k|_0\left|\D_1\left(\frac{\D_1a_{k2}}{|\D_1\boe|^3}\right)\right|_0|\bar{\D}^\alpha\eta_ia_{i2}|_1(t)+\sigma|\bar{\D}^\alpha\eta_k|_0\left|\frac{\D_1a_{k2}}{|\D_1\boe|^3}\D_1a_{i2}\right|_{L^4}|\bar{\D}^\alpha\eta_i|_{L^4}(t)\nonumber\\
	&\ls M_0+\sigma\left(|\D_1^2\boe|_1+|\D_1^2\boe
	|_{L^4}^2\right)|\bar{\D}^\alpha\boe|_0|\bar{\D}^\alpha\D_1\eta_ia_{i2}|_0(t)\\
	&\quad+\sigma(|\D_1^3\boe|_0+|\D_1^2\boe
	|_{L^4}^2)|\bar{\D}^\alpha\boe|_0(|\bar{\D}^\alpha\boe|_{L^4}|\D_1^2\boe|_{L^4}+|\bar{\D}^\alpha\boe|_0)\nonumber\\
	&\quad+\sigma|\bar{\D}^\alpha\eta_k|_0\left|\D_1a_{k2}\D_1a_{i2}\right|_{\hal}|\bar{\D}^\alpha\eta_k|_\hal\nonumber\\
	&\ls M_0+\sigma\|\D_1^3\DD\boe\|_0\|\bar{\D}^\alpha\DD\boe\|^\frac{1}{2}\|\bar{\D}^\alpha\boe\|_0^\frac{1}{2}|\bar{\D}^\alpha\D_1\eta_ia_{i2}|_0\nonumber\\
	&\quad+\sigma\|\D_1^3\DD\boe\|_0\|\bar{\D}^\alpha\DD\boe\|^\frac{3}{2}\|\bar{\D}^\alpha\boe\|_0^\frac{1}{2}\|\D_1^2\boe\|_1+P(\|\boe\|_{\MH^m}^2)\nonumber\\
	&\ls M_0+\delta\sigma|\bar{\D}^\alpha\D_1\eta_ia_{i2}|_0^2+\delta\|\bar{\D}^\alpha\DD\eta\|_0^2+TP\left(\sup_{t\in[0,T]}\mathfrak{E}(t)\right).\nonumber
\end{align}	
From \eqref{BB}, one has
\begin{align*}
\sigma\D_1^2\eta_ka_{k2}=-|\D_1\boe|^3q+\tilde{\rho}_0|\D_1\boe|J+2\mu\vep a_{k2}A_{il}\D_lv_ka_{i2}|\D_1\boe|+\lambda\vep\Div_Av|\D_1\boe|^3.
\end{align*}
Then, for any $\beta$ with $|\beta|\leq m-1$
	\begin{align}\label{4.24}
		|\bar{\D}^\beta\D_1^2\eta_ka_{k2}|_{L^2_t(H^\frac{1}{2})}^2&\ls |[\bar{\D}^\beta, a_{k2}]\D_1^2\eta_k|_{L^2_t(H^\frac{1}{2})}^2+|\bar{\D}^\beta(\tilde{\rho}_0|\D_1\eta|J)|_{L^2_t(H^\frac{1}{2})}^2+|\bar{\D}^\beta(|\D_1\boe|^3q)|_{L^2_t(H^\frac{1}{2})}^2\nonumber\\
		&\quad+\vep^2|\bar{\D}^\beta(a_{k2}A_{il}\D_lv_ka_{i2}|\D_1\boe|)|_{L^2_t(H^\frac{1}{2})}^2+\vep^2|\bar{\D}^\beta(\Div_Av|\D_1\boe|^3)|_{L^2_t(H^\frac{1}{2})}^2\nonumber\\
		&\ls\|[\bar{\D}^\beta, a_{k2}]\D_1^2\eta_k\|_{L^2_t(H^1)}^2+\|\bar{\D}^\beta(\tilde{\rho}_0|\D_1\boe|J)\|_{L^2_t(H^1)}^2+\|\bar{\D}^\beta(|\D_1\boe|^3q)\|_{L^2_t(H^1)}^2\nonumber\\
		&\quad+\vep^2\|\bar{\D}^\beta(a_{k2}A_{il}\D_lv_ka_{i2}|\D_1\boe|)\|_{L^2_t(H^1)}^2+\vep^2\|\bar{\D}^\beta(\Div_Av|\D_1\boe|^3)\|_{L^2_t(H^1)}^2\nonumber\\
		&\ls \int_{0}^{t}\|\nabla\eta\|_{\MH^m}^2+\|\nabla q\|_{\MH^{m-1}}^2+\vep^2\|\DD\bv\|_{\MH^m}^2\ls P\left(\sup_{t\in[0,T]}\mathfrak{E}(t)\right).
	\end{align}
which combing with \eqref{co123} yields that,  for $\alpha_1\geq 1$, one has
	\begin{align}\label{Rb32}
		\left|\int_{0}^{t}R_b^{3,2}\right|&\ls\sigma\int_{0}^{t}|\bar{\D}^\alpha \bv|_{-\frac{1}{2}}|\D_1^2\boe|_1|a_{i2}\bar{\D}^\alpha \D_1\eta_i|_{\frac{1}{2}}\nonumber\\
		&\ls\sigma|a_{i2}\bar{\D}^\alpha \D_1\eta_i|_{L^2_t(H^\frac{1}{2})}|\bar{\D}^{m-1}\bv|_{L^4_t(H^\frac{1}{2})}|\D_1^2\boe|_{L^4_t(H^1)}\\
		&\ls T^\frac{1}{4} P\left(\sup_{t\in[0,T]}\mathfrak{E}(t)\right).\nonumber
	\end{align}
It follows from \eqref{co123} that
\begin{align}\label{Rb33}
\left|\int_{0}^{t}R_b^{3,3}\right|&\ls\int_{0}^{t}|\bar{\D}^\alpha\bar{\D}\boe|_{-\frac{1}{2}}|\D_1^2\bar{\D}\boe|_\frac{1}{2}|\bar{\D}^\alpha\eta_i a_{i2}|_1+\int_{0}^{t}|\bar{\D}^\alpha\bar{\D}\boe|_{-\frac{1}{2}}|\D_1^2\boe\D_1\bar{\D}\boe|_1|\bar{\D}^\alpha\boe|_\frac{1}{2}\nonumber\\
&\ls \int_{0}^{t}\|\bar{\D}^m\boe\|_1P(\|\boe\|_{\MH^m})(|\bar{\D}^\alpha\D_1\eta_i a_{i2}|_0+\|\bar{\D}^\alpha\boe\|_1)\ls T^\frac{1}{4}P\left(\sup_{t\in[0,T]}\mathfrak{E}(t)\right).
\end{align}
and
\begin{align}\label{Rb34}
\left|\int_{0}^{t}R_b^{3,4}\right|&\ls\int_{0}^{t}\left|\left[\bar{\D}^\alpha, \frac{1}{|\D_1\boe|^3}\right]\D_1\boe\right|_\frac{1}{2}|\D_1^2\boe \D_1\boe|_1|\bar{\D}^\alpha \bv|_{-\frac{1}{2}}\ls T^\frac{1}{4}P\left(\sup_{t\in[0,T]}\mathfrak{E}(t)\right).
\end{align}
Thus, \eqref{Rb31}-\eqref{Rb34} yields
\begin{equation}\label{Rb3-est}
	\left|\int_{0}^{t}R_b^3\right|\ls  M_0+\delta\sup_{t\in[0,T]}\mathfrak{E}(t)+T^\frac{1}{4}P\left(\sup_{t\in[0,T]}\mathfrak{E}(t)\right).
\end{equation}
For $R_{b4}$, we can obtain from \eqref{co123}, trace estimate and \eqref{co2} that
\begin{equation}\label{Rb4-est}
	\left|\int_{0}^{t}R_b^4\right|\ls\int_{0}^{t}\left|\left[\bar{\D}^\alpha,\D_1^2\eta_k,\frac{a_{k2}}{|\D_1\boe|^3}\right]\right|_\frac{1}{2}|a_{i2}|_1|\bar{\D}^\alpha v_i|_{-\frac{1}{2}}\ls T^\frac{1}{4}P\left(\sup_{t\in[0,T]}\mathfrak{E}(t)\right).
\end{equation}
The term $R_{b5}$ can be bounded by similar arguments in \eqref{Rb31} and \eqref{4.24} as follows 
\begin{align}\label{Rb5-est}
		\left|\int_{0}^{t}R_b^5\right|&\ls\int_{0}^{t}\sigma|\bar{\D}^\alpha\D_1\eta_k a_{k2}|_\frac{1}{2}|\D_1^2\boe|_1|\bar{\D}^\alpha \bv|_{-\frac{1}{2}}+\left|\int_{0}^{t}\int_\Gamma\sigma\bar{\D}^\alpha\D_1\eta_k\frac{\D_1a_{k2}}{|\D_1\boe|^3}a_{i2}\bar{\D}^\alpha v_i\right|\nonumber\\
		&\ls M_0+\delta\sup_{t\in[0,T]}\mathfrak{E}(t)+T^\frac{1}{4}P\left(\sup_{t\in[0,T]}\mathfrak{E}(t)\right).
\end{align}
For $R_{b6}$, we utilize \eqref{4.24}, \eqref{co123} and \eqref{tre} to yield that
\begin{align}\label{Rb6-est}
	\left|\int_{0}^{t}R_b^6\right|&\ls\left|\sigma\int_{0}^{t}\int_\Gamma\frac{\D_t a_{k2}}{|\D_1\boe|^3}\bar{\D}^\alpha\D_1\eta_k\bar{\D}^\alpha\D_1\eta_ia_{i2}\right|+\left|\sigma\int_{0}^{t}\int_\Gamma\frac{\D_1\eta_i\D_t\D_1\eta_i}{|\D_1\eta|^5}|\bar{\D}^\alpha\D_1\eta_ka_{k2}|^2\right|\nonumber\\
	&\ls \sigma\int_{0}^{t}|\bar{\D}^\alpha\D_1\boe|_{-\frac{1}{2}}\left|\frac{\D_t a_{k2}}{|\D_1\boe|^3}\right|_1|\bar{\D}^\alpha\D_1\eta_ia_{i2}|_\frac{1}{2}+\sigma|\bar{\D}^\alpha\D_1\eta_ia_{i2}|_{L^\infty_t(L^2_x)}\int_{0}^{t}|\D_t\D_1\boe|_{L^\infty}\\
	&\ls T^\frac{1}{2}P\left(\sup_{t\in[0,T]}\mathfrak{E}(t)\right).\nonumber
\end{align}
Therefore, plugging \eqref{II-est} into \eqref{tang-est}, integrating the resulted equation over $[0,t]$, and then substituting \eqref{R1-est}-\eqref{R4-est}, \eqref{Rq1-est}, \eqref{Rq2-est},\eqref{Rvep1},\eqref{Rvep2}, \eqref{Rb1-est},\eqref{Rb2-est},\eqref{Rb3-est}-\eqref{Rb6-est}, and using the estimate \eqref{kon} and Korn's inequality \eqref{Korn's ineq}, we can complete the proof of this lemma.
\end{proof}

\subsection{Fully temporal derivative estimates}
In view of argument in the proof of lemma \ref{Tan-est}, we used the duality argument taking advantage of the estimate $|\bar{\D}^\alpha\boe|_{-\frac{1}{2}}\ls |\bar{\D}^{m-1}\boe|_\frac{1}{2}$ if $|\alpha|=m$ and $\alpha_1\geq 1$, which is invalid for $|\D_t^m\boe|_{-\frac{1}{2}}$. Moreover, we are lack of the estimate of $|\D_t^m\D_1\eta_ka_{k2}|_{L^2_t(H^\frac{1}{2})}$, while the similar estimate \eqref{4.24} is the key to derive \eqref{Rb6-est} in the proof of lemma \ref{tang-est}. To overcome these difficulties, we introduce the following Alinhac's good unknows 
\begin{equation}\label{V-Q}
	\mathcal{V}_i^\vep=\D_t^mv_i^\vep-\D_t^m\eta_k^\vep A_{kl}^\vep\D_lv_i^\vep,\quad \mathcal{Q}^\vep=\D_t^mq^\vep-\D_t^m\eta_k^\vep A_{kl}^\vep\D_lq^\vep,
\end{equation}
which can be used to cancel the highest nonlinear term in the proof of fully temporal derivative estimates. In order to derive the equation of $\mathcal{V}_i$ and $\mathcal{Q}$, we take the fully temporal derivative of $a_{ij}^\vep\D_j q^\vep$ and $a_{ij}^\vep\D_j v_i^\vep$ and  get from \eqref{Geo-iden-1} that
\begin{align*}
	&\D_t^m(a_{ij}^\vep\D_j q^\vep)=a_{ij}^\vep\D_j \mathcal{Q}^\vep+[\D_t^m,a_{ij}^\vep,\D_jq^\vep]+\mathcal{C}_i(q^\vep),\\
	&\D_t^m(a_{ij}^\vep\D_j v_i^\vep)=a_{ij}^\vep\D_j\mathcal{V}_i^\vep+[\D_t^m,a_{ij}^\vep,\D_jv_i^\vep]+\mathcal{C}_i(v_i^\vep),
\end{align*}
where $\mathcal{C}_i(f)$ is given by
\begin{equation}\label{com-q}
	\mathcal{C}_i(f)=\D_t^m\eta^\vep_ka_{ij}^\vep\D_j(A_{kl}^\vep\D_lf)-J^\vep[\D_t^{m-1},A_{il}^\vep A_{kj}^\vep]\D_l\D_t\eta^\vep_k\D_jf+[\D_t^m,J^\vep]A^\vep_{ij}\D_jf.
\end{equation}
Then,  by applying $\D_t^m$ to $\eqref{FEL-final}_2$ and the equation   
\begin{equation*}
	a_{ij}\D_jv_i=J_t=-\frac{J^{\gamma+1}}{\gamma \tilde{\rho}_0^\gamma}\D_t q,
\end{equation*}
we have
\begin{align}\label{Eqn of good unknowns}
	&\tilde{\rho}_0^\vep\D_t\mathcal{V}_i^\vep+a_{ij}^\vep\D_j\mathcal{Q}^\vep-2\mu\vep\D_l\D_t^m((S_Av)^\vep_{ik}a^\vep_{kl})-\lambda\vep\D_j\D_t^m(a_{ij}^\vep(\Div_Av)^\vep)-\D_j(\tilde{\rho}_0^\vep\D_j\D_t^m\eta_i^\vep)\nonumber\\
	&=-\tilde{\rho}_0^\vep\D_t(\D_t^m\eta^\vep_kA^\vep_{kl}\D_lv^\vep_i)-[\D_t^m,a_{ij}^\vep,\D_jq^\vep]-\mathcal{C}_i(q^\vep),
\end{align}
and
\begin{equation}\label{com-v}
	a_{ij}^\vep\D_j\mathcal{V}^\vep_i=-\frac{(J^\vep)^{\gamma+1}}{\gamma (\tilde{\rho}_0^\vep)^\gamma}\D_t\mathcal{Q}^\vep+\frac{(J^\vep)^{\gamma+1}}{\gamma (\tilde{\rho}_0^\vep)^\gamma}\D_t(\D_t^m\eta_k^\vep A_{kl}^\vep\D_l q^\vep)-\left[\D_t^m,\frac{(J^\vep)^{\gamma+1}}{\gamma (\tilde{\rho}_0^\vep)^\gamma}\right]\D_tq^\vep-[\D_t^m,a_{ij}^\vep,\D_jv_i^\vep]-\mathcal{C}_i(v_i^\vep).
\end{equation}

Before giving the full temporal derivative estimates, we first derive the following estimate of the commutator $\mathcal{C}_i(q^\vep)$ and $\mathcal{C}_i(v_i^\vep)$.
\begin{lemma}
	For any $m\geq 4$ and $t\in [0,T^\vep]$, it holds that
	\begin{equation}\label{commutator-est}
		\int_{0}^{t}\|\mathcal{C}_i(q^\vep)\|_0^2+\|\mathcal{C}_i(v_i^\vep)\|_0^2ds\ls T^\vep P\left(\sup_{t\in[0,T^\vep]}\mathfrak{E}^\vep(t)\right).
	\end{equation}
\end{lemma}
\begin{pf}
	We only need to give the estimate of $\mathcal{C}_i(q)$ since the other is the same. By using \eqref{Piola}, \eqref{pro-est}, the Sobolev embedding and \eqref{A-priori-assum}, one has
		\begin{align*}
			&\|\D_t^m\eta_ka_{ij}\D_j(A_{kl}\D_l q)\|_{L^2_{t,x}}^2=\|\D_t^m\eta_k\D_j(a_{ij}A_{kl}\D_l q)\|_{L^2_{t,x}}^2\nonumber\\
			&\ls \|\D_t^2\eta_k\|_{L^\infty_{t,x}}^2\|a_{ij}A_{kl}\D_l q\|_{L^2_t(\MH^{m-1})}^2+\|\D_t^2\eta_k\|_{L^2_t(\MH^{m-1})}^2\|a_{ij}A_{kl}\D_l q\|_{L^\infty_{t,x}}^2\\
			&\ls\|\D_t^2\boe\|_{L^2_t(\MH^{m-1})}^2(\|\DD \boe\|_{L^\infty_{t,x}}^2\|\DD q\|_{L^2(\MH^{m-1})}^2+\|\DD q\|_{L^\infty_{t,x}}^2\|\DD \boe
			\|_{L^2(\MH^{m-1})}^2)\nonumber\\
			&\ls TP\left(\sup_{t\in[0,T]}\mathfrak{E}(t)\right).\nonumber
		\end{align*}
Similarly, it follows from \eqref{A-priori-assum} and \eqref{co1} that
	\begin{align*}
		&\|J[\D_t^{m-1},A_{il}A_{kj}]\D_l\D_t\eta_k\D_jq\|_{L^2_{t,x}}^2\ls\|\DD q\|_{L^\infty_{t,x}}^2\|[\D_t^{m-1},A_{il}A_{kj}]\D_l\D_t\eta_k\|_{L^2_{t,x}}^2\nonumber\\
		&\ls \|\DD q\|_{L^\infty_{t,x}}^2\|\D_l\D_t\eta_k\|_{L^\infty_{t,x}}^2\|\D_t(A_{il}A_{kl})\|_{L^2_t(\MH^{m-2})}^2+\|\DD q\|_{L^\infty_{t,x}}^2\|\D_t(A_{il}A_{kl})\|_{L^\infty}^2\|\D_l\D_t\eta_k\|_{L^2_t(\MH^{m-2})}^2\\
		&\ls TP\left(\sup_{t\in[0,T]}\mathfrak{E}(t)\right)\nonumber
	\end{align*}
and
	\begin{align*}
		\|[\D_t^m,J]A_{ij}\D_jq\|_{L^2_{t,x}}^2&\ls \|\D_j q\|_{L^\infty_{t,x}}^2\|[\D_t^m,J]A_{ij}\|_{L^2_{t,x}}^2\nonumber\\
		&\ls\|\D_j q\|_{L^\infty_{t,x}}^2(\|A_{ij}\|_{L^\infty_{t,x}}^2\|\D_tJ\|_{L^2(\MH^2)}^2+\|\D_tJ\|_{L^\infty}^{m-1}\|A_{ij}\|_{L^2_t(\MH^{m-1})}^2)\\
		&\ls TP\left(\sup_{t\in[0,T]}\mathfrak{E}(t)\right).\nonumber
	\end{align*}
Therefore, combing the above estimates, we have
\begin{equation*}
	\|\mathcal{C}_i(q)\|_{L^2_t(L^2)}\ls TP\left(\sup_{t\in[0,T]}\mathfrak{E}(t)\right).
\end{equation*}

\end{pf}

Now, we derive the following estimate of fully temporal derivatives.
\begin{lemma}\label{Full-time-derivative}
	For any $m\geq 4$, $t\in[0,T_\vep]$, it holds that
\begin{align}\label{full-time-est}
	&\int_{0}^{t}\|\D_t^m\bv^\vep\|_0^4+\|\D_t^mq^\vep\|_0^4+\|\DD\D_t^m\boe^\vep\|_0^4+|\D_t^m\D_1\boe^\vep\cdot \bn^\vep|_0^4+\vep^2\int_{0}^{t}\left(\int_{0}^{s}\|\DD\D_t^m\bv^\vep\|^2_0\right)^2\nonumber\\
	&\ls M_0+\delta\sup_{t\in[0,T_\vep]}\mathfrak{E}^\vep(t)+T_\vep^\frac{1}{4}P\left(\sup_{t\in[0,T_\vep]}\mathfrak{E}^\vep(t)\right).
\end{align}
\end{lemma}
\begin{pf}
		Multiplying \eqref{Eqn of good unknowns}by  $\mathcal{V}_i$ and integrating over $\Omega$ yield
	\begin{align}\label{Full-time-deri}
		&\frac{1}{2}\frac{d}{dt}\int_{\Omega}\tilde{\rho}_0|\mathcal{V}|^2dx+\int_{\Omega}a_{ij}\D_j\mathcal{Q}\mathcal{V}_idx-\int_\Omega\D_j(\tilde{\rho}_0\D_j\D_t^3\eta_i)\mathcal{V}_i\nonumber\\
		&\quad-2\mu\vep\int_\Omega\D_l\D_t^m(S_A(v)_{ik}a_{kl})\mathcal{V}_i-\lambda\vep\int_\Omega\D_j\D_t^m(a_{ij}\Div_Av)\mathcal{V}_i\\
		&=\underbrace{-\int_\Omega\tilde{\rho}_0\D_t(\D_t^m\eta_kA_{kl}\D_lv_i)\mathcal{V}_idx}_{\mathcal{R}_\eta^1}\underbrace{-\int_\Omega[\D_t^m,a_{ij},\D_jq]\mathcal{V}_idx}_{\mathcal{R}_\eta^2}\underbrace{-\int_\Omega\mathcal{C}_i(q)\mathcal{V}_idx}_{\mathcal{R}_c^1}.\nonumber		
	\end{align}
	By integration by parts and \eqref{com-v}, one has
	\begin{align}\label{Ftd-1}
		&\int_{\Omega}a_{ij}\D_j\mathcal{Q}\mathcal{V}_idx-\int_\Omega\D_j(\tilde{\rho}_0\D_j\D_t^m\eta_i)\mathcal{V}_i-2\mu\vep\int_\Omega\D_l\D_t^m(S_A(v)_{ik}a_{kl})\mathcal{V}_i-\lambda\vep\int_\Omega\D_j\D_t^m(a_{ij}\Div_Av)\mathcal{V}_i\nonumber\\
		&=\int_\Gamma(a_{i2}\mathcal{Q}-\tilde{\rho}_0\D_t^m\D_2\eta_i-2\mu\vep\D_t^m(S_A(v)_{ik}a_{k2})-\lambda\vep\D_t^m(a_{i2}\Div_Av))\mathcal{V}_i-\int_\Omega a_{ij}\mathcal{Q}\D_j\mathcal{V}_idx\nonumber\\
		&\quad+\int_{\Omega}\tilde{\rho}_0\D_j\D_t^m\eta_i\D_j\mathcal{V}_idx+2\mu\vep\int_\Omega J|\boldsymbol{S}_{\bA}(\D_t^m \bv)|^2dx+\lambda\vep\int_\Omega J(\Div_A(\D_t^m v))^2dx\nonumber\\
		&\quad+\underbrace{2\mu\vep\int_\Omega[\D_t^m,a_{kl}](S_Av)_{ik}\D_l\mathcal{V}_i\,dx}_{-\mathcal{R}_\vep^1}+\underbrace{\lambda\vep\int_\Omega[\D_t^m,a_{ij}](\Div_Av)\D_j\mathcal{V}_idx}_{-\mathcal{R}_\vep^2}\nonumber\\
		&\quad+\underbrace{2\mu\vep\int_\Omega[\D_t^m,A_{ij}]\D_jv_kJS_A(\mathcal{V})_{ik}dx}_{-\mathcal{R}_\vep^3}+\underbrace{\lambda\vep\int_\Omega [\D_t^m,A_{ij}]\D_j v_iJ\Div_A(\mathcal{V})}_{-\mathcal{R}_\vep^4}\\
		&\quad+\underbrace{2\mu\vep\int_\Omega S_A(\D_t^mv)_{ik}a_{kl}\D_l(\D_t^m\eta_rA_{rs}\D_sv_i)dx}_{-\mathcal{R}_\vep^5}+\underbrace{\lambda\vep\int_\Omega\Div_A(\D_t^mv)a_{ij}\D_j(\D_t^m\eta_rA_{rs}\D_sv_i)dx}_{-\mathcal{R}_\vep^6}\nonumber\\
		&=\int_\Gamma(a_{i2}\mathcal{Q}-\tilde{\rho}_0\D_t^m\D_2\eta_i-2\mu\vep\D_t^m(S_A(v)_{ik}a_{k2})-\lambda\vep\D_t^m(a_{i2}\Div_Av))\mathcal{V}_i-\int_\Omega a_{ij}\mathcal{Q}\D_j\mathcal{V}_idx\nonumber\\
		&\quad+\frac{1}{2}\frac{d}{dt}\int_\Omega\tilde{\rho}_0|\D^m_t\DD\eta|^2dx+2\mu\vep\int_\Omega J|\boldsymbol{S}_{\bA}(\D_t^m \bv)|^2dx+\lambda\vep\int_\Omega J(\Div_A(\D_t^m v))^2dx\nonumber\\
		&\quad-\underbrace{\int_\Omega\tilde{\rho}_0\D_t^m\D_j\eta_i\D_j(\D_t^m\eta_kA_{kl}\D_lv_i)dx}_{\mathcal{R}_\eta^3}-\sum_{i=1}^{6}\mathcal{R}_\vep^i.\nonumber
	\end{align}
It follows from \eqref{com-v} that
\begin{align}\label{Ftd-2}
	-\int_\Omega a_{ij}\mathcal{Q}\D_j\mathcal{V}_idx&=\int_\Omega\frac{J^{\gamma+1}}{\gamma\tilde{\rho}_0^\gamma}\D_t\mathcal{Q}\mathcal{Q}dx-\underbrace{\int_\Omega\frac{J^{\gamma+1}}{\gamma\tilde{\rho}_0^\gamma}\D_t(\D_t^m\eta_kA_{kl}\D_l q)\mathcal{Q}dx}_{\mathcal{R}_q^1}+\underbrace{\int_\Omega\left[\D_t^m, \frac{J^{\gamma+1}}{\gamma\tilde{\rho}_0^\gamma}\right]\D_tq\mathcal{Q}dx}_{-\mathcal{R}_q^2}\nonumber\\
	&\quad+\underbrace{\int_\Omega[\D_t^m,a_{ij},\D_j v_i]\mathcal{Q}dx}_{-\mathcal{R}_q^3}+\underbrace{\int_\Omega\mathcal{C}_i(v_i)\mathcal{Q}dx}_{-\mathcal{R}_c^2}\\
	&=\frac{1}{2}\frac{d}{dt}\int_\Omega\frac{J^{\gamma+1}}{\gamma\tilde{\rho}_0^\gamma}\mathcal{Q}^2-\underbrace{\frac{\gamma+1}{2\gamma}\int_\Omega\frac{J^\gamma\D_tJ}{\tilde{\rho}_0^\gamma}\mathcal{Q}^2}_{\mathcal{R}_q^4}-\mathcal{R}_q^1-\mathcal{R}_q^2-\mathcal{R}_q^3-\mathcal{R}_c^2.\nonumber
\end{align}
	Plugging \eqref{Ftd-1} and \eqref{Ftd-2} into \eqref{Full-time-deri} gives
	\begin{align}\label{FTD}
		&\frac{1}{2}\frac{d}{dt}\int_\Omega\left(\tilde{\rho}_0|\mathcal{V}|^2+\frac{J^{\gamma+1}}{\gamma\tilde{\rho}_0^\gamma}\mathcal{Q}^2+\tilde{\rho}_0|\D^m_t\DD\eta|^2\right)dx+2\mu\vep\int_\Omega J|\boldsymbol{S}_{\bA}(\D_t^m \bv)|^2dx+\lambda\vep\int_\Omega J(\Div_A(\D_t^m v))^2dx\nonumber\\
		&\quad+\underbrace{\int_\Gamma(a_{i2}\mathcal{Q}-\tilde{\rho}_0\D_t^m\D_2\eta_i-2\mu\vep\D_t^m(S_A(v)_{ik}a_{k2})-\lambda\vep\D_t^m(a_{i2}\Div_Av))\mathcal{V}_i}_{\mathcal{R}_b}\\
		&=\sum_{i=1}^{3}\mathcal{R}_\eta^i+\sum_{i=1}^{4}\mathcal{R}_q^i+\mathcal{R}_c^1+\mathcal{R}_c^2+\sum_{i=1}^{6}\mathcal{R}_\vep^i.\nonumber
	\end{align}
Similar to \eqref{II-est}, it follows from the third equation of \eqref{FEL-final} that
	\begin{align}\label{RB}
		\mathcal{R}_b&=-\sigma\int_\Gamma\D_t^m\left(\frac{\D_1^2\eta_ka_{k2}}{|\D_1\boe|^3}\right)a_{i2}\mathcal{V}_i-\int_\Gamma\mathfrak{B}\D_t^ma_{i2}\mathcal{V}_i-\int_\Gamma[\D_t^m,\mathfrak{B},a_{i2}]\mathcal{V}_i-\int_\Gamma\D_t^m\eta_kA_{kl}\D_lqa_{i2}\mathcal{V}_i\nonumber\\
		&=-\sigma\int_\Gamma\frac{\D_t^m\D_1^2\eta_ka_{k2}}{|\D_1\boe|^3}a_{i2}\mathcal{V}_i-\int_\Gamma\D_t^m\left(\frac{a_{k2}}{|\D_1\boe|^3}\right)\D_1^2\eta_ka_{i2}\mathcal{V}_i-\int_\Gamma\mathfrak{B}\D_t^ma_{i2}\mathcal{V}_i\nonumber\\
		&\quad-\int_\Gamma[\D_t^m,\mathfrak{B},a_{i2}]\mathcal{V}_i-\int_\Gamma \left[\D_t^m,\D_1^2\eta_k,\frac{a_{k2}}{|\D_1\boe|^3}\right]a_{i2}\mathcal{V}_i-\int_\Gamma\D_t^m\eta_kA_{kl}\D_lqa_{i2}\mathcal{V}_i\nonumber\\
		&=\sigma\int_\Gamma\frac{\D_t^m\D_1\eta_ka_{k2}}{|\D_1\boe|^3}a_{i2}\D_1\mathcal{V}_i-\int_\Gamma\mathfrak{B}\D_t^ma_{i2}\mathcal{V}_i-\int_\Gamma[\D_t^m,\mathfrak{B},a_{i2}]\mathcal{V}_i-\int_\Gamma \left[\D_t^m,\D_1^2\eta_k,\frac{a_{k2}}{|\D_1\boe|^3}\right]a_{i2}\mathcal{V}_i\nonumber\\
		&\quad+\int_\Gamma\D_t^m\D_1\eta_k\D_1\left(\frac{a_{k2}a_{i2}}{|\D_1\boe|^3}\right)\mathcal{V}_i-\int_\Gamma\D_t^m\left(\frac{a_{k2}}{|\D_1\boe|^3}\right)\D_1^2\eta_ka_{i2}\mathcal{V}_i-\int_\Gamma\D_t^m\eta_kA_{kl}\D_lqa_{i2}\mathcal{V}_i\nonumber\\
		&=\sigma\int_\Gamma\frac{\D_t^m\D_1\eta_ka_{k2}}{|\D_1\boe|^3}\D_t(a_{i2}\D_1\D_t^m\eta_i)-\sigma\int_\Gamma\frac{\D_t^m\D_1\eta_ka_{k2}}{|\D_1\boe|^3}\left(\D_ta_{i2}\D_1\D_t^m\eta_i+a_{i2}\D_t^m\D_1\eta_kA_{kl}\D_lv_i\right)\nonumber\\
		&\quad-\int_\Gamma\mathfrak{B}\D_t^ma_{i2}\mathcal{V}_i-\int_\Gamma[\D_t^m,\mathfrak{B},a_{i2}]\mathcal{V}_i-\int_\Gamma \left[\D_t^m,\D_1^2\eta_k,\frac{a_{k2}}{|\D_1\boe|^3}\right]a_{i2}\mathcal{V}_i\nonumber\\
		&\quad-\int_\Gamma\left(\D_t^m\left(\frac{a_{k2}}{|\D_1\boe|^3}\right)\D_1^2\eta_ka_{i2}-\D_t^m\D_1\eta_k\D_1\left(\frac{a_{k2}a_{i2}}{|\D_1\boe|^3}\right)\right)\mathcal{V}_i\\
		&\quad-\sigma\int_\Gamma\frac{\D_t^m\D_1\eta_ka_{k2}}{|\D_1\boe|^3}a_{i2}\D_t^m\eta_k\D_1(A_{kl}\D_lv_i)-\int_\Gamma\D_t^m\eta_kA_{kl}\D_lqa_{i2}\mathcal{V}_i\nonumber\\	
		&=\frac{\sigma}{2}\frac{d}{dt}\int_\Gamma\frac{|\D_t^m\D_1\eta_ka_{k2}|^2}{|\D_1\boe|^3}-\int_\Gamma\mathfrak{B}\D_t^ma_{i2}\mathcal{V}_i-\int_\Gamma[\D_t^m,\mathfrak{B},a_{i2}]\mathcal{V}_i\nonumber\\
		&\quad-\int_\Gamma \left[\D_t^m,\D_1^2\eta_k,\frac{a_{k2}}{|\D_1\boe|^3}\right]a_{i2}\mathcal{V}_i-\int_\Gamma\left(\frac{1}{2}\D_t\left(\frac{1}{|\D_1\boe|^3}\right)+\frac{\D_tJ}{J}\right)|\D_t^m\D_1\eta_ka_{k2}|^2\nonumber\\
		&\quad-\int_\Gamma\left(\D_t^m\left(\frac{a_{k2}}{|\D_1\boe|^3}\right)\D_1^2\eta_ka_{i2}-\D_t^m\D_1\eta_k\D_1\left(\frac{a_{k2}a_{i2}}{|\D_1\boe|^3}\right)\right)\mathcal{V}_i\nonumber\\
		&\quad-\sigma\int_\Gamma\frac{\D_t^m\D_1\eta_ka_{k2}}{|\D_1\boe|^3}a_{i2}\D_t^m\eta_k\D_1(A_{kl}\D_lv_i)-\int_\Gamma\D_t^m\eta_kA_{kl}\D_lqa_{i2}\mathcal{V}_i\nonumber\\
		&=\frac{\sigma}{2}\frac{d}{dt}\int_\Gamma\frac{|\D_t^m\D_1\eta_ka_{k2}|^2}{|\D_1\boe|^3}-\sum_{i=1}^{7}\mathcal{R}_b^i,\nonumber
	\end{align}
where we have used the fact
\begin{equation*}
	\D_ta_{i2}\D_1\D_t^m\eta_i+a_{i2}\D_t^m\D_1\eta_kA_{kl}\D_lv_i=\D_tJJ^{-1}a_{i2}\D_1\D_t^m\eta_i.
\end{equation*}
Substituting \eqref{RB} into \eqref{FTD} gives
\begin{align}\label{Ful-tim-est}
		&\frac{1}{2}\frac{d}{dt}\int_\Omega\left(\tilde{\rho}_0|\mathcal{V}|^2+\frac{J^{\gamma+1}}{\gamma\tilde{\rho}_0^\gamma}\mathcal{Q}^2+\tilde{\rho}_0|\D^m_t\DD\eta|^2\right)dx+\frac{\sigma}{2}\frac{d}{dt}\int_\Gamma\frac{|\D_t^m\D_1\eta_ka_{k2}|^2}{|\D_1\boe|^3}\nonumber\\
		&+2\mu\vep\int_\Omega J|\boldsymbol{S}_{\bA}(\D_t^m \bv)|^2dx+\lambda\vep\int_\Omega J(\Div_A(\D_t^m v))^2dx\\
		&=\sum_{i=1}^{7}\mathcal{R}_b^i+\sum_{i=1}^{3}\mathcal{R}_\eta^i+\sum_{i=1}^{4}\mathcal{R}_q^i+\mathcal{R}_c^1+\mathcal{R}_c^2+\sum_{i=1}^{6}\mathcal{R}_\vep^i.\nonumber
\end{align}
Next, we estimate boundary terms $\mathcal{R}_b^i$ on the right hand side of \eqref{Ful-tim-est} one by one. For $\mathcal{R}_b^1$, similarly to \eqref{Rb1}, one has
\begin{align}
	\mathcal{R}_b^1&=\underbrace{\hal\frac{d}{dt}\int_\Gamma\mathfrak{B}\D_t^ma_{i2}\D_t^m\eta_i}_{\mathcal{R}_b^{1,1}}\underbrace{-\hal\int_\Gamma\D_t\mathfrak{B}\D_t^ma_{i2}\D_t^m\eta_i-\int_\Gamma\mathfrak{B}\D_t^ma_{i2}\D_t^m\eta_kA_{kl}\D_lv_i}_{\mathcal{R}_b^{1,2}}\nonumber\\
	&\quad-\underbrace{\hal\int_\Gamma\D_1\mathfrak{B}(\D_t^mv_2\D_t^m\eta_1-\D_t^mv_1\D_t^m\eta_2)}_{\mathcal{R}_b^{1,3}}.
\end{align}
It follows from the same argument as \eqref{Rb11} that
\begin{equation}\label{RB11}
	\left|\int_{0}^{t}\mathcal{R}_b^{1,1}\right|\ls M_0+\delta(|\D_t^m\D_1\eta_ia_{i2}(t)|_0^2+\|\D_t^m\DD\boe(t)\|_0^2)+TP\left(\sup_{t\in[0,T]}\mathfrak{E}(t)\right).
\end{equation}
By using \eqref{co123}, \eqref{BB} and \eqref{est-B}, one has
	\begin{align}\label{RB12}
		\left|\int_{0}^{t}\mathcal{R}_b^{1,2}\right|&\ls\int_{0}^{t}|\D_t^m\D_1\boe|_{-\frac{1}{2}}|\D_t^m\boe|_\frac{1}{2}\left(|\D_t\mathfrak{B}|_1+|\mathfrak{B}A_{kl}\D_lv_i|_1\right)\nonumber\\
		&\ls \|\D_t^m\boe\|_{L^4_t(H^1)}^2\left(|\D_t\mathfrak{B}|_{L^2_t(H^1)}+|\mathfrak{B}A_{kl}\D_lv_i|_{L^2_t(H^1)}\right)\\
		&\ls P(\sup_{t\in[0,T]}\mathfrak{E}(t)).\nonumber 
	\end{align}
Since that $\mathcal{R}_b^{1,3}$ involves only time derivatives of highest order, we can not use the duality argument as \eqref{Rb12} and \eqref{RB12}. Instead, we use the following decomposition,
\begin{equation}\label{decomp-2}
	\frac{-a_{i2}\D_1\eta_j+\D_1\eta_ia_{j2}}{|\D_1\boe|^2}=\left(\begin{array}{ll}0&1\\
		-1&0
	\end{array}\right),
\end{equation}
which is different to \eqref{decomp}, and use integration by parts to write the troublesome term into an integral in $\Omega$. That is,
\begin{align}
	\mathcal{R}_b^{1,3}&=\frac{1}{2}\int_\Gamma\frac{\D_1\mathfrak{B}}{|\D_1\boe|^2}(-a_{i2}\D_1\eta_j+\D_1\eta_ia_{j2})\D_t^mv_i\D_t^m\eta_j\nonumber\\
	&=-\int_\Gamma\frac{\D_1\mathfrak{B}}{|\D_1\boe|^2}\D_1\eta_j\D_t^m\eta_ja_{i2}\D_t^mv_i+\hal\frac{d}{dt}\int_\Gamma\frac{\D_1\mathfrak{B}}{|\D_1\boe|^2}\D_1\eta_i\D_t^m\eta_ia_{j2}\D_t^m\eta_j\nonumber\\
	&\quad-\hal\int_\Gamma\D_t\left(\frac{\D_1\mathfrak{B}}{|\D_1\boe|^2}\D_1\eta_ia_{j2}\right)\D_t^m\eta_i\D_t^m\eta_j\nonumber\\
	&=-\int_\Omega\frac{\D_1\mathfrak{B}}{|\D_1\boe|^2}\D_1\eta_j\D_t^m\eta_ja_{ik}\D_t^m\D_kv_i-\int_\Omega\D_k\left(\frac{\D_1\mathfrak{B}}{|\D_1\boe|^2}\D_1\eta_j\D_t^m\eta_j\right)a_{ik}\D_t^mv_i\nonumber\\
	&\quad+\hal\frac{d}{dt}\int_\Gamma\frac{\D_1\mathfrak{B}}{|\D_1\boe|^2}\D_1\eta_i\D_t^m\eta_ia_{j2}\D_t^m\eta_j-\hal\int_\Gamma\D_t\left(\frac{\D_1\mathfrak{B}}{|\D_1\boe|^2}\D_1\eta_ia_{j2}\right)\D_t^m\eta_i\D_t^m\eta_j\nonumber\\
	&=-\frac{d}{dt}\int_\Omega\frac{\D_1\mathfrak{B}}{|\D_1\boe|^2}\D_1\eta_j\D_t^m\eta_ja_{ik}\D_t^m\D_k\eta_i\,dx+\int_\Omega\D_t\left(\frac{\D_1\mathfrak{B}}{|\D_1\boe|^2}\D_1\eta_j\D_t^m\eta_ja_{ik}\right)\D_t^m\D_k\eta_i\,dx\\
	&\quad-\int_\Omega\D_k\left(\frac{\D_1\mathfrak{B}}{|\D_1\boe|^2}\D_1\eta_j\D_t^m\eta_j\right)a_{ik}\D_t^mv_i\,dx+\hal\frac{d}{dt}\int_\Gamma\frac{\D_1\mathfrak{B}}{|\D_1\boe|^2}\D_1\eta_i\D_t^m\eta_ia_{j2}\D_t^m\eta_j\nonumber\\
	&\quad-\hal\int_\Gamma\D_t\left(\frac{\D_1\mathfrak{B}}{|\D_1\boe|^2}\D_1\eta_ia_{j2}\right)\D_t^m\eta_i\D_t^m\eta_j=:\sum_{i=1}^{5}\mathcal{R}_b^{1,3,i}\nonumber
\end{align}
By using H\"{o}lder's inequality and \eqref{est-B-2}, one has
\begin{align*}
	\left|\int_{0}^{t}\mathcal{R}_b^{1,3,1}\right|&=\left|\int_\Omega\frac{\D_1\mathfrak{B}}{|\D_1\boe|^2}\D_1\eta_j\D_t^m\eta_ja_{ik}\D_t^m\D_k\eta_i\bigg|_0^t\right|\ls M_0+\|\D_t^m\DD\boe(t)\|_0\|\D_t^m\boe(t)\|_0\|\D_1\mathfrak{B}(t)\|_{L^\infty}\nonumber\\
	&\ls M_0+\delta\|\D_t^m\DD\boe(t)\|_0^2+\|\D_t^m\boe(t)\|_0^2\|\D_1\mathfrak{B}(t)\|_2\\
	&\ls M_0+\delta\|\D_t^m\DD\boe(t)\|_0^2+P\left(\sup_{t\in[0,T]}\mathfrak{E}(t)\right).\nonumber
\end{align*}                                                                                                                                                        
It follows from the Sobolev embedding theorem and \eqref{est-B-2} that
\begin{align*}
	\left|\int_{0}^{t}\mathcal{R}_b^{1,3,2}+\mathcal{R}_b^{1,3,3}\right|&\ls\int_{0}^{t}\|\D\D_1\mathfrak{B}\|_{L^4}\|\D_t^m\boe\|_{L^4}\|\D_t^m\D\boe\|_0+\int_{0}^{t}\|\D_1\mathfrak{B}\|_{L^\infty}\|\D_t^m\D\boe\|_0^2\nonumber\\
	&\quad+\int_{0}^{t}\|\D_1\mathfrak{B}\|_{L^\infty}\|\D\DD\boe\|_{L^4}\|\D_t^m\boe\|_{L^4}\|\D_t^m\D\boe\|_0\\
	&\ls \left(\|\D_1\mathfrak{B}\|_{L^2_t(\MH^2)}+\|\DD\boe\|_{L^\infty_t(\MH^2)}\|\D_1\mathfrak{B}\|_{L^2_t(H^2)}\right)\|\D_t^m\boe\|_{L^4_t(\MH^1)}^2\nonumber\\
	&\ls P\left(\sup_{t\in[0,T]}\mathfrak{E}(t)\right).\nonumber
\end{align*} 
In view of \eqref{A-priori-assum}, \eqref{est-B} and \eqref{tre}, one has 
\begin{align*}
	\left|\int_{0}^{t}\mathcal{R}_b^{1,3,4}\right|&\ls \left|\int_\Gamma\frac{\D_1\mathfrak{B}}{|\D_1\boe|^2}\D_1\eta_i\D_t^m\eta_ia_{j2}\D_t^m\eta_j\bigg|_0^t\right|\nonumber\\
	&\ls M_0+|\D_1\mathfrak{B}|_{L^4}|\D_t^m\boe|_{L^4}|\D_t^m\boe|_0\ls M_0+|\D_1\mathfrak{B}|_\hal|\D_t^m\boe|_\hal|\D_t^m\boe|_0\nonumber\\
	&\ls M_0+\|\D_1\mathfrak{B}\|_1\|\D_t^m\boe\|_1(\|\D_t^m\boe\|_0^\hal\|\D_t^m\DD\boe\|_0^\hal+\|\D_t^m\boe\|_0)\\
	&\ls M_0+\delta\|\D_t^m\DD\boe(t)\|_0^2+P(\|\boe\|_{\MH^m}^2,\|\mathfrak{B}\|_2^2)\nonumber\\
	&\ls M_0+\delta\|\D_t^m\DD\boe(t)\|_0^2+TP\left(\sup_{t\in[0,T]}\mathfrak{E}(t)\right).\nonumber
\end{align*}                                                                                          
and
\begin{align*}
	\left|\int_{0}^{t}\mathcal{R}_b^{1,3,5}\right|&\ls\int_{0}^{t}\left(|\D_t\D_1\mathfrak{B}|_{L^4}+|\D_1\mathfrak{B}\D_t\D_1\boe|_{L^4}\right)|\D_t^m\boe|_{L^4}|\D_t^m\boe|_0\nonumber\\
	&\ls \int_{0}^{t}\left(|\D_t\D_1\mathfrak{B}|_\hal+|\D_1\mathfrak{B}\D_t\D_1\boe|_\hal\right)|\D_t^m\boe|_\hal|\D_t^m\boe|_0\nonumber\\
	&\ls \int_{0}^{t}\left(\|\D_t\D_1\mathfrak{B}\|_1+\|\D_1\mathfrak{B}\D_t\D_1\boe\|_1\right)\|\D_t^m\boe\|_1(\|\D_t^m\boe\|_0^\hal\|\D_t^m\DD\boe\|_0^\hal+\|\D_t^m\boe\|_0)\\
	&\ls TP\left(\sup_{t\in[0,T]}\mathfrak{E}(t)\right).\nonumber
\end{align*}
As a consequence, we have
\begin{equation}\label{RB13}
	\left|\int_{0}^{t}\mathcal{R}_b^{1,3}\right|\ls M_0+\delta\|\D_t^m\boe(t)\|_0^2+P\left(\sup_{t\in[0,T]}\mathfrak{E}(t)\right).
\end{equation}
which, combining with \eqref{RB11} and \eqref{RB12}, yields that
\begin{equation}\label{RB1}
	\left|\int_{0}^{t}\mathcal{R}_b^1\right|\ls M_0+\delta(|\D_t^m\D_1\eta_ia_{i2}(t)|_0^2+\|\D_t^m\DD\boe(t)\|_0^2)+P\left(\sup_{t\in[0,T]}\mathfrak{E}(t)\right)
\end{equation}
For $\mathcal{R}_b^2$, we also can not use the duality argument as in \eqref{Rb2-est} since it only involves time derivatives of highest order. Instead, we use integration by parts to get that
	\begin{align}\label{RB3}
		\mathcal{R}_b^2&=\int_\Omega[\D_t^m,\mathfrak{B},a_{ij}]\D_j\mathcal{V}_i+\int_\Omega\D_j[\D_t^m,\mathfrak{B},a_{ij}]\mathcal{V}_i\nonumber\\
		&=\frac{d}{dt}\int_\Omega[\D_t^m,\mathfrak{B},a_{ij}]\D_j\D_t^m\eta_i-\int_\Omega\D_t[\D_t^m,\mathfrak{B},a_{ij}]\D_j\D_t^m\eta_i+\int_\Omega\D_j[\D_t^m,\mathfrak{B},a_{ij}]\mathcal{V}_i\\
		&\quad-\int_\Omega[\D_t^m,\mathfrak{B},a_{ij}]\D_j(\D_t^m\eta_kA_{kl}\D_lv_i)=:\mathcal{R}_b^{2,1}+\mathcal{R}_b^{2,1}+\mathcal{R}_b^{2,3}+\mathcal{R}_b^{2,4}.\nonumber
	\end{align}
By using the Sobolev embedding, \eqref{est-B-2} and \eqref{L-infL-2}, one has
\begin{align}
	\left|\int_{0}^{t}\mathcal{R}_b^{2,1}\right|&\ls \left|\int_\Omega[\D_t^m,\mathfrak{B},a_{ij}]\D_j\D_t^m\eta_i\bigg|_0^t\right|\nonumber\\
	&\ls M_0+(\|\D_t\DD\boe\|_{L^\infty_{t,x}}\|\D_t^{m-1}\mathfrak{B}\|_0+\|\D_t\mathfrak{B}\|_{L^\infty}\|\D_t^{m-1}\DD\boe\|_0)\|\D_t^m\DD\boe\|_0\nonumber\\
	&\quad+\sum_{\ell=2}^{m-2}\|\D_t^{m-\ell}\mathfrak{B}\|_{L^4}\|\D_t^\ell\DD\boe\|_{L^4}\|\D_t^m\boe\|_0\\
	&\ls M_0+\delta\|\D_t^m\DD\boe(t)\|_0^2+C_\delta\|\boe\|_{\MH^m}^2\|\mathfrak{B}\|_{\MH^{m-1}}^2.\nonumber\\
	&\ls M_0+\delta\|\D_t^m\DD\boe(t)\|_0^2+TP\left(\sup_{t\in[0,T]}\mathfrak{E}(t)\right).\nonumber
\end{align}
By using H\"{o}lder's inequality, one has
\begin{align}
	\left|\int_0^t\mathcal{R}_b^{2,2}+\mathcal{R}_b^{2,3}\right|&\ls\int_{0}^{t}\left(\|[\D_t^m,\D\mathfrak{B},\DD\boe]_0+\|[\D_t^m,\mathfrak{B},\D\DD\boe]\|_0\right)(\|\D_t^m\D\boe\|_0+\|\D_t^m\boe\|_0\|\DD\bv\|_{L^\infty_{t,x}})\nonumber\\
	&\ls T^\frac{1}{4}(\|[\D_t^m,\D\mathfrak{B},\DD\boe]\|_{L^2_t(L^2)}+\|[\D_t^m,\mathfrak{B},\D\DD\boe]\|_{L^2_t(L^2)})\cdot\nonumber\\
	&\quad\cdot(\|\D_t^m\D\boe\|_{L^4_t(L^2)}+\|\D_t^m\boe\|_{L^\infty_t(L^2)}\|\DD\bv\|_{L^4_t(H^2)})\\
	&\ls T^\frac{1}{4}P\left(\sup_{t\in[0,T]}\mathfrak{E}(t)\right).\nonumber
\end{align}
where we have used 
\begin{align*}
	\|[\D_t^m,\D\mathfrak{B},\DD\boe]\|_{L^2_t(L^2)}&\ls \|\D^2\mathfrak{B}\|_{L^\infty_{t,x}}\|\D\DD\boe\|_{L^2_t(\MH^{m-2})}+\|\D\DD\boe\|_{L^\infty_{t,x}}\|\D^2\mathfrak{B}\|_{L^2_t(\MH^{m-2})}\\
	&\ls \|\D^2\mathfrak{B}\|_{L^2_t(\MH^{m-2})}\|\D\DD\boe\|_{L^2_t(\MH^{m-2})}\ls P\left(\sup_{t\in[0,T]}\mathfrak{E}(t)\right).
\end{align*}
and 
\begin{align*}
	&\|[\D_t^m,\mathfrak{B},\D\DD\boe]\|_{L^2_t(L^2)}\ls \sum\limits_{1\leq\ell\leq m-1}\|\D_t^{m-\ell}\mathfrak{B}\D_t^\ell\D\DD\boe\|_{L^2_t(L^2)}\\
	&\ls\|\D_t\mathfrak{B}\|_{L^\infty_{t,x}}\|\D\DD\boe\|_{L^2_t(\MH^{m-1})}+\|\D\DD\boe\|_{L^\infty_{t,x}}\|\D_t\mathfrak{B}\|_{L^2_t(\MH^{m-1})}\ls P\left(\sup_{t\in[0,T]}\mathfrak{E}(t)\right).
\end{align*}
Similar, 
\begin{equation}
	\left|\int_{0}^{t}\mathcal{R}_b^{3,4}\right|\ls T^\frac{1}{4}P\left(\sup_{t\in[0,T]}\mathfrak{E}(t)\right).
\end{equation}
As a consequence,
\begin{equation}\label{RB2}
	\left|\int_{0}^{t}\mathcal{R}_b^2\right|\ls M_0+\delta\|\D_t^m\DD\boe(t)\|_0^2+T^\frac{1}{4}P\left(\sup_{t\in[0,T]}\mathfrak{E}(t)\right).
\end{equation}
To control $\mathcal{R}_b^3$, we can not use either the duality argument as \eqref{Rb4-est} due to loss of spatial derivatives, or the similar idea  as \eqref{RB3} since that $\left[\D_t^m,\D_1^2\eta_k,a_{k2}/|\D_1\boe|^3\right]$
has one more derivative than $[\D_t^m,\mathfrak{B}, a_{i2}]$. Instead, we deal with the troublesome term involving $$\D_t^{m-1}\D_1^2\eta_k\D_t(a_{k2}/|\D_1\boe|^3)$$
through combining it with $\mathcal{R}_q^3$. In view of \eqref{BB}, we have
	\begin{align*}
		\mathcal{R}_q^3&=-m\int_\Omega\D_t^mq\D_t a_{ij}\D_t^{m-1}\D_jv_idx-\sum\limits_{2\leq\ell\leq m-1}C_m^\ell\int_\Omega\D_t^mq\D_t^\ell a_{ij}\D_t^{m-\ell}\D_j v_idx\\
		&\quad+\int_\Omega[\D_t^m,a_{ij},\D_jv_i]\D_t^m\eta_kA_{kl}\D_lqdx\\
		&=-m\int_\Gamma\D_t^mq\D_t a_{i2}\D_t^{m-1}v_i+m\int_\Omega\D_j\D_t^mq\D_t a_{ij}\D_t^{m-1}v_idx\nonumber\\
		&\quad-\sum\limits_{2\leq\ell\leq m-1}C_m^\ell\int_\Omega\D_t^mq\D_t^\ell a_{ij}\D_t^{m-\ell}\D_j v_idx+\int_\Omega[\D_t^m,a_{ij},\D_jv_i]\D_t^m\eta_kA_{kl}\D_lqdx\\
		&=m\int_\Gamma\D_t^{m}\left(\frac{\D_1^2\eta_ka_{k2}}{|\D_1\boe|^3}\right)\D_t a_{i2}\D_t^{m-1}v_i-m\int_\Gamma\D_t^{m}\mathfrak{B}\D_ta_{ij}\D_t^{m-1}v_i\\
	    &\quad+m\int_\Omega\D_j\D_t^mq\D_t a_{ij}\D_t^{m-1}v_idx-\sum\limits_{2\leq\ell\leq m-1}C_m^\ell\int_\Omega\D_t^mq\D_t^\ell a_{ij}\D_t^{m-\ell}\D_j v_idx\\
	    &\quad+\int_\Omega[\D_t^m,a_{ij},\D_jv_i]\D_t^m\eta_kA_{kl}\D_lqdx=:\mathcal{R}_q^{3,1}+\cdots+\mathcal{R}_q^{3,5}
	\end{align*}
which, summing with  $-\mathcal{R}_b^3$, yields
	\begin{align}\label{MRB3}
	\mathcal{R}_b^3+\mathcal{R}_q^{3,1}&=m\int_\Gamma\D_t\left(\frac{a_{k2}}{|\D_1\boe|^3}\right)\D_t^{m-1}\D_1^2\eta_ka_{i2}\mathcal{V}_i+\sum\limits_{2\leq\ell\leq m-1}C_m^\ell\int_\Gamma\D_t^{m-\ell}\D_1^2\eta_k\D_t^\ell\left(\frac{a_{k2}}{|\D_1\boe|^3}\right)a_{i2}\mathcal{V}_i\nonumber\\
	&+m\int_\Gamma\D_t a_{i2}\D_t^{m-1}v_i\frac{\D_t^m\D_1^2\eta_ka_{k2}}{|\D_1\boe|^3}+m\int_\Gamma\left[\D_t^m,\frac{a_{k2}}{|\D_1\boe|^3}\right]\D_1^2\eta_k\D_ta_{i2}\D_t^{m-1}v_i\nonumber\\
	&=\frac{d}{dt}\int_\Gamma\frac{m}{|\D_1\boe|^3}\D_t^{m-1}\D_1^2\eta_ka_{k2}\D_t^{m-1}v_i \D_t a_{i2}-m\int_\Gamma\frac{\D_t^{m-1}\D_1^2\eta_k\D_t^mv_i}{|\D_1\boe|^3}\left(\D_ta_{i2}a_{k2}-a_{i2}\D_ta_{k2}\right)\nonumber\\
	&-m\int_\Gamma\D_t^{m-1}\D_1^2\eta_k\D_t^{m-1}v_i\D_t\left(\frac{\D_ta_{i2}a_{k2}}{|\D_1\boe|^3}\right)+m\int_\Gamma\D_t\left(\frac{1}{|\D_1\boe|^3}\right)\D_t^{m-1}\D_1^2\eta_ka_{k2}\D_t^mv_ia_{i2}\nonumber\\
	&+\sum\limits_{2\leq\ell\leq m-1}C_m^\ell\int_\Gamma\D_t^{m-\ell}\D_1^2\eta_k\D_t^\ell\left(\frac{a_{k2}}{|\D_1\boe|^3}\right)a_{i2}\mathcal{V}_i+m\int_\Gamma\left[\D_t^m,\frac{a_{k2}}{|\D_1\boe|^3}\right]\D_1^2\eta_k\D_ta_{i2}\D_t^{m-1}v_i\nonumber\\
	&-m\int_\Gamma\D_t\left(\frac{a_{k2}}{|\D_1\boe|^3}\right)\D_t^{m-1}\D_1^2\eta_ka_{i2}\D_t^m\eta_j A_{jl}\D_lv_i\nonumber\\
	&=\mathcal{R}_b^{3,1}+\cdots+\mathcal{R}_b^{3,7}.
	\end{align}
By using H\"{o}lder's inequality, \eqref{fra-emb}, trace theorem and \eqref{tan-est-1}, one has
	\begin{align}\label{RB41}
		\left|\int_{0}^{t}\mathcal{R}_b^{3,1}\right|&\ls M_0+|\D_t^{m-1}\D_1^2\eta_ka_{k2}|_0|\D_t^{m-1}v|_{L^4}|\D_t\D_1\eta|_{L^4}\nonumber\\
		&\ls M_0+|\D_t^{m-1}\D_1^2\eta_ka_{k2}|_0\|\D_t^m\eta\|_1\|\D_t\D_1\eta\|_1\\
		&\ls M_0+\delta\|\DD\D_t^m\eta\|_0+P\left(\sup_{t\in[0,T]}\mathfrak{E}(t)\right).\nonumber
	\end{align}
The term $\mathcal{R}_b^{3,2}$ vanishes when $i=k$, but reduces to the following equation, when $i\neq k$:
\begin{align*}
	\mathcal{R}_b^{3,2}&=-m\int_\Gamma\frac{1}{|\D_1\boe|^3}\left(\D_t^{m-1}\D_1^2\eta_2\D_t^mv_1-\D_t^{m-1}\D_1^2\eta_1\D_t^mv_2\right)(\D_1a_{12}a_{22}-a_{12}\D_ta_{22})\nonumber\\
	&=-m\frac{d}{dt}\int_\Gamma\frac{\D_t\D_1\eta_ia_{i2}}{|\D_1\boe|^3}\D_t^{m-1}\D_1a_{k2}\D_t^m\eta_k+m\int_\Gamma\D_t^{m-1}\D_1a_{k2}\D_t^m\eta_k\D_t\left(\frac{\D_t\D_1\eta_ia_{i2}}{|\D_1\boe|^3}\right)\\
	&\quad-m\int_\Gamma\D_t^ma_{k2}\D_t^m\eta_k\D_1\left(\frac{\D_t\D_1\eta_ia_{i2}}{|\D_1\boe|^3}\right)=\mathcal{R}_b^{3,2,1}+\mathcal{R}_b^{3,2,2}+\mathcal{R}_b^{3,2,3}.
\end{align*}
where we have used the following identities
\begin{align*}
	&\D_t^{m-1}\D_1^2\eta_2\D_t^mv_1-\D_t^{m-1}\D_1^2\eta_1\D_t^mv_2=-\D_t^{m-1}\D_1a_{k2}\D_t^mv_k,\\&\D_1a_{12}a_{22}-a_{12}\D_ta_{22}=-\D_t\D_1\eta_ia_{i2},\quad\D_t^ma_{k2}\D_t^m\D_1\eta_k=0.
\end{align*} 
Then,  it follows from \eqref{co123}, \eqref{fra-emb} and trace theorem and \eqref{tan-est-1} that
\begin{align*}
	\left|\int_{0}^{t}\mathcal{R}_b^{3,2,1}\right|&\ls M_0+|\D_t^{m-1}\D_1a_{k2}|_{-\frac{1}{2}}|\D_t^m\eta_k|_\hal|\D_t\D_1\eta_ia_{i2}|_1\nonumber\\
	&\ls M_0+|\D_t^{m-1}a_{k2}|_\hal|\D_t^m\eta_k|_\hal(|\D_t\D_1^2\eta_ia_{i2}|_0+|\D_t\D_1\eta_i\D_1a_{i2}|+|\D_t\D_1\eta_ia_{i2}|_0)\nonumber\\
	&\ls M_0+\|\D_t^{m-1}\D_1\boe\|_1\|\D_t^m\boe\|_1(|\D_t\D_1^2\eta_ia_{i2}|_0+|\D_t\D_1\eta_ia_{i2}|_0+\|\bar{\D}\D_1\boe\|_1^2)\\
	&\ls M_0+\delta\|\DD\D_t^m\boe\|_0^2+P\left(\sup_{t\in[0,T]}\mathfrak{E}(t)\right).\nonumber
\end{align*}
Similarly, one has
\begin{align*}
\left|\int_{0}^{t}\mathcal{R}_b^{3,2,2}+\mathcal{R}_b^{3,2,3}\right|&\ls \int_{0}^{t}\|\D_t^{m-1}\bar{\D}\boe\|_1\|\D_t^m\boe\|_1(|\bar{\D}^2\D_1^2\eta_ia_{i2}|_0+|\bar{\D}^3\boe\bar{\D}^2\boe|_0+|\bar{\D}^2\D_1\eta_ia_{i2}|_0+|\bar{\D}^2\boe|_{L^4}^2)\nonumber\\
&\ls T^\frac{1}{2}P\left(\sup_{t\in[0,T]}\mathfrak{E}(t)\right).	
\end{align*}
Therefore, we have
\begin{equation}\label{RB42}
	\left|\int_{0}^{t}\mathcal{R}_b^{3,2}\right|\ls M_0+\delta\|\DD\D_t^m\eta\|_0+P\left(\sup_{t\in[0,T]}\mathfrak{E}(t)\right).
\end{equation}
It follows from \eqref{co123}, \eqref{fra-emb} and trace theorem and \eqref{tan-est-1} that
\begin{align}\label{RB43}
	\left|\int_{0}^{t}\mathcal{R}_b^{3,3}\right|&\ls \int_{0}^{t}|\D_t^{m-1}\D_1^2\boe|_{-\frac{1}{2}}|\D_t^m\boe|_\frac{1}{2}|\D_1\bv\D_1\bv|_1+|\D_t^{m-1}\D_1^2\eta_ka_{k2}|_0|\D_t^{m-1}\bv|_{L^4}|\D_t^2\D_1\boe|_{L^4}\nonumber\\
	&\ls \int_{0}^{t}\|\D_t^{m-1}\D_1\boe\|_1\|\D_t^m\boe\|_1(\|\D_1\bv\|_1+\|\D_1^2\bv\|_1)\|\D_1\bv\|_1+|\D_t^{m-1}\D_1^2\eta_ka_{k2}|_0\|\D_t^m\boe\|_1\|\D_t^2\D_1\boe\|_1\nonumber\\
	&\ls T^\frac{1}{4}(\|\D_t^{m-1}\D_1\boe\|_{L^\infty_t(H^1)}\|\boe\|_{L^\infty_t(\MH^3)}+|\D_t^{m-1}\D_1^2\eta_ka_{k2}|_{L^\infty(L^2)})\|\D_t^m\boe\|_{L^4_t(H^1)}\|\D_1\bv\|_{L^2_t(\MH^2)}\nonumber\\
	&\ls T^\frac{1}{4}P\left(\sup_{t\in[0,T]}\mathfrak{E}(t)\right).	 
\end{align}
For $\mathcal{R}_b^{3,4}$, we can not use the same argument as \eqref{RB43} since we can not bound $|\D_t^m\bv|_\hal$. Instead, we use integration by parts to get
	\begin{align*}
		\mathcal{R}_b^{3,4}&=-m\int_\Gamma\D_t\left(\frac{1}{|\D_1\boe|^3}\right)\D_t^{m-1}\D_1\eta_ka_{k2}\D_t^{m+1}\D_1\eta_ia_{i2}\\
		&\quad-m\int_\Gamma\D_1\left(\D_t\left(\frac{1}{|\D_1\boe|^3}\right)a_{k2}a_{i2}\right)\D_t^{m-1}\D_1\eta_k\D_t^mv_i\\
		&=-m\frac{d}{dt}\int_\Gamma\D_t\left(\frac{1}{|\D_1\boe|^3}\right)\D_t^{m-1}\D_1\eta_ka_{k2}\D_t^m\D_1\eta_ia_{i2}+m\int_\Gamma\D_t\left(\frac{1}{|\D_1\boe|^3}\right)|\D_t^m\D_1\eta_ka_{k2}|^2\\
		&\quad+m\int_\Gamma\D_t\left(\D_t\left(\frac{1}{|\D_1\boe|^3}\right)a_{k2}a_{i2}\right)\D_t^{m-1}\D_1\eta_k\D_t^m\D_1\eta_i\\
		&\quad-m\int_\Gamma\D_1\left(\D_t\left(\frac{1}{|\D_1\boe|^3}\right)a_{k2}a_{i2}\right)\D_t^{m-1}\D_1\eta_k\D_t^mv_i\\
		&=\mathcal{R}_b^{3,4,1}+\cdots\mathcal{R}_b^{3,4,4}.
	\end{align*}
Then, we can obtain from H\"{o}lder's inequality, \eqref{fra-emb} and the trace theorem that
\begin{align}\label{RB441}
		\left|\int_{0}^{t}\mathcal{R}_b^{3,4,1}\right|&\ls M_0+|\D_t^m\D_1\eta_ia_{i2}|_0|\D_t^{m-1}\D_1\boe|_{L^4}|\D_t\D_1\boe|_{L^4}\nonumber\\
		&\ls M_0+|\D_t^m\D_1\eta_ia_{i2}|_0\|\D_t^{m-1}\D_1\boe\|_1\|\D_t\D_1\boe\|_1\\
		&\ls M_0+\delta|\D_t^m\D_1\eta_ia_{i2}|_0^2+P\left(\sup_{t\in[0,T]}\mathfrak{E}(t)\right).\nonumber
\end{align}
and
	\begin{align}\label{RB442}
		\left|\int_{0}^{t}\mathcal{R}_b^{3,4,2}\right|&\ls \int_{0}^{t}|\D_t\D_1\eta|_{L^\infty}|\D_t^m\D_1\eta_ia_{i2}|_0^2\nonumber\\
		&\ls|\D_t^m\D_1\eta_ia_{i2}|_{L^4_t(L^2)}^2\|\DD\eta\|_{L^2(\MH^{m-1})}\\
		&\ls T^\frac{1}{2}P\left(\sup_{t\in[0,T]}\mathfrak{E}(t)\right). \nonumber
	\end{align}
Moreover, it follows from H\"{o}lder's inequality, \eqref{fra-emb}, \eqref{co123} and the trace theorem that
	\begin{align}\label{RB443}
		\left|\int_{0}^{t}\mathcal{R}_b^{3,4,3}\right|&\ls\left|\int_{0}^{t}\int_\Gamma\D_t\left(\D_t\left(\frac{1}{|\D_1\boe|^3}\right)a_{k2}\right)\D_t^{m-1}\D_1\eta_k\D_t^m\D_1\eta_ia_{i2}\right|\nonumber\\
		&\quad+\left|\int_{0}^{t}\int_\Gamma\D_t\left(\frac{1}{|\D_1\boe|^3}\right)\D_t^{m-1}\D_1\eta_ka_{k2}\D_t^m\D_1\eta_i\D_ta_{i2}\right|\nonumber\\
		&\ls \int_{0}^{t}|\D_t^m\D_1\eta_ia_{i2}|_0|\D_t^{m-1}\D_1\boe|_{L^4}(|\D_t^2\D_1\boe|_{L^4}+|\D_t\D_1\boe\D_t\D_1\boe|_{L^4})\nonumber\\
		&\quad+\int_{0}^{t}|\D_t^{m-1}\D_1\eta_ia_{i2}|_1|\D_t^m\D_1\boe|_{-\frac{1}{2}}|\D_t\D_1\boe\D_1\boe|_\frac{1}{2}|\D_t\D_1\boe\D_1\boe|_1\\
		&\ls\int_{0}^{t}|\D_t^m\D_1\eta_ia_{i2}|_0\|\D_t^{m-1}\D_1\boe\|_1(\|\D_t^2\D_1\boe\|_1+\|\D_t\D_1\boe\D_t\D_1\boe\|_1)\nonumber\\
		&\quad+\int_{0}^{t}|\D_t^{m-1}\D_1\eta_ia_{i2}|_1\|\D_t^m\boe\|_1\|\D_t\D_1\boe\D_1\boe\|_1|\D_t\D_1\boe\D_1\boe|_1\nonumber\\
		&\ls T^\frac{1}{4}P\left(\sup_{t\in[0,T]}\mathfrak{E}(t)\right).\nonumber
	\end{align}
To bound $\mathcal{R}_b^{3,4,4}$, similar to \eqref{RB3}, we use integration by parts to reduce this boundary integral to the volume one as
\begin{align*}
	\mathcal{R}_b^{3,4,4}=&-m\int_\Omega\D_1\left(\D_t\left(\frac{1}{|\D_1\boe|^3}\right)a_{k2}a_{ij}\right)\D_t^{m-1}\D_1\eta_k\D_t^m\D_jv_i\\
	&+m\int_\Omega\D_1\left(\D_t\left(\frac{1}{|\D_1\boe|^3}\right)a_{k2}a_{ij}\right)\D_t^{m-1}\D_1\D_j\eta_k\D_t^mv_i\\
	&+m\int_\Omega\D_j\D_1\left(\D_t\left(\frac{1}{|\D_1\boe|^3}\right)a_{k2}a_{ij}\right)\D_t^{m-1}\D_1\eta_k\D_t^mv_i\\
	=&-m\frac{d}{dt}\int_\Omega\D_1\left(\D_t\left(\frac{1}{|\D_1\boe|^3}\right)a_{k2}a_{ij}\right)\D_t^{m-1}\D_1\eta_k\D_t^m\D_j\eta_i\\
	&+m\int_\Omega\D_t\D_1\left(\D_t\left(\frac{1}{|\D_1\boe|^3}\right)a_{k2}a_{ij}\right)\D_t^{m-1}\D_1\eta_k\D_t^m\D_j\eta_i\\
	&+m\int_\Omega\D_j\D_1\left(\D_t\left(\frac{1}{|\D_1\boe|^3}\right)a_{k2}a_{ij}\right)\D_t^{m-1}\D_1\eta_k\D_t^mv_i\\
	&+m\int_\Omega\D_1\left(\D_t\left(\frac{1}{|\D_1\boe|^3}\right)a_{k2}a_{ij}\right)\D_t^{m-1}\D_1\D_j\eta_k\D_t^mv_i\\
	=&:\mathcal{R}_b^{3,4,4,1}+\cdots+\mathcal{R}_b^{3,4,4,4}
\end{align*}
By using H\"{o}lder's inequality, one has
\begin{align*}
	\left|\int_{0}^{t}\mathcal{R}_b^{3,4,4,1}\right|&\ls \left|\int_\Omega\D_1\left(\D_t\left(\frac{1}{|\D_1\boe|^3}\right)a_{k2}a_{ij}\right)\D_t^{m-1}\D_1\eta_k\D_t^m\D_j\eta_i\bigg|_0^t\right|\nonumber\\
	&\ls M_0+(\|\D_t\D_1^2\boe\|_{L^4}+\|\D_t\D_1\boe\D_1\DD\boe\|_{L^4})\|\D_t^{m-1}\D_1\boe\|_{L^4}\|\D_t^m\DD\boe\|_0\\
	&\ls M_0+\delta\|\D_t^m\DD\boe(t)\|_0+P\left(\sup_{t\in[0,T]}\mathfrak{E}(t)\right).\nonumber
\end{align*}
\begin{align*}
	\left|\int_{0}^{t}\mathcal{R}_b^{3,4,4,2}+\mathcal{R}_b^{3,4,4,3}\right|&\ls\int_{0}^{t}(\|\D_t\D_1^2\D\boe\|_{L^4}+\|\D^3\boe\D^2\boe\|_{L^4}+\|\D^2\boe\|_{L^6}^3)\|\D_t^{m-1}\D_1\boe\|_{L^4}\|\D_t^m\DD\boe\|_0 \nonumber\\
	&\ls\|\D_t^m\DD\boe\|_{L^4_t(L^2)}\|\D_t^{m-1}\D_1\boe\|_{L^4_t(H^1)}(\|\DD\boe\|_{L^2_t(\MH^4)}+T^\frac{1}{2}P(\|\boe\|_{L^\infty_t(\MH^m)}))\nonumber\\
	&\ls T^\frac{1}{4}P\left(\sup_{t\in[0,T]}\mathfrak{E}(t)\right).
\end{align*}
and
\begin{align*}
	\left|\int_{0}^{t}\mathcal{R}_b^{3,4,4,4}\right|&\ls \int_{0}^{t}(\|\D_t\D_1^2\boe\|_{L^\infty}+\|\bar{\D}\DD\boe\|_{L^\infty}^2)\|\D_t^{m-1}\D_1\DD\boe\|_0\|\D_t^m\bv\|_0\nonumber\\
	&\ls \|\D_t^m\bv\|_{L^4_t(L^2)}\|\D_t^{m-1}\D_1\DD\boe\|_{L^4_t(L^2)}(\|\DD\boe\|_{L^2_t(\MH^4)}+T^\frac{1}{2}P(\|\boe\|_{L^\infty_t(\MH^4)}))\\
	&\ls T^\frac{1}{4}P\left(\sup_{t\in[0,T]}\mathfrak{E}(t)\right).\nonumber
\end{align*}
As a consequence, we obtain that
\begin{equation}\label{RB444}
		\left|\int_{0}^{t}\mathcal{R}_b^{3,4,4}\right|\ls M_0+\delta\|\D_t^m\DD\boe(t)\|_0^2+P\left(\sup_{t\in[0,T]}\mathfrak{E}(t)\right).
\end{equation}
Plugging \eqref{RB441}-\eqref{RB444} yields that
\begin{equation}\label{RB44}
	\left|\int_{0}^{t}\mathcal{R}_{b}^{3,4}\right|\ls M_0+\delta(\|\D_t^m\DD\boe(t)\|_0^2+|\D_t^m\D_1\eta_ia_{i2}(t)|_0^2)+P\left(\sup_{t\in[0,T]}\mathfrak{E}(t)\right).
\end{equation}
The remaining terms $\mathcal{R}_b^{3,5},\mathcal{R}_b^{3,6},\mathcal{R}_b^{3,7}$ are easy to bounded by the similar argument used for $\mathcal{R}_b^{3,3},\mathcal{R}_b^{3,4}$. We omit the details and conclude that
\begin{equation}\label{RB45-47}
	\left|\int_{0}^{t}\mathcal{R}_b^{3,5}+\mathcal{R}_b^{3,6}+\mathcal{R}_b^{3,7}\right|\ls M_0+\delta(\|\D_t^m\DD\boe\|_0^2+|\D_t^m\D_1\eta_ia_{i2}|_0^2)+P\left(\sup_{t\in[0,T]}\mathfrak{E}(t)\right).
\end{equation}
Substituting \eqref{RB41}, \eqref{RB42}, \eqref{RB43}, \eqref{RB44}, \eqref{RB45-47} into \eqref{MRB3} yields that
\begin{equation*}
	\left|\int_{0}^{t}\mathcal{R}_b^3+\mathcal{R}_q^{3,1}\right|\ls M_0+\delta(\|\D_t^m\DD\boe(t)\|_0^2+|\D_t^m\D_1\eta_ia_{i2}(t)|_0^2)+P\left(\sup_{t\in[0,T]}\mathfrak{E}(t)\right).
\end{equation*}
Moreover, it is not difficult to prove that
\begin{equation*}
	\left|\int_{0}^{t}\mathcal{R}_q^{3,2}+\cdots+\mathcal{R}_q^{3,5}\right|\ls P\left(\sup_{t\in[0,T]}\mathfrak{E}(t)\right).
\end{equation*}
So that, one gets
\begin{equation}\label{RB3-est}
	\left|\int_{0}^{t}\mathcal{R}_b^3+\mathcal{R}_q^3\right|\ls M_0+\delta(\|\D_t^m\DD\boe(t)\|_0^2+|\D_t^m\D_1\eta_ia_{i2}(t)|_0^2)+P\left(\sup_{t\in[0,T]}\mathfrak{E}(t)\right).
\end{equation}
For $\mathcal{R}_b^4$, it follows from \eqref{tre} that
\begin{align}\label{RB4}
	\left|\int_{0}^{t}\mathcal{R}_b^{4}\right|&\ls \int_{0}^{t}|\D_t\D_1\boe|_{L^\infty}|\D_t^m\D_1\eta_ka_{k2}|_0^2\ls \int_{0}^{t}|\D_t\D_1\boe|_1|\D_t^m\D_1\eta_ka_{k2}|_0^2\nonumber\\
	&\ls|\D_t\D_1\boe|_{L^2_t(H^1)}|\D_t^m\D_1\eta_ka_{k2}|_{L^4_t(L^2)}^2\ls T^\frac{1}{2}P\left(\sup_{t\in[0,T]}\mathfrak{E}(t)\right).
\end{align}
For $\mathcal{R}_b^5$, by tedious but not difficult calculation, we have
	\begin{align*}
		\mathcal{R}_b^5&=\int_\Gamma\frac{\D_t^ma_{k2}\D_1^2\eta_k-\D_t^m\D_1\eta_k\D_1a_{k2}}{|\D_1\boe|^3}a_{i2}\mathcal{V}_i-\int_\Gamma\frac{\D_t^m\D_1\eta_ka_{k2}}{|\D_1\boe|^3}\D_1a_{i2}\mathcal{V}_i\nonumber\\
		&\quad+3\int_\Gamma\frac{\D_t^m\D_1\eta_k\D_1a_{j2}a_{j2}a_{k2}-\D_t^ma_{j2}a_{j2}a_{k2}\D_1^2\eta_k}{|\D_1\boe|^5}a_{i2}\mathcal{V}_i\nonumber\\
		&\quad+\underbrace{\int_\Gamma\left(\left[\D_t^m,a_{k2},\frac{1}{|\D_1\boe|^3}\right]-3\left[\D_t^{m-1},\D_t a_{j2},\frac{a_{j2}}{|\D_1\eta|^5}\right]a_{k2}\right)\D_1^2\eta_k a_{i2}\mathcal{V}_i}_{\mathcal{R}_b^{5,2}}\\
		&\quad-\underbrace{3\int_\Gamma\D_ta_{j2}\D_t^{m-1}\left(\frac{a_{j2}}{|\D_1\eta|^5}\right)\D_1^2\eta_k a_{k2}a_{i2}\mathcal{V}_i}_{\mathcal{R}_b^{5,3}}\nonumber\\
		&=\underbrace{\int_\Gamma\frac{\D_1^2\eta_k a_{k2}}{|\D_1\boe|^3}\D_t^m a_{i2}\mathcal{V}_i}_{\mathcal{R}_b^{5,1}}+\mathcal{R}_b^{5,2}+\mathcal{R}_b^{5,3},\nonumber
	\end{align*}
where we have used  \eqref{decomp} and the following identity
\begin{equation*}
	\D_t^m\D_1\eta_k\D_1a_{j2}a_{j2}a_{k2}-\D_t^ma_{j2}a_{j2}a_{k2}\D_1^2\eta_k=|\D_1\boe|^2\D_t^m\D_1\eta_k\D_1a_{k2}.
\end{equation*}
By the similar argument as in the estimate of $\mathcal{R}_b^1$, one has
\begin{equation*}
	\left|\int_{0}^{t}\mathcal{R}_b^{5,1}\right|\ls \delta\|\DD\D_t^m\boe(t)\|_0^2+P\left(\sup_{t\in[0,T]}\mathfrak{E}(t)\right).
\end{equation*}
By the similar argument as in the estimate of $\mathcal{R}_b^2$ and $\mathcal{R}_b^3$, it is not difficult to get
\begin{equation*}
	\left|\int_{0}^{t}\mathcal{R}_b^{5,2}\right|+\left|\int_{0}^{t}\mathcal{R}_b^{5,3}\right|\ls P\left(\sup_{t\in[0,T]}\mathfrak{E}(t)\right).
\end{equation*}
Thus, we obtain that
\begin{equation}\label{RB5}
	\left|\int_{0}^{t}\mathcal{R}_b^5\right|\ls\delta\|\DD\D_t^m\boe\|_0^2+P\left(\sup_{t\in[0,T]}\mathfrak{E}(t)\right). 
\end{equation}
For $\mathcal{R}_b^6$, it follows from \eqref{tre} that
\begin{align}\label{RB6}
	\left|\int_{0}^{t}\mathcal{R}_b^6\right|\ls \int_{0}^{t}|\D_t^m\D_1\eta_ka_{k2}|_0|\D_t^m\boe|_0(|\D_1\DD\bv|_{L^\infty}+|\D_1\DD\boe|_{L^\infty}|\DD\bv|_{L^\infty})\ls P\left(\sup_{t\in[0,T]}\mathfrak{E}(t)\right).
\end{align}
For $\mathcal{R}_b^7$, we obtain from integration by parts and the Piola identity \eqref{Piola} that
\begin{align}\label{RB7-0}
	\mathcal{R}_b^7&=\int_\Gamma\D_t^m\eta_kA_{kl}\D_lqa_{i2}\D_t^mv_i-\int_\Gamma\D_t^m\eta_kA_{kl}\D_lqa_{i2}\D_t^m\eta_kA_{kl}\D_l v_i\nonumber\\
	&=\int_\Omega\D_t^m\eta_kA_{kl}\D_lqa_{ij}\D_j\D_t^mv_idxds+\int_\Omega\D_j(\D_t^m\eta_kA_{kl}\D_lq)a_{ij}\D_t^mv_idxds\nonumber\\
	&\quad-\int_\Gamma\D_t^m\eta_kA_{kl}\D_lqa_{i2}\D_t^m\eta_kA_{kl}\D_l v_i\\
	&=\underbrace{\frac{d}{dt}\int_\Omega\D_t^m\eta_kA_{kl}\D_lqa_{ij}\D_j\D_t^m\eta_idx}_{\mathcal{R}_b^{7,1}}-\underbrace{\int_\Omega\D_t(\D_t^m\eta_kA_{kl}\D_lqa_{ij})\D_j\D_t^m\eta_i dxds}_{\mathcal{R}_b^{7,2}}\nonumber\\
	&\quad+\underbrace{\int_\Omega\D_j(\D_t^m\eta_kA_{kl}\D_lq)a_{ij}\D_t^mv_idxds}_{\mathcal{R}_b^{7,3}}-\underbrace{\int_\Gamma\D_t^m\eta_kA_{kl}\D_lqa_{i2}\D_t^m\eta_kA_{kl}\D_l v_i}_{\mathcal{R}_b^{7,4}}.\nonumber
\end{align}
By interpolation inequality, one has
\begin{align}\label{RB7-1}
	\left|\int_{0}^{t}\mathcal{R}_b^{7,1}\right|&\ls M_0+(\|\DD\D_t^m\boe\|_0^\frac{1}{2}\|\D_t^m\boe\|_0^\frac{1}{2}+\|\D_t^m\boe\|_0)(\|\DD^2q\|_0^\frac{1}{2}\|\DD q\|_{0}^\frac{1}{2}+\|\DD q\|_0)\|\DD\D_t^m\eta\|_0\nonumber\\
	&\ls\delta\|\DD\D_t^m\boe\|_0^2+M_0+P(\|\boe\|_{\MH^m}^2,\|q\|_{\MH^2}^2)\\
	&\ls\delta\|\DD\D_t^m\boe\|_0^2+M_0+TP\left(\sup_{t\in[0,T]}\mathfrak{E}(t)\right).\nonumber
\end{align}
In view of \eqref{A-priori-assum} and the Sobolev embedding, one has
\begin{align}\label{RB7-2}
	\left|\int_{0}^{t}\mathcal{R}_b^{7,2}\right|+\left|\int_{0}^{t}\mathcal{R}_b^{7,3}\right|&\ls \int_{0}^{t}\|\D_t^m \bv\|_0\|\DD q\|_{L^\infty}\|\D_t^m\DD\boe\|_0\nonumber\\
	&+\int_{0}^{t}\|\D_t^m\boe\|_{L^4}\|A\|_{L^\infty}\|\D\DD q\|_{L^4}(\|\D_t^m\DD\boe\|_0+\|\D_t^m\bv\|_0)\nonumber\\
	&\quad+\int_{0}^{t}\|\D_t^m\boe\|_{L^4}\|\DD q\|_{L^\infty}\|\D\DD\boe\|_{L^4}(\|\D_t^m\DD\boe\|_0+\|\D_t^m\bv\|_0)\\
	&\ls \|\DD\boe\|_{L^\infty(\MH^2)}(\|\D_t^m \bv\|_{L^4_t(L^2_x)}+\|\D_t^m\DD\boe\|_{L^4_t(L^2_x)})\|\D_t^m\boe\|_{L^4_t(H^1_x)}\|\DD q\|_{L^2_t (\MH^2)}\nonumber\\
	&\ls P\left(\sup_{t\in[0,T]}\mathfrak{E}(t)\right).\nonumber
\end{align}
On the other hand,
\begin{align}\label{RB7-3}
	\left|\int_{0}^{t}\mathcal{R}_b^{7,4}\right|&\ls\int_{0}^{t}|\D_t^m\boe|_{L^4}^2|A_{kl}\D_lqa_{i2}A_{jk}\D_k v_i|_0\ls \int_{0}^{t}\|\D_t^m\boe\|_{H^1}^2\|\D_lq\D_kv_i\|_{H^1}\nonumber\\
	&\ls \|\D_t^m\boe\|_{L^4_t(H^1)}^2\|\DD q\|_{L^2(\MH^2)}\|\DD \bv\|_{L^2(\MH^2)}\ls P\left(\sup_{t\in[0,T]}\mathfrak{E}(t)\right).
\end{align}
Plugging \eqref{RB7-1}-\eqref{RB7-3} into \eqref{RB7-0} yields that
\begin{equation}\label{RB7}
	\left|\int_{0}^{t}\mathcal{R}_b^7\right|\ls\delta\|\DD\D_t^m\boe\|_0^2+M_0+TP\left(\sup_{t\in[0,T]}\mathfrak{E}(t)\right).
\end{equation}
Finally, the volume integral terms, $\mathcal{R}_\eta^i$, $\mathcal{R}_q^i$, $\mathcal{R}_c^i$ and $\mathcal{R}_\vep^i$ can be controlled as follows. It is clear that
\begin{align}\label{REta-1}
	\left|\int_{0}^{t}\mathcal{R}_\eta^1\right|&\ls\int_{0}^{t}\|\D_t^m\bv\|_0^2\|\DD\bv\|_{L^\infty}+\|\D_t^m\bv\|_0\|\D_t^m\boe\|_0\|\DD\bv\|_{L^\infty}^2+\|\D_t^m\boe\|_0^2\|\DD\bv\|_{L^\infty}^3\nonumber\\
	&\quad+\int_{0}^{t}\|\D_t^m\boe\|_{L^4}\|\DD\D_t\bv\|_{L^4}(\|\D_t^m\bv\|_0+\|\D_t^m\boe\|_0\|\DD\bv\|_{L^\infty})\nonumber\\
	&\ls \|\D_t^m\bv\|_{L^4_t(L^2)}^2\|\DD\bv\|_{L^2_t(H^2)}+\|\D_t^m\bv\|_{L^4_t(L^2)}\|\D_t^m\boe\|_{L^4_t(L^2)}\|\DD\bv\|_{L^4_t(H^2)}^2\\
	&\quad+\|\D_t^m\boe\|_{L^4_t(L^2)}^2\|\DD\bv\|_{L^6_t(H^2)}^3+\|\D_t^m\boe\|_{L^4_t(H^1)}\|\DD\D_t\bv\|_{L^2_t(H^1)}\|\D_t^m\bv\|_{L^4_t(L^2)}\nonumber\\
	&\quad+\|\D_t^m\boe\|_{L^4_t(H^1)}\|\DD\D_t\bv\|_{L^2_t(H^1)}\|\D_t^m\boe\|_{L^\infty_t(L^2)}\|\DD\bv\|_{L^4_t(H^2)}\nonumber\\
	&\ls T^\frac{1}{2}P\left(\sup_{t\in[0,T]}\mathfrak{E}(t)\right).\nonumber
\end{align}
It follows from \eqref{co2}, the Sobolev embedding, and \eqref{m-J-q} that
\begin{align}\label{RQ1}
	\left|\int_{0}^{t}R_\eta^2\right|&\ls \|[\D_t^m,a_{ik},\D_kq]\|_{L^2_t(L^2_x)}\|\mathcal{V}_i\|_{L^2_t(L^2_x)}\nonumber\\
	&\ls (\|\D\DD q\|_{L^\infty_{t,x}}\|\D\DD\boe\|_{L^2_t(\MH^{m-2})}+\|\D\DD\boe\|_{L^\infty_{t,x}}\|\D\DD q\|_{L^2_t(\MH^{m-2})})\cdot\nonumber\\
	&\quad\cdot(\|\D_t^m\bv\|_{L^4_t(L^2)}T^\frac{1}{4}+\|\D_t^m\boe\|_{L^4(L^2)}\|\DD\bv\|_{L^4(L^\infty)})\\
	&\ls T^\frac{1}{4}\|\DD q\|_{L^2_t(\MH^{m-1})}\|\DD\boe\|_{L^2_t(\MH^m)}\|\D_t^m\bv\|_{L^4_t(L^2)}\ls T^\frac{3}{4}P\left(\sup_{t\in[0,T]}\mathfrak{E}(t)\right).\nonumber
\end{align}
Similar to \eqref{REta-1}, one has 
\begin{equation}
\left|\int_{0}^{t}\mathcal{R}_\eta^3+\mathcal{R}_q^1\right|\ls T^\frac{1}{2}P\left(\sup_{t\in[0,T]}\mathfrak{E}(t)\right).
\end{equation}
Similar to \eqref{RQ1}, we can obtain from \eqref{co1}, the Sobolev embedding, and \eqref{m-J-q} that
\begin{align}
	\left|\int_{0}^{t}\mathcal{R}_q^2\right|&\ls T^\frac{3}{4}P\left(\sup_{t\in[0,T]}\mathfrak{E}(t)\right).\nonumber
\end{align}
 It is obvious that
\begin{align}\label{RQ4}
        \left|\int_{0}^{t}\mathcal{R}_q^4\right|\ls \|J_t\|_{L^\infty_{t,x}}\|\D_t^mq\|_{L^2_t(L^2_x)}\ls\|J_t\|_{L^2_t(\MH^{m-1})}\|\D_t^mq\|_{L^4_t(L^2_x)}T^\frac{1}{4}\ls T^\frac{3}{4}P\left(\sup_{t\in[0,T]}\mathfrak{E}(t)\right).
\end{align}
By using \eqref{commutator-est}, one has
\begin{equation}
	\left|\int_{0}^{t}\mathcal{R}_c^1+\mathcal{R}_c^2\right|\ls\|\mathcal{C}(\bv)\|_{L^2_t(L^2)}\|\mathcal{Q}\|_{L^2_t(L^2)}+\|\mathcal{C}(q)\|_{L^2_t(L^2)}\|\mathcal{V}\|_{L^2_t(L^2)}\ls TP\left(\sup_{t\in[0,T]}\mathfrak{E}(t)\right).
\end{equation}
Similar to \eqref{Rvep1} and \eqref{Rvep2}, one has
\begin{equation}\label{RVEP}
	\sum_{i=1}^{6}\left|\int_{0}^{t}\mathcal{R}_\vep^i\right|\ls \delta\vep\|\DD\D_t^m\bv\|_{L^2_t(L^2)}^2+TP\left(\sup_{t\in[0,T]}\mathfrak{E}(t)\right)
\end{equation}
Therefore, integrating \eqref{Ful-tim-est} with respect to time, substituting the above estimates at hand into the resultant equation, using the estimate \eqref{kon} and Korn's inequality \eqref{Korn's ineq}, we have
\begin{align*}
	&\|\mathcal{V}(t)\|_0^2+\|\mathcal{Q}(t)\|_0^2+\|\DD\D_t^m\boe(t)\|_0^2+|\D_t^m\D_1\eta_ka_{k2}(t)|_0^2+\vep\int_{0}^{t}\|\DD\D_t^m\bv\|^2_0\nonumber\\
	&\ls M_0+\delta\sup_{t\in[0,T]}\mathfrak{E}(t)+P\left(\sup_{t\in[0,T]}\mathfrak{E}(t)\right).
\end{align*}
Squaring the above inequality and integrating over time once again yield
\begin{align*}
	&\int_{0}^{t}\|\mathcal{V}\|_0^4+\|\mathcal{Q}\|_0^4+\|\DD\D_t^m\boe\|_0^4+|\D_t^m\D_1\eta_ka_{k2}|_0^4+\vep^2\int_{0}^{t}\left(\int_{0}^{s}\|\DD\D_t^m\bv\|^2_0\right)^2\nonumber\\
	&\ls M_0+TP\left(\sup_{t\in[0,T]}\mathfrak{E}(t)\right).
\end{align*}
By the definition of $\mathcal{V}$ and $\mathcal{Q}$ and the fundamental theorem of calculus, we obatin
\begin{align*}
	\int_{0}^{t}\|\D_t^m\bv\|_0^4+\|\D_t^mq\|_0^4&\ls \int_{0}^{t}\|\mathcal{V}\|_0^4+\|\mathcal{Q}\|_0^4+\int_{0}^{t}\|\D_t^m\boe\|_0^4(\|A_{jk}\D_kv\|_{L^\infty}^4+\|A_{jk}\D_kq\|_{L^\infty}^4)\\
	&\ls M_0+TP\left(\sup_{t\in[0,T]}\mathfrak{E}(t)\right).
\end{align*}
Therefore, we can complete the proof.
\end{pf}

\subsection{Normal derivative estimates}
To derive normal derivative estimate, we can first obtain  from $\eqref{FEL-final}_2$, \eqref{f-eqn} and \eqref{Geo-iden-1}  that
\begin{align*}
	&-\tilde{\rho}_0J\Delta\eta_i-\gamma (\tilde{\rho}_0J^{-1})^{\gamma}a_{ik}a_{rs}\D_k\D_s\eta_r-\mu\vep a_{kl}a_{kj}\D_j\D_l v_i-(\mu+\lambda)\vep  a_{kl}a_{ij}\D_l\D_j v_k\\
	&=-\gamma (\tilde{\rho}_0J^{-1})^{\gamma-1}a_{ik}+\D_j\tilde{\rho}_0J\D_j\eta_i-\tilde{\rho}_0J\D_tv_i+\mu\vep Ja_{kl}\D_l A_{kj}\D_jv_i\\
	&\quad+\mu\vep Ja_{kl}\D_lA_{kj}\D_jv_i+\lambda\vep Ja_{ij}\D_jA_{kl}\D_lv_k,
\end{align*}
where we have used the following fact $\D_k f=-\tilde{\rho}_0J^{-2}a_{rs}\D_k\D_s\eta_r+\D_k\tilde{\rho}_0J^{-1}.$
As a consequence, we have
\begin{equation}\label{normal-eqn}
	-\mathcal{A}_{ij}\D_2^2\eta_j-\mu\vep a_{k2}a_{k2}\D_2^2v_i-(\mu+\lambda)\vep a_{i2}a_{j2}\D_2^2v_j= \mathcal{F}_i+\mathcal{G}_i,
\end{equation}
where 
\begin{equation*}
	\mathcal{A}_{ij}=\tilde{\rho}_0J\delta_{ij}+\gamma (\tilde{\rho}_0J^{-1})^\gamma a_{i2}a_{j2}
\end{equation*}
\begin{align*}
	\mathcal{F}_i=&\sum\limits_{l\neq 2,\text{ or }j\neq 2}\left(\mu\vep a_{kl}a_{kj}\D_j\D_l v_i+(\mu+\lambda)\vep  a_{kl}a_{ij}\D_l\D_j v_k\right)\nonumber\\
	&+\mu\vep Ja_{kl}\D_l A_{kj}\D_jv_i+\mu\vep Ja_{kl}\D_lA_{kj}\D_jv_i+\lambda\vep Ja_{ij}\D_jA_{kl}\D_lv_k
\end{align*}
and 
\begin{align*}
	\mathcal{G}_i=&\gamma (\tilde{\rho}_0J^{-1})^\gamma\left(a_{i1}a_{r1}\D_1^2\eta_r+a_{i1}a_{r2}\D_{12}^2\eta_r+a_{i2}a_{r1}\D_{21}^2\eta_r\right)\nonumber\\
	&-\tilde{\rho}_0J\D_t v_i-\tilde{\rho}_0J\D_1^2\eta_i-\gamma(\tilde{\rho}_0J^{-1})^{\gamma-1}a_{ik}\D_k \tilde{\rho}_0+J\DD\tilde{\rho}_0\cdot\DD\eta_i.		
\end{align*}
It is clear that
\begin{equation*}
	\det \mathcal{A}=\tilde{\rho}_0J+\gamma (\tilde{\rho}_0J^{-1})^\gamma(a_{12}^2+a_{22}^2)>\tilde{\rho}_0J>0,
\end{equation*}
so that $\mathcal{A}$ is invertible. Moreover, it is obvious that $\mathcal{A}$ is also symmetric. These imply that we can estimate the normal derivatives of $\eta$  by using \eqref{normal-eqn} so that we have the following lemma
\begin{lemma}\label{Nor-est}
	For any $t\in[0,T_\vep]$, $m\geq 4$, it holds that
	\begin{align}\label{nor-est}
		&\|\boe^\vep(t)\|_{\MH^m}^2+\|\DD \boe^\vep(t)\|_{\MH^{m-1}}^2+\|\sqrt{\vep}\DD^2\boe^\vep(t)\|_{\MH^{m-1}}^2+\int_{0}^{t}\|\DD\boe^\vep\|_{\MH^m}^2+\|\bv^\vep\|_{\MH^m}^2+\|\vep\DD\bv^\vep\|_{\MH^m}^2\nonumber\\
		&\ls M_0+\delta\sup_{t\in[0,T_\vep]}\mathfrak{E}^\vep(t)+T_\vep^\frac{1}{4}P\left(\sup_{t\in[0,T_\vep]}\mathfrak{E}^\vep(t)\right).
	\end{align}
\end{lemma}
\begin{pf}
	 Applying $\bar{\D}^\beta$ with $|\beta|\leq m-1$ to \eqref{normal-eqn} yields that
	 	\begin{align}\label{2nd-ord-eqn}	    
	 		&-\mathcal{A}_{ij}\bar{\D}^\beta\D_2^2\eta_j-\mu\vep a_{k2}a_{k2}\bar{\D}^\beta\D_2^2v_i-(\mu+\lambda)\vep a_{i2}a_{j2}\bar{\D}^\beta\D_2v_j\nonumber\\
	 		&=[\bar{\D}^\beta,\mathcal{A}_{ij}]\D_2^2\eta_j+\mu\vep[\bar{\D}^\beta,a_{k2}a_{k2}]\D_2^2v_i+(\mu+\lambda)\vep[\bar{\D}^\beta,a_{i2}a_{j2}]\D_2^2v_j+\bar{\D}^\beta\mathcal{F}_i+\bar{\D}^\beta\mathcal{G}_{i}.
	 	\end{align}
	 Then, by taking inner product between the above equation with $-\D_t^\beta\D_2^2 \eta_i$, one has
	 \begin{align*}
	 	&\int_\Omega\mathcal{A}_{ij}\bar{\D}^\beta\D_2^2\eta_j\bar{\D}^\beta\D_2^2\eta_i+\frac{\mu}{2}\frac{d}{dt}\int_\Omega\vep|a_{\cdot2}|^2|\bar{\D}^\beta\D_2^2\boe|^2+\frac{\mu+\lambda}{2}\frac{d}{dt}\int_\Omega\vep|\bar{\D}^\beta\D_2^2\eta_ia_{i2}|^2\nonumber\\
	 	&=\int_\Omega\mu\vep a_{k2}\D_t a_{k2}|\bar{\D}^\beta\D_2^2\boe|^2+(\mu+\lambda)\vep a_{i2}\D_ta_{j2}\bar{\D}^\beta\D_2^2\eta_i\bar{\D}^\beta\D_2^2\eta_j\nonumber\\
	 	&\quad-\int_\Omega[\bar{\D}^\beta,\mathcal{A}_{ij}]\D^2\eta_j\bar{\D}^\beta\D_2^2\eta_i-\mu\int_\Omega\vep[\bar{\D}^\beta,a_{k2}a_{k2}]\D_2^2v_i\bar{\D}^\beta\D_2^2\eta_i\\
	 	&\quad-(\mu+\lambda)\int_\Omega\vep[\bar{\D}^\beta,a_{i2}a_{j2}]\D_2^2v_j\bar{\D}^2\D_2^2\eta_i-\int_\Omega\bar{\D}^\beta\mathcal{F}_i\bar{\D}^\beta\D_2^2\eta_i-\int_\Omega\bar{\D}^\beta\mathcal{G}_i\bar{\D}^\beta\D_2^2\eta_i.\nonumber
	 \end{align*}
%
Then, we can obtain from integrating over time and using the fact $\mu>0,\mu+\lambda>0$ that
\begin{align}\label{nor-est-1}
	\|\bar{\D}^\beta\D_2^2\boe\|_{L^2_t(L^2)}^2+\|\sqrt{\vep}\bar{\D}^\beta\D_2^2\boe\|_0^2(t)\ls& \int_{0}^{t}\|\D_t\D_1\boe\|_{L^\infty}\|\sqrt{\vep}\bar{\D}^\beta\D_2^2\boe\|_0^2+\|[\bar{\D}^\beta,\mathcal{A}_{\cdot j}]\D_2^2\eta_j\|_{L^2_t(L^2)}^2\\
	&+\vep^2\|[\bar{\D}^\beta, a_{\cdot 2}a_{\cdot 2}]\D_2^2\bv\|_{L^2_t(L^2)}^2+\|\bar{\D}^\beta\boldsymbol{\mathcal{F}}\|_{L^2_t(L^2)}^2+\|\bar{\D}^\beta\boldsymbol{\mathcal{G}}\|_{L^2_t(L^2)}^2.\nonumber
\end{align}
It follows from Sobolev embedding and H\"{o}lder's inequality that
\begin{align}\label{nor-est-1.1}
	\int_{0}^{t}\|\D_t\D_1\boe\|_{L^\infty}\|\sqrt{\vep}\bar{\D}^\beta\D_2^2\boe\|_0^2\ls \|\sqrt{\vep}\bar{\D}^\beta\D_2^2\boe\|_{L^\infty_t(L^2)}^2T^\frac{1}{2}\|\DD\boe\|_{L^2_t(\MH^3)}\ls T^\frac{1}{2}P\left(\sup_{t\in[0,T]}\mathfrak{E}(t)\right)
\end{align}
It follows from \eqref{pro-est}, \eqref{Js-est} and \eqref{m-J-q} that
\begin{align}\label{A-com-est}
		\|[\bar{\D}^\beta,\mathcal{A}_{\cdot j}]\D_2^2\eta_j\|_{L^2_t(L^2)}^2
		&\ls\|\bar{\D}\mathcal{A}_{\cdot j}\|_{L^\infty_{t,x}}^2\|\D_2\eta_j\|_{L^2_t(\MH^{m-1})}^2+\|\D_2\eta_j\|_{L^\infty_{t,x}}^2\|\bar{\D}\mathcal{A}_{\cdot j}\|_{L^2_{t}(\MH^{m-1})}^2\nonumber\\
		&\ls\|\D_2\eta_j\|_{L^2_t(\MH^{m-1})}^2\|\bar{\D}\mathcal{A}_{\cdot j}\|_{L^2_{t}(\MH^{m-1})}^2\nonumber\\
		&\ls T\|\DD \boe\|_{L^\infty_t(\MH^{m-1})}(\|J\|_{L^2_t(\MH^m)}+\|\D_1\boe\|_{L^2_t(\MH^m)})\\
		&\ls TP\left(\sup_{t\in[0,T]}\mathfrak{E}(t)\right),\nonumber
\end{align}
and
\begin{align}\label{nor-est-1.2}
	\vep^2\|[\bar{\D}^\beta, a_{\cdot 2}a_{\cdot 2}]\D_2^2\bv\|_{L^2_t(L^2)}^2&\ls\vep^2\left(\|\bar{\D}(a_{\cdot 2}a_{\cdot 2})\|_{L^\infty_{t,x}}^2\|\D_2^2\boe\|_{L^2_t(\MH^{m-1})}^2+\|\D_2^2\boe\|_{L^\infty_{t,x}}^2\|\bar{\D}(a_{\cdot 2}a_{\cdot 2})\|_{L^2_t(\MH^{m-1})}^2\right)\nonumber\\
	&\ls \vep\|\bar{\D}\D_1\boe\|_{L^2_t(\MH^{m-1})}^2T\|\sqrt{\vep}\DD^2\boe\|_{L^\infty_t(\MH^{m-1})}^2\\
	&\ls TP\left(\sup_{t\in[0,T]}\mathfrak{E}(t)\right).\nonumber
\end{align}
Similarly, by using \eqref{pro-est},\eqref{A-priori-assum} and \eqref{tan-est-1}, we have
\begin{align}\label{nor-est-1.3}
	\|\bar{\D}^2\boldsymbol{\mathcal{F}}\|_{L^2_t(L^2)}&\ls \vep^2\left(\|\boldsymbol{a}\boldsymbol{a}\|_{L^\infty_{t,x}}^2\|\bar{\D}^\beta\D_1\DD\bv\|_{L^2_t(L^2)}^2+\|[\bar{\D}^\beta,\boldsymbol{a}\boldsymbol{a}]\D_1\DD\bv\|_{L^2_t(L^2)}^2\right)\nonumber\\
	&\quad+\vep^2\left(\|\DD^2\boe\|_{L^\infty_{t,x}}^2\|\DD \bv\|_{L^2_t(\MH^{m-1})}^2+\|\DD\bv\|_{L^\infty_{t,x}}^2\|\DD^2\boe\|_{L^2_t(\MH^{m-1})}^2\right)\nonumber\\
	&\ls \vep^2\|\bar{\D}^\beta\D_1\DD\bv\|_{L^2_t(L^2)}^2+\vep\|\bar{\D}\DD\boe\|_{L^2_t(\MH^{m-1})}^2T\|\sqrt{\vep}\DD^2\boe\|_{L^\infty_t(\MH^{m-1})}^2\\
	&\ls M_0+\delta\sup_{t\in[0,T]}\mathfrak{E}(t)+T^\frac{1}{4}P\left(\sup_{t\in[0,T]}\mathfrak{E}(t)\right).\nonumber
\end{align}
and
\begin{align}\label{nor-est-1.4}
	\|\bar{\D}^\beta\boldsymbol{\mathcal{G}}\|_{L^2_t(L^2)}^2&\ls\|\bar{\D}^\beta\D_1\DD\boe\|_{L^2_t(L^2)}^2+\|\bar{\D}^\beta\D_t v\|_{L^2_t(L^2)}^2+\left\|\left[\bar{\D}^\beta,\left(\tilde{\rho}_0J^{-1}\right)^\gamma \boldsymbol{a}\boldsymbol{a}\right]\D_1\DD\boe\right\|_{L^2_t(L^2)}^2\nonumber\\
	&\quad+\left\|\left[\bar{\D}^\beta,\tilde{\rho}_0J^{-1}\right]\D_t\DD\bv\right\|_{L^2_t(L^2)}^2+\left\|\left[\bar{\D}^\beta,\tilde{\rho}_0J^{-1}\right]\D_1^2\boe\right\|_{L^2_t(L^2)}^2\nonumber\\
	&\quad+\left\|\left(\tilde{\rho}_0J^{-1}\right)^{\gamma-1}\right\|_{L^\infty_{t,x}}^2\|a_{ik}\D_k\tilde{\rho}_0\|_{L^2_t(\MH^{m-1})}^2+\|a_{ik}\D_k\tilde{\rho}_0\|_{L^\infty_{t,x}}^2\left\|\left(\tilde{\rho}_0J^{-1}\right)^{\gamma-1}\right\|_{L^2_t(\MH^{m-1})}^2\nonumber\\
	&\quad+\|J\DD\tilde{\rho}_0\|_{L^\infty_{t,x}}^2\|\DD\boe\|_{L^2_t(\MH^{m-1})}^2+\|\DD\boe\|_{L^\infty_{t,x}}\|J\DD\tilde{\rho}_0\|_{L^2_t(\MH^{m-1})}^2\\
	&\ls T\|\DD\boe\|_{L^\infty_t(\MH^{m,1}_{tan})}^2+T^\frac{1}{2}\| \bv\|_{L^4_t(\MH^m_{tan})}^2+\|J\|_{L^2_t(\MH^{m-1})}+\|\DD\boe\|_{L^2_t(\MH^{m-1})}\nonumber\\
	&\quad+\|\bar{\D}\DD\boe\|_{L^\infty_{t,x}}^2\left(\|\boe\|_{L^2_t(\MH^{m})}^2+\|J\|_{L^2_t(\MH^{m-1})}\right)+\|\DD\tilde{\rho}_0\|_{L^2_t(\MH^{m-1})}\nonumber\\
	&\ls T^\frac{1}{2}P\left(\sup_{t\in[0,T]}\mathfrak{E}(t)\right).\nonumber
\end{align}
%
Therefore, plugging \eqref{nor-est-1.1}-\eqref{nor-est-1.4} into \eqref{nor-est-1} yields
\begin{equation}\label{2nd-nor}
		\|\D_2^2\boe\|_{L^2_t(\MH^{m-1}_{tan})}^2+\|\sqrt{\vep}\D_2^2\boe\|_{L^\infty_t(\MH^{m-1}_{tan})}^2\ls M_0+\delta\sup_{t\in[0,T]}\mathfrak{E}(t)+T^\frac{1}{4}P\left(\sup_{t\in[0,T]}\mathfrak{E}(t)\right).
\end{equation}
On the other hand, by taking inner product between \eqref{2nd-ord-eqn} with $-\vep\bar{\D}^\beta\D_2^2v_i$, integrating over $\Omega\times[0,t]$ and using \eqref{A-com-est}-\eqref{2nd-nor}, we can obtain that
\begin{align}\label{2nd-nor-1}
	\vep^2\|\D_2^2\bv\|_{L^2_t(\MH^{m-1}_{tan})}^2\ls M_0+\delta\sup_{t\in[0,T]}\mathfrak{E}(t)+T^\frac{1}{4}P\left(\sup_{t\in[0,T]}\mathfrak{E}(t)\right).
\end{align}
Next, for any $|\beta|\leq m-1-\ell,\ell\in \mathbb{N}$, applying $\bar{\D}^{\beta}\D_2^\ell$ to \eqref{normal-eqn}, taking inner product with $-\bar{\D}^\beta\D_2^{\ell+2}\eta_i$, and using the same argument as in the proof of \eqref{2nd-nor} and \eqref{2nd-nor-1}, we can successively obtain the following estimate for $\ell=1,2,\cdots, m-1$ 
\begin{align*}
&\|\bar{\D}^\beta\D_2^{\ell+2}\boe\|_{L^2_t(L^2)}^2+\|\sqrt{\vep}\bar{\D}^\beta\D_2^{\ell+2}\boe\|_{L^\infty_t(L^2)}^2+\|\vep\bar{\D}^\beta\D_2^{\ell+2}\bv\|_{L^2_t(L^2)}\nonumber\\
&\ls M_0+\delta\sup_{t\in[0,T]}\mathfrak{E}(t)+T^\frac{1}{4}P\left(\sup_{t\in[0,T]}\mathfrak{E}(t)\right)	
\end{align*}
As a consequence, one has
\begin{align*}
	&\|\D_2^2\boe\|_{L^2_t(\MH^{m-1})}^2+\|\sqrt{\vep}\D_2^2\boe\|_{L^\infty_t(\MH^{m-1})}^2+\|\vep\D_2^2\bv\|_{L^2_t(\MH^{m-1})}\nonumber\\
	&\ls M_0+\delta\sup_{t\in[0,T]}\mathfrak{E}(t)+T^\frac{1}{4}P\left(\sup_{t\in[0,T]}\mathfrak{E}(t)\right).
\end{align*}
which, combining with \eqref{tan-est-1} and \eqref{full-time-est}, implies
\begin{align*}
	&\|\sqrt{\vep}\DD^2\boe(t)\|_{\MH^{m-1}}^2+\int_{0}^{t}\|\DD\boe\|_{\MH^m}^2+\|\bv^\vep\|_{\MH^m}^2+\|\vep\DD\bv\|_{\MH^m}^2\\
	&\ls M_0+\delta\sup_{t\in[0,T]}\mathfrak{E}(t)+T^\frac{1}{4}P\left(\sup_{t\in[0,T]}\mathfrak{E}(t)\right).
\end{align*}
Furthermore, by using \eqref{L-infL-2}, we have
\begin{align*}
	\|\boe(t)\|_{\MH^m}^2+\|\DD \boe(t)\|_{\MH^{m-1}}^2&\ls\|\boe_0\|_{\MH^m}^2+\|\DD \boe_0\|_{\MH^{m-1}}^2+T\left(\|\D_t\boe\|_{L^2_t(\MH^{m})}^2+\|\DD\D_t\boe\|_{L^2_t(\MH^{m-1})}^2\right)\\
	&\ls M_0+\delta\sup_{t\in[0,T]}\mathfrak{E}(t)+T^\frac{1}{4}P\left(\sup_{t\in[0,T]}\mathfrak{E}(t)\right).
\end{align*}
Therefore, we complete the proof of lemma.
\end{pf}
\subsection{Proof of Proposition \ref{uniform estimates}} We now collect the estimates derived above and verify the \textit{a priori} assumption \eqref{A-priori-assum}. Indeed, we obtain from lemmas \ref{basic-energy-est}, \ref{Tan-est}, \ref{Full-time-derivative} and \ref{Nor-est} that
\begin{equation*}
	\sup_{t\in[0,T_\vep]}\mathfrak{E}^\vep(t)\leq M_0+\delta\sup_{t\in[0,T_\vep]}\mathfrak{E}^\vep(t)+{T_\vep}^\frac{1}{4}P\left(\sup_{t\in[0,T_\vep]}\mathfrak{E}^\vep(t)\right).
\end{equation*}
As a result, it holds for any $t\in[0,T_\vep]$
\begin{equation*}
	|J^\vep(t)-J_0^\vep|\leq \left|\int_{0}^{t}J_t^\vep\right|\leq {T_\vep}^\frac{1}{2}\|J_t^\vep\|_{L^2_T(L^\infty)}\ls{T_\vep}^\frac{1}{2}\sup_{t\in[0,T_\vep]}\mathfrak{E}^\vep(t).
\end{equation*}
Similarly, we also have
\begin{equation*}
	|\D_j\eta^\vep_i(t)-\D_j\eta^\vep_{0i}|\ls {T_\vep}^\frac{1}{2}\sup_{t\in[0,T_\vep]}\mathfrak{E}^\vep(t).
\end{equation*}
Therefore, by taking $\delta$ sufficiently small, there exist a $T$ independent of $\vep$ such that \eqref{A-priori-assum} is satisfied and
\begin{equation*}
	\sup_{t\in[0,T]}\mathfrak{E}^\vep(t)\leq 2M_0.
\end{equation*}

\section{Proof of Theorems}
\subsection{Proof of Theorem \ref{Theorem 1}}
For any fixed $\vep>0$, we can construct the local classical solutions $(\boe^\vep,q^\vep, \bv^\vep)$ to the free boundary problem of compressible viscoelastic fluid system \eqref{FEL-final} by the approach similar to \cite{xu2013} if the initial data $(\tilde{\rho}^\vep_0,\boldsymbol{\eta}_0^\vep,\bv_0^\vep)$ satisfy \eqref{compat-cond}-\eqref{v-uni-bdd}. We omit the detail here for simplicity.  From the uniform estimates of $(\boe^\vep,q^\vep, \bv^\vep)$ obtained in Proposition \ref{uniform estimates}, we can find a $T_0>0$ independent of $\vep$ such that $(\boe^\vep,q^\vep, \bv^\vep)$ satisfy $\sup_{t\in[0,T]}\mathfrak{E}^\vep(t)\leq C_1,$ which complete the proof of theorem \ref{Theorem 1}.
\subsection{Proof of Theorem \ref{vanishing viscosity}}
It follows from Theorem \ref{Theorem 1} that $\boe^\vep$ is uniformly bounded in $L^\infty(0,T_0;\MH^m)$, $\DD\boe^\vep$ is uniformly bounded in $L^2(0,T_0; \MH^m)$, and $\D_t\boe^\vep$ is uniformly bounded in $L^\infty(0,T_0;\MH^{m-1})\cap L^2(0,T_0;\MH^{m-1})$.  Then,  by using Aubin-Lions compactness theorem (c.f. \cite{simon1986}), we have $\boe^\vep$ is compact in $C([0,T_0];\MH^{m-1})$. In particular, there exist a sequence $\vep_n\rightarrow 0^+$ and a $\boe$ such that $\boe^{\vep_n}\rightarrow \boe$ in $C([0,T_0];\MH^{m-1})$ and $\bv^{\vep_n}\rightarrow\bv$ in $C([0,T_0];\MH^{m-2})$ as $\vep_n\rightarrow 0^+$. These convergence allow us to pass the limit in \eqref{FEL-final} and obtain $(\boe,\bv)$ solving the elastodynamic equations \eqref{FEL-final-2}. Moreover, by applying the lower semi-continuity of norms to bounds in Theorem 1, we can obtain that $(\boe,\bv)$ satisfy the following uniform regularities of  
\begin{align*}
	&\|\boe\|_{L^\infty(0,T_0;\MH^m(\Omega))}+\|\DD\boe\|_{L^\infty(0,T_0;\MH^{m,1}_{tan}(\Omega))}+|\bar{\D}^{m-1}\D_1\boe\cdot\bn|_{L^\infty(0,T_0;L^2(\Gamma))}+\|\DD\boe\|_{L^2(0,T_0;\MH^m(\Omega))}\nonumber\\
	&+\|\D_t^m\bv\|_{L^4(0,T_0;L^2(\Omega))}+\|\D_t^m\DD\boe\|_{L^4(0,T_0;L^2(\Omega))}+|\D_1\D_t^m\boe\cdot\bn|_{L^4(0,T_0;L^2(\Gamma))}\leq C.
\end{align*}
 By standard energy method, it is not difficult to prove the uniqueness of classical solutions $(\boe,\bv)$ to \eqref{FEL-final-2} with these regularities. Therefore, we have that the whole family $(\boe^\vep,\bv^\vep)$ converge to $(\boe,\bv)$.

\	
	

	\

	\bibliographystyle{abbrv}
    
    \bibliography{reference}

%
%
%
%
%
%
%
%

\end{document}